\documentclass[a4paper,12pt]{article}
\usepackage{amsmath}
\usepackage{latexsym}
\numberwithin{equation}{section}

\def\R{{\bf R}}

\def\N{{\bf N}}

\def\d{\displaystyle}
\def\e{{\varepsilon}}
\def\o{\overline}
\def\wt{\widetilde}
\def\wh{\widehat}

\def\p{\partial}
\def\v#1{\mbox{\boldmath $#1$}}

\newtheorem{thm}{Theorem}[section]

\newtheorem{lem}{Lemma}[section]
\newtheorem{prop}{Proposition}[section]
\newtheorem{rem}{Remark}[section]

\title{Almost global solutions of semilinear wave equations with the critical exponent
\\
in high dimensions
\footnote{This work is partially supported by
the Grant-in-Aid for Scientific Research (C) (No. 24540183),
Japan Society for the Promotion of Science.}
\vskip5pt
{\small In memory of Professor Rentaro Agemi}}
\author{
Hiroyuki Takamura
\footnote{Department of Complex and Intelligent Systems,
Faculty of Systems Information Science,
Future University Hakodate,
116-2 Kamedanakano-cho,
Hakodate, Hokkaido 041-8655, Japan.
e-mail: takamura@fun.ac.jp.}
\quad
and
\quad
Kyouhei Wakasa
\footnote{The 2nd year of the doctor course,
Department of Mathematics, Hokkaido University, Sapporo, 060-0810, Japan. 
e-mail: wakasa@math.sci.hokudai.ac.jp.
}
}
\date{
\[
\begin{array}{ll}
\mbox{\footnotesize{\bf Keywords:}}
& \mbox{\footnotesize semilinear wave equation, high dimensions, critical exponent, lifespan}\\
\mbox{\footnotesize{\bf MSC2010:}}
& \mbox{\footnotesize primary 35L71, 35E15, secondary 35A01, 35A09, 35B33, 35B44}\\
\end{array}
\]
}

\begin{document}
\maketitle
\begin{abstract}
\footnotesize{
We are interested in the \lq\lq almost" global-in-time existence of classical solutions
in the general theory for nonlinear wave equations.
All the three such cases are known to be sharp
due to blow-up results in the critical case for model equations.
However, it is known that we have a possibility to get the global-in-time existence
for two of them in low space dimensions
if the nonlinear term is of derivatives of the unknown function
and satisfies so-called null condition, or non-positive condition.
But another one for the quadratic term in four space dimensions is out of the case
as the nonlinear term should include a square of the unknown function itself.
\par
In this paper, we get one more example guaranteeing the sharpness
of the almost global-in-time existence in four space dimensions.
It is also the first example of the blow-up of classical solutions
for non-single and indefinitely signed term in high dimensions.
Such a term arises from the neglect of derivative-loss factors
in Duhamel's formula for positive and single nonlinear term.
This fact may help us to describe a criterion to get the global-in-time existence
in this critical situation.}
\end{abstract}


\section{Introduction}
\par
First we shall outline the general theory on the initial value problem
for fully nonlinear wave equations,
\begin{equation}
\label{GIVP}
\left\{
\begin{array}{l}
u_{tt}-\Delta u=H(u,Du,D_xDu) \quad \mbox{in}\quad\R^n\times[0,\infty),\\
u(x,0)=\e f(x),\ u_t(x,0)=\e g(x),
\end{array}
\right.
\end{equation}
where $u=u(x,t)$ is a scalar unknown function of space-time variables,
\[
\begin{array}{l}
Du=(u_{x_0},u_{x_1},\cdots,u_{x_n}),\ x_0=t,\\
D_xDu=(u_{x_ix_j},\ i,j=0,1,\cdots,n,\ i+j\ge1),
\end{array}
\]
$f,g\in C^\infty_0(\R^n)$ and $\e>0$ is a \lq\lq small" parameter.
We note that it is impossible to construct a general theory
for \lq\lq large" $\e$ due to blow-up results.
For example, see Glassey \cite{G73}, Levine \cite{Le74}, or Sideris \cite{Si84b}.
Let
\[
\wh{\lambda}=(\lambda;\ (\lambda_i),i=0,1,\cdots,n;
\ (\lambda_{ij}),i,j=0,1,\cdots,n,\ i+j\ge1). 
\]
Suppose that the nonlinear term $H=H(\wh{\lambda})$ is a sufficiently smooth function with
\[
H(\wh{\lambda})=O(|\wh{\lambda}|^{1+\alpha})
\]
in a neighborhood of $\wh{\lambda}=0$, where $\alpha\ge1$ is an integer.
Let us define the lifespan $\wt{T}(\e)$ of classical solutions of (\ref{GIVP}) by
\[
\begin{array}{ll}
\wt{T}(\e)=\sup\{t>0\ :\ \exists
&\mbox{a classical solution $u(x,t)$ of (\ref{GIVP})}\\
&\mbox{for arbitrarily fixed data, $(f,g)$.}\}.
\end{array}
\]
When $\wt{T}(\e)=\infty$, the problem (\ref{GIVP}) admits
a global-in-time solution,
while we only have a local-in-time solution on $[0,\wt{T}(\e))$
when $\wt{T}(\e)<\infty$.
For local-in-time solutions, one can measure the long time stability
of a zero solution by orders of $\e$.
Because the uniqueness of the solution of (\ref{GIVP})
may yield that $\lim_{\e\rightarrow+0}\wt{T}(\e)=\infty$.
Such an uniqueness theorem can be found in Appendix of John \cite{J90} for example.
From now on, we omit \lq\lq -in-time" and simply use \lq\lq global" and \lq\lq local".
\par
In Chapter 2 of Li and Chen \cite{LC92},
we have long histories on the estimate for $\wt{T}(\e)$.
The lower bounds of $\wt{T}(\e)$ are summarized in the following table.
Let $a=a(\e)$ satisfy
\begin{equation}
\label{a}
a^2\e^2\log(a+1)=1
\end{equation}
and $c$ stands for a positive constant independent of $\e$.
Then,
due to the fact that it is impossible to obtain an $L^2$ estimate for $u$ itself
by standard energy methods, we have
\begin{center}
\begin{tabular}{|c||c|c|c|}
\hline
$\wt{T}(\e)\ge$  & $\alpha=1$ & $\alpha=2$ & $\alpha\ge3$\\
\hline
\hline
$n=2$ &
$\begin{array}{l}
ca(\e)\\
\quad\mbox{in general case},\\
c\e^{-1}\\
\quad\mbox{if}\ \d\int_{\R^2}g(x)dx=0,\\
c\e^{-2}\\
\quad\mbox{if}\ \partial^2_uH(0)=0
\end{array}
$
&
$\begin{array}{l}
c\e^{-6}\\
\quad\mbox{in general case},\\
c\e^{-18}\\
\quad\mbox{if}\ \partial^3_uH(0)=0,\\
\exp(c\e^{-2})\\
\quad\mbox{if}\ \partial^3_uH(0)=\partial^4_uH(0)=0\
\end{array}
$
& $\infty$ \\
\hline
$n=3$ &
$\begin{array}{l}
c\e^{-2}\\
\quad\mbox{in general case},\\
\exp(c\e^{-1})\\
\quad\mbox{if}\ \partial^2_uH(0)=0
\end{array}
$
& $\infty$ & $\infty$ \\\hline
$n=4$ &
$\begin{array}{l}
\exp(c\e^{-2})\\
\quad\mbox{in general case},\\
\infty\\
\quad\mbox{if}\ \partial^2_uH(0)=0
\end{array}$
& $\infty$ & $\infty$ \\
\hline
$n\ge5$ & $\infty$ & $\infty$ & $\infty$\\
\hline
\end{tabular} 
\end{center}
The result for $n=1$ is that
\begin{equation}
\label{lifespan_lower_onedim}
\wt{T}(\e)\ge
\left\{
\begin{array}{lllll}
\e^{-\alpha/2}\quad&\mbox{in general case},&\\
\d c\e^{-\alpha(1+\alpha)/(2+\alpha)}\quad&\mbox{if}\ \d \int_{\R}g(x)dx=0,&\\
c\e^{-\alpha} \quad &\mbox{if}\ \partial^\beta_uH(0)=0
\ \mbox{for}\ 1+\alpha\le\forall\beta\le2\alpha.&
\end{array}
\right.
\end{equation}
For references on these results, see Li and Chen \cite{LC92}.
We shall skip to refer them here.
But we note that two parts in this table are different
from the one in Li and Chen \cite{LC92}.
One is the general case in $(n,\alpha)=(4,1)$.
In this part, the lower bound of $\wt{T}(\e)$ is
$\exp(c\e^{-1})$ in Li and Chen \cite{LC92}.
But later, it has been improved by Li and Zhou \cite{LZ95}.
Another is the case for $\partial^3_uH(0)=0$ in $(n,\alpha)=(2,2)$. 
This part is due to Katayama \cite{Kata01}.
But it is missing in Li and Chen \cite{LC92}.
Its reason is closely related to the sharpness of results in the general theory.
The sharpness is achieved by the fact that there is no possibility to improve
the lower bound of $\wt{T}(\e)$ in sense of order of $\e$
by blow-up results for special equations and special data.
It is expressed in the upper bound of $\wt{T}(\e)$
with the same order of $\e$ as in the lower bound.
On this matter, Li and Chen \cite{LC92} says that
all these lower bounds are known to be sharp except for $(n,\alpha)=(4,1)$.
But before this article,
Li \cite{Li91} says that $(n,\alpha)=(2,2)$ has also open sharpness
while the case for $\partial^3_uH(0)=0$ is still missing.
Li and Chen \cite{LC92} might have dropped the open sharpness in $(n,\alpha)=(2,2)$
by conjecture that $\partial^4_uH(0)=0$ is a technical condition.
No one disagrees with this observation
because the model case of $H=u^4$ has a global solution
in two space dimensions, $n=2$.
See the next section below.
However, Zhou and Han \cite{ZH12} have obtained
this final sharpness in $(n,\alpha)=(2,2)$ by studying $H=u_t^2u+u^4$.
This puts Katayama's result into the table after 20 years from Li and Chen \cite{LC92}.
We note that Godin \cite{Go93} has showed the sharpness of the case
for $\partial^3_uH(0)=\partial^4_uH(0)=0$ in $(n,\alpha)=(2,2)$ by studying $H=u_t^3$.
This result has reproved by Zhou and Han \cite{ZH11_two}.
\par
We now turn back to another open sharpness of the general case in $(n,\alpha)=(4,1)$.
It has been obtained by our previous work, Takamura and Wakasa \cite{TW11},
by studying model case of $H=u^2$.
This part had been open more than 20 years
in the analysis on the critical case for model equations, $u_{tt}-\Delta u=|u|^p\ (p>1)$.
We mention to whole histories on this equation precisely in the next section. 
In this way, the general theory and its optimality have been completed. 
\par
After the completion of the general theory,
we are interested in \lq\lq almost" global existence,
namely, the case where $\wt{T}(\e)$ has an lower bound of the exponential function
of $\e$ with a negative power.
Such a case appears in $(n,\alpha)=(2,2),(3,1),(4,1)$
in the table of the general theory.
It is remarkable that Klainerman \cite{Kl86} and
Christodoulou \cite{Chri86} have independently found a special structure
on $H=H(Du,D_xDu)$ in $(n,\alpha)=(3,1)$
which guarantees the global existence.
This algebraic condition on nonlinear terms of derivatives
of the unknown function is so-called \lq\lq null condition".
It has been also established independently by Godin \cite{Go93} for $H=H(Du)$ 
and Katayama \cite{Kata95} for $H=H(Du,D_xDu)$ in $(n,\alpha)=(2,2)$.
The null condition had been supposed to be not sufficient
for the global existence in $(n,\alpha)=(2,2)$.
For this direction, Agemi \cite{A} proposed \lq\lq non-positive condition"
in this case for $H=H(Du)$.
Thsi conjecture has been verified
by Hoshiga \cite{Ho08} and Kubo \cite{K07} independently.  
It might be necessary and sufficient condition to the global existence.
On the other hand, the situation in $(n,\alpha)=(4,1)$
is completely different from $(n,\alpha)=(2,2),(3,1)$
because $H$ has to include $u^2$.
\par
In this paper, we get the first attempt to clarify a criterion
on $H$ guaranteeing the global existence
by showing another blow-up example of $H$.
More precisely, we have an almost global existence and its optimality
for an equation of the form
\begin{equation}
\label{integral_term}
\begin{array}{ll}
u_{tt}-\Delta u=u^2
&\d-\frac{1}{\pi^2}\int_{0}^{t}d\tau
\int_{|\xi|\le 1}\frac{(u_tu)(x+(t-\tau)\xi,\tau)}{\sqrt{1-|\xi|^2}}d\xi\\
&\d-\frac{\e^2}{2\pi^2}\int_{|\xi|\le 1}\frac{f(x+t\xi)^2}{\sqrt{1-|\xi|^2}}d\xi
\end{array}
\end{equation}
in $\R^4\times[0,\infty)$.
We note that the third term in the right-hand side of (\ref{integral_term})
can be neglected by simple modification.
One can say that this result is the first example of the blowing-up
of a classical solution to nonlinear wave equation with \lq\lq non-single" 
and \lq\lq indefinitely signed" term in high dimensions. 
(\ref{integral_term}) arises from a neglect
of derivative loss factors in Duhamel's term
for positive and single nonlinear term, $u^2$.
We introduce this observation in the next section
with more general situation about on space dimensions.
Therefore, one can conclude that
derivative loss factors in Duhamel's term due to high dimensions
do not contribute to any order of $\e$ in the estimate of the lifespan.
\par
Finally, we note that, in contrast with (\ref{integral_term}),
another equation of the form
\[
\begin{array}{ll}
u_{tt}-\Delta u=u^2
&\d-\frac{1}{2\pi^2}\int_0^t d\tau\int_{|\omega|=1}
(u_tu)(x+(t-\tau)\omega,\tau)dS_\omega\\
&\d-\frac{\e}{4\pi^2}\int_{|\omega|=1}(\e f^2+\Delta f+2\omega\cdot\nabla g)
(x+t\omega)dS_\omega
\end{array}
\]
admits a global classical solution in $\R^4\times[0,\infty)$.
Its details will appear in our forthcoming paper.
\par
This paper is organized as follows.
Main theorems and whole histories on closely related model equations,
$u_{tt}-\Delta u=|u|^p$, are stated
in the next section including how to derive our problem.
They are discussed in all high space dimensions and for the nonlinear term
with fractional powers to describe what the critical power depends on. 
In section 3, we introduce a weighted $L^\infty$ space
in which the solution will be constructed by contraction mapping.
In section 4, we show {\it a priori} estimate for the existence part. 
The lower bounds of the lifespan in odd space dimensions or even space dimensions
are obtained in section 5 or section 6 respectively.
Upper bounds of the lifespan in odd space dimensions are
obtained in section 7 for the critical case and in section 8 for the subcritical case.
Similarly to them, upper bounds of the lifespan in even space dimensions are
obtained in section 9 for the critical case and in section 10 for the subcritical case.
\par
The essential part of this work has been completed when the second author
was in the 2nd year of the master course,
Graduate School of Systems Information Science,
Future University Hakodate.


\section{Model equations and main theorems}
\par
Before deriving our problem (\ref{integral_term}),
we shall introduce the whole histories on the following model problems.
By such problems, the optimality of the general theory
is guaranteed especially in the case where nonlinear term $H$ includes
the lower order of $u$ itself.
This may help us to know the difficulty to analyze
the quadratic terms in four space dimensions.
\par
We first consider an initial value problem,
\begin{equation}
\label{IVP}
\left\{
\begin{array}{l}
u_{tt}-\Delta u=|u|^p,\quad \mbox{in}\quad \R^n\times[0,\infty),\\
u(x,0)=\e f(x),\ u_t(x,0)=\e g(x)
\end{array}
\right.
\end{equation}
assuming that $\e>0$ is \lq\lq small" again.
Let us define a lifespan $T(\e)$ of a solution of (\ref{IVP}) by
\[
T(\e)=\sup\{t>0\ :\ \exists\mbox{a solution $u(x,t)$ of (\ref{IVP})
for arbitrarily fixed $(f,g)$.}\},
\]
where \lq\lq solution" means classical one when $p\ge2$.
When $1<p<2$, it means weak one, but sometimes the one of 
a solution of associated integral equations to (\ref{IVP})
by standard Strichartz's estimate.
See Georgiev, Takamura and Zhou \cite{GTZ06} for example on such an argument. 
\par
When $n=1$, we have $T(\e)<\infty$ for any power $p>1$ by Kato \cite{Kato80}.
When $n\ge2$, we have the following Strauss' conjecture on (\ref{IVP})
by Strauss \cite{St81}. 
\[
\begin{array}{lll}
T(\e)=\infty & \mbox{if $p>p_0(n)$ and $\e$ is \lq\lq small"}
& \mbox{(global-in-time existence)},\\
T(\e)<\infty & \mbox{if $1<p\le p_0(n)$}
& \mbox{(blow-up in finite time)},
\end{array}
\]
where $p_0(n)$ is so-called Strauss' exponent defined
by positive root of the quadratic equation,
\begin{equation}
\label{gamma}
\gamma(p,n)=2+(n+1)p-(n-1)p^2=0.
\end{equation}
That is,
\begin{equation}
\label{p_0(n)}
p_0(n)=\frac{n+1+\sqrt{n^2+10n-7}}{2(n-1)}. 
\end{equation}
We note that $p_0(n)$ is monotonously decreasing in $n$ and $p_0(4)=2$.
This conjecture had been verified by many authors of partial results.
All the references on the final result in each part
can be summarized in the following table.
\begin{center}
\begin{tabular}{|c||c|c|c|}
\hline
& $p<p_0(n)$ & $p=p_0(n)$ & $p>p_0(n)$ \\
\hline
\hline
$n=2$ & Glassey \cite{G81a} & Schaeffer \cite{Sc85} & Glassey \cite{G81b}\\
\hline
$n=3$ & John \cite{J79} & Schaeffer \cite{Sc85} & John \cite{J79}\\
\hline
$n\ge4$ & Sideris \cite{Si84} &  
$
\begin{array}{l}
\mbox{Yordanov $\&$ Zhang \cite{YZ06}}\\
\mbox{Zhou \cite{Z07}, indep.}
\end{array}
$
&
$
\begin{array}{l}
\mbox{Georgiev $\&$ Lindblad}\\
\mbox{$\&$ Sogge \cite{GLS97}}\\
\end{array}
$
\\
\hline
\end{tabular} 
\end{center}
\par
In the blow-up case of $1<p\le p_0(n)$,
we are interested in the estimate of the lifespan $T(\e)$.
From now on, $c$ and $C$ stand for positive constants but independent of $\e$.
When $n=1$, we have the following estimate of the lifespan $T(\e)$ for any $p>1$.
\begin{equation}
\label{lifespan_est_one_dim.}
\left\{
\begin{array}{cl}
c\e^{-(p-1)/2}\le T(\e)\le C\e^{-(p-1)/2}
& \mbox{if}\quad\d\int_{\R}g(x)dx\neq0,\\
c\e^{-p(p-1)/(p+1)}\le T(\e)\le C\e^{-(p-1)/(p+1)}
&\mbox{if} \quad\d\int_{\R}g(x)dx=0.
\end{array}
\right.
\end{equation}
This result has been obtained by Y.Zhou \cite{Z92_one}.
We note that its order of $\e$ coincides with the general theory
when $p=1+\alpha\ (\alpha=1,2,3,\cdots)$.
Moreover, Lindblad \cite{L90} has obtained more precise result for $p=2$,
\begin{equation}
\label{lifespan_est_one_dim._det}
\left\{
\begin{array}{ll}
\d \exists \lim_{\e\rightarrow+0}\e^{1/2}T(\e)>0
&\mbox{for}\quad\d\int_{\R}g(x)dx\neq0,\\
\d \exists \lim_{\e\rightarrow+0}\e^{2/3}T(\e)>0
&\mbox{for}\quad\d\int_{\R}g(x)dx=0.
\end{array}
\right.
\end{equation}
Similarly to this, Lindblad \cite{L90} has also obtained the following result
for $(n,p)=(2,2)$.
\begin{equation}
\label{lifespan_est_two_quad}
\left\{
\begin{array}{ll}
\d \exists \lim_{\e\rightarrow+0}a(\e)^{-1}T(\e)>0
&\mbox{for}\quad\d\int_{\R^2}g(x)dx\neq0\\
\d \exists \lim_{\e\rightarrow+0}\e T(\e)>0
&\mbox{for}\quad\d\int_{\R^2}g(x)dx=0,
\end{array}
\right.
\end{equation}
where $a(\e)$ is the one in (\ref{a}).
\par
When $1<p<p_0(n)\ (n\ge3)$ or $2<p<p_0(2)\ (n=2)$,
we have the following conjecture.
\begin{equation}
\label{lifespan1}
c\e^{-2p(p-1)/\gamma(p,n)}\le T(\e)\le C\e^{-2p(p-1)/\gamma(p,n)},
\end{equation}
where $\gamma(p,n)$ is defined by (\ref{gamma}).
We note that (\ref{lifespan1}) coincides
with the second line in (\ref{lifespan_est_one_dim.})
if we define $\gamma(p,n)$ by (\ref{gamma}) even for $n=1$.
All the results verifying this conjecture are summarized in the following table.
\begin{center}
\begin{tabular}{|c||c|c|c|}
\hline
& lower bound of $T(\e)$ & upper bound of $T(\e)$\\
\hline
\hline
$n=2$ & Zhou \cite{Z93} & Zhou \cite{Z93}\\
\hline
$n=3$ & Lindblad \cite{L90} & Lindblad \cite{L90} \\
\hline  
$n\ge4$ 
&
$
\begin{array}{l}
\mbox{Lindblad $\&$ Sogge \cite{LS96}}\\
\mbox{: $n\le 8$ or radially symmetric sol.}\\
\end{array}
$
&
Sideris \cite{Si84}\\
\hline
\end{tabular} 
\end{center}
We note that, for $n=2,3$,
\[
 \exists \lim_{\e\rightarrow+0}\e^{2p(p-1)/\gamma(p,n)}T(\e)>0.
\]
is established.
Moreover, the upper bound in the the case where $n\ge 4$
easily follows from the rescaling method applied to the proof in Sideris \cite{Si84}
which proves $T(\e)<\infty$.
Such an argument can be found in Georgiev, Takamura and Zhou \cite{GTZ06}.
\par
On the other hand, when $p=p_0(n)$,
we have the following conjecture.
\begin{equation}
\label{lifespan2}
\exp\left(c\e^{-p(p-1)}\right)\le T(\e)\le\exp\left(C\e^{-p(p-1)}\right).
\end{equation}
All the results verifying this conjecture are also summarized in the following table.
\begin{center}
\begin{tabular}{|c||c|c|c|}
\hline
& lower bound of $T(\e)$ & upper bound of $T(\e)$\\
\hline
\hline
$n=2$ & Zhou \cite{Z93} & Zhou \cite{Z93}\\
\hline
$n=3$ & Zhou \cite{Z92_three} & Zhou \cite{Z92_three} \\
\hline  
$n\ge4$ 
&
$
\begin{array}{l}
\mbox{Lindblad $\&$ Sogge \cite{LS96}}\\
\mbox{: $n\le 8$ or radially symm. sol.}
\end{array}
$
&
Takamura $\&$ Wakasa \cite{TW11}\\
\hline
\end{tabular} 
\end{center}
In this way, we realize that one of the last open problem had been
the upper bound of $T(\e)$ for the critical case $p=p_0(n)$
in high dimensions $n\ge4$.
This difficulty is due to the lack of the positivity of the fundamental solution
of linear wave equations which is caused by so-called \lq\lq derivative loss".
\par
We are now in a position to derive our problem (\ref{integral_term})
in the general situation by neglecting such derivative loss factors.
From now on, we assume that $n\ge2$ and
write $n=2m,\ 2m+1\ (m=1,2,3,\cdots)$.
Let us consider the following initial value problem.
\begin{equation}
\label{u_depend}
\left\{
\begin{array}{l}
u_{tt}-\Delta u=F(u) \quad \mbox{in}\quad\R^n\times[0,\infty),\\
u(x,0)=\e f(x),\ u_t(x,0)=\e g(x),
\end{array}
\right.
\end{equation}
where $f\in C^{m+3}(\R^n),\ g\in C^{m+2}(\R^n)$ and $F\in C^{m+1}(\R)$.
Then, any solution $u$ of (\ref{u_depend}) has to satisfy
\begin{equation}
\label{IE_u_depend}
u(x,t)=\e u^0(x,t)+\int_0^tR(F(u(\cdot,\tau))|x,t-\tau)d\tau.
\end{equation}
Here we set
\[
u^0(x,t)=\p_t R(f|x,t)+R(g|x,t)
\]
and
\[
R(\phi|x,t)=
\frac{1}{(2m-1)!!}\left(\frac{1}{t}\frac{\p}{\p t}\right)^{m-1}
\left\{t^{2m-1}M(\phi|x,t)\right\},
\]
where
\begin{equation}
\label{M}
M(\phi|x,r)=
\left\{
\begin{array}{ll}
\d\frac{1}{\omega_{n}}\int_{|\omega|=1}\phi(x+r\omega)dS_{\omega}
&\mbox{for}\ n=2m+1,\\
\d\frac{2}{\omega_{n+1}}\int_{|\xi|\le 1}\frac{\phi(x+r\xi)}{\sqrt{1-|\xi|^2}}d\xi
&\mbox{for}\ n=2m.
\end{array}
\right.
\end{equation}
We note that $\omega_n$ stands for a measure of the unit sphere in $\R^n$, i.e. 
\[
\omega_{n}=\frac{2\pi^{n/2}}{\Gamma\left(n/2\right)}=
\left\{
\begin{array}{ll}
\d\frac{2(2\pi)^m}{(2m-1)!!} & \mbox{for}\ n=2m+1,\\
\d\frac{2\pi^m}{(m-1)!} & \mbox{for}\ n=2m.
\end{array}
\right.
\]
This representation formula is well-known.
See pp.681-692 in Courant and Hilbert \cite{CH62} for example.
\par
If we neglect derivative loss factors, namely all the derivatives of $F(u)$,
from the Duhamel's term in (\ref{IE_u_depend}),
then we get a new integral equation,
\begin{equation}
\label{IE0}
u(x,t)=\e u^0(x,t)+L(F(u))(x,t)\quad\mbox{for}\ (x,t)\in\R^n\times[0,\infty),
\end{equation}
where we set
\begin{equation}
\label{L}
L(F(u))(x,t)=\frac{1}{2m-1}\int_{0}^{t}(t-\tau)M(F(u(\cdot,\tau))|x,t-\tau)d\tau.
\end{equation}
A simple computation yields that if $u$ is a $C^2$ solution of (\ref{IE0})
with $f\in C^{m+3}(\R^n),\ g\in C^{m+2}(\R^n)$ and $F\in C^2(\R)$,
then $u$ satisfies
\begin{equation}
\label{NWIVP0}
\left\{
\begin{array}{ll}
\begin{array}{ll}
u_{tt}-\Delta u=
& F(u)-G(x,t)\\
&\d-\frac{2(m-1)}{2m-1}M(F(\e f)|x,t)
\end{array}
& \mbox{in}\ \R^n\times[0,\infty),\\
u(x,0)=\e f(x),\ u_t(x,0)=\e g(x),
& x\in \R^n
\end{array}
\right.
\end{equation}
in the classical sense, where we set
\begin{equation}
\label{H}
G(x,t)=\frac{2(m-1)}{2m-1}\int_0^tM(F'(u(\cdot,\tau))u_t(\cdot,\tau)|x,t-\tau)d\tau.
\end{equation}
In this way, (\ref{integral_term}) follows from setting
$n=4\ (m=2)$ and $F(u)=u^2$ in (\ref{NWIVP0}).
\begin{rem}
The uniqueness of the solution of (\ref{NWIVP0}) is open
while Agemi, Kubota and Takamura \cite{AKT94} has a wrong comment on this fact
after (1.8) on 242p. in \cite{AKT94}.
The restricted uniqueness theorem such as in Appendix 1 in John \cite{J90}
cannot be applicable here because (99a) in \cite{J90} does not hold for this case.
\end{rem}
\par
Moreover, assuming $F\in C^{m+1}(\R)$ if $f\not\equiv0$,
and replacing $\e u^0(x,t)$ by a classical solution $v=v(x,t)$ of
\begin{equation}
\label{v}
\left\{
\begin{array}{ll}
v_{tt}-\Delta v=\d\frac{2(m-1)}{2m-1}M(F(\e f)|x,t)
& \mbox{in}\ \R^n\times[0,\infty),\\
v(x,0)=\e f(x),\ v_t(x,0)=\e g(x)
& x\in \R^n
\end{array}
\right.
\end{equation}
in (\ref{IE0}), we have that a $C^2$ solution $u$ of
\begin{equation}
\label{IE}
u(x,t)=v(x,t)+L(F(u))(x,t)\quad\mbox{for}\ (x,t)\in\R^n\times[0,\infty)
\end{equation}
satisfies
\begin{equation}
\label{NWIVP}
\left\{
\begin{array}{ll}
u_{tt}-\Delta u=F(u)-G(x,t)
& \mbox{in}\ \R^n\times[0,\infty),\\
u(x,0)=\e f(x),\ u_t(x,0)=\e g(x),
& x\in \R^n
\end{array}
\right.
\end{equation}
in the classical sense.
We note that $G\equiv0$ when $n=2,3\ (m=1)$.
This observation appears in pp.254-255 of Agemi, Kubota and Takamura \cite{AKT94}.
\par
From now on, in order to clarify the critical situation
in (\ref{NWIVP0}) and (\ref{NWIVP}),
we may set $F(u)=|u|^p$, or $|u|^{p-1}u$ $(p>1)$.
But, as we see, the existence of a classical solution $v$ of
(\ref{v}) with $f\not\equiv0$ is guaranteed by additional regularity $F\in C^{m+1}$.
In such case we have to regularize $F(u)$ at $u=0$.
Let us define a lifespan $\o{T}(\e)$ by
\[
\begin{array}{ll}
\o{T}(\e)=\sup\{t>0\ :\ \exists & \mbox{a solution $u(x,t)$ of (\ref{IE0}) or (\ref{IE})}\\
&\mbox{for arbitrarily fixed data, $(f,g)$.}\},
\end{array}
\]
where \lq\lq solution" means a $C^2$ solution for $p\ge2$,
or the $C^1$ solution for $1<p<2$.
Agemi, Kubota and Takamura \cite{AKT94} have obtained that
$\o{T}(\e)=\infty$ for $p>p_0(n)$ if $\e$ is small enough,
where $p_0(n)$ is Strauss' exponent.
Their theorem is written for (\ref{IE}) only,
but it is trivial to be available also for (\ref{IE0}).
\par
Our purpose in this paper is to establish the same results for $\o{T}(\e)$
as in (\ref{lifespan1}) and (\ref{lifespan2}) when $n\ge4$ and $1<p\le p_0(n)$. 
They are divided into two theorems below.
We note that one can expect to get a $C^2$ solution
only for $n=4$ and $p=p_0(4)=2$ in this situation.
Except for this case, we assume on $F$ that
\begin{equation}
\left\{
\begin{array}{l}
\label{hypo_F}
\mbox{there exists a constant $A>0$ such that}\\ 
\mbox{$F\in C^1(\R)$ satisfies that $F(0)=F'(0)=0$ and}\\
\mbox{$|F'(s)-F'(\wt{s})|\le pA|s-\wt{s}|^{p-1}$ for $1<p<2$.}
\end{array}
\right.
\end{equation}
Note that (\ref{hypo_F}) implies that $|F(s)|\le A|s|^p$ for $1<p<2$.
We also assume on the data that
\begin{equation}
\left\{
\begin{array}{l}
\label{hypo_data}
\mbox{both $f\in C_0^{m+3}(\R^n)$ and $g\in C^{m+2}_{0}(\R^n)$ do not}\\
\mbox{vanish identically and have compact support}\\
\mbox{contained in $\{x\in\R^n\ :\ |x|\le k\}$ with some constant $k>0$.}
\end{array}
\right.
\end{equation}
Then, we have the following existence theorem for large time interval. 
\begin{thm}
\label{thm:main1}
Let $n\ge4$ and $1<p\le p_0(n)$.
Assume $F(s)=As^2$ when $n=4$ and $p=p_0(4)=2$ or {\rm (\ref{hypo_F})} otherwise, 
where $A$ is a positive constant. 
Suppose that {\rm (\ref{hypo_data})} is fulfilled.
Moreover, assume that $F\in C^{m+1}(\R)$ for {\rm (\ref{IE})}.
Then there exists a positive constant $\e_0=\e_0(f,g,n,p,k)$ such that
the lifespan $\o{T}(\e)$ satisfies
\begin{equation}
\label{lifespan_lower_main}
\begin{array}{ll}
\d \o{T}(\e)\ge c\e^{-2p(p-1)/\gamma(p,n)} 
\quad\mbox{if}\ 1<p<p_0(n),\\
\d \o{T}(\e)\ge \exp\left(c\e^{-p(p-1)}\right)
\quad\mbox{if}\ p=p_0(n)
\end{array}
\end{equation}
for any $\e$ with $0<\e \le\e_0$, where $c$ is a positive constant independent of $\e$.
\end{thm}
This theorem follows from the similar weighted $L^\infty$ iteration method
to Agemi, Kubota and Takamura \cite{AKT94}.
Its basic argument has been introduced by John \cite{J79}
in the simplest case for $n=3$.
The proof is divided into sections 3, 4, 5 ,6 and 7 below.
\par
For the counter part,
the following assumptions on the data are required.
\begin{equation}
\label{blowup_asm}
\left\{
\begin{array}{l}
\mbox{Let $f\equiv 0$, $g(x)=g(|x|)$ and $g\in C^{1}_{0}([0,\infty))$ satisfy that there}\\
\mbox{exist positive constants $k_0$ and $k_1$ with $0<k_0<k_1<k$}\\
\mbox{such that the following three conditions hold.}\\
\mbox{(i)\quad supp $g\subset\{x\in\R^n\ :\ |x|\le k.\}$}\\
\mbox{(ii)\quad $g(|x|)\ge 0$ for $k_0<|x|<k$ and
$\d \int_{(k_1+k)/2}^{k}\hspace{-20pt}\lambda^{[n/2]}g(\lambda)d\lambda>0$,}\\
\mbox{(iii)\quad $k_0$ is sufficiently close to $k$ to satisfy}\\
\mbox{\qquad $\d P_{m}(z)>\frac{1}{2}$ and $\d T_{m}(z)>\frac{1}{2}$
for all $\d z>\frac{k_0}{k}$,}\\ 
\mbox{\qquad where $P_m$ or $T_m$ denote Legendre or Tschebyscheff}\\
\mbox{\qquad polynomials of degree $m$ respectively.}
\end{array}
\right.
\end{equation}
Then, we have the following blow-up theorem. 
\begin{thm}
\label{thm:main2}
Let $n\ge4$ and $F(u)=|u|^p$ with $1<p\le p_0(n)$.
Assume that {\rm (\ref{blowup_asm})}.
Then there exist a positive constant $\e_1=\e_1(g,n,p,k)$ such that
the lifespan $\o{T}(\e)$ satisfies
\begin{equation}
\label{lifespan_main}
\begin{array}{ll}
\d \o{T}(\e)\le C\e^{-2p(p-1)/\gamma(p,n)}
\quad\mbox{if}\ 1<p<p_0(n),\\
\d \o{T}(\e)\le\exp\left(C\e^{-p(p-1)}\right)
\quad\mbox{if}\ p=p_0(n)
\end{array}
\end{equation}
for any $\e$ with $0<\e \le\e_1$, where $C$ is a positive constant independent of $\e$.
\end{thm}
The proof of this theorem is an iteration argument of point-wise estimates
which is basically introduced by John \cite{J79} for $n=3$.
But for the critical case,
we have to reduce the proof to the argument of Zhou \cite{Z92_three}
which compares the solution with a blowing-up solution
of nonlinear ordinary differential equation of the second order.
We also have to employ the slicing method of the blow-up set
which is introduced by Agemi, Kurokawa and Takamura \cite{AKT00}
due to technical difficulties in high dimensions. 
See sections 7, 8, 9 and 10 below.

\section{Weighted $L^\infty$ space}
First, we shall state the following two lemmas on $v$
which play key roles in proofs of 
Theorem \ref{thm:main1} and Theorem \ref{thm:main2}.
The first one is Huygens' principle in odd space dimensions.
\begin{lem}[Agemi, Kubota and Takamura \cite{AKT94}]
\label{lem:huygens}
Let $n=5,7,9,\cdots$.
Under the same assumption as in Theorem \ref{thm:main1},
there exists a classical solution of {\rm (\ref{v})} which satisfies
\begin{equation}
\label{huygens}
\mbox{\rm supp}\ v\subset\{x\in\R^n\ :\  t-k\le|x|\le t+k\}.
\end{equation}
\end{lem}
See 253p. in \cite{AKT94} for the proof of this lemma.
\par
Next, we shall introduce the decay estimate for $v$. First, we 
write $v$ in the form 
\begin{equation}
\label{v_sum}                                                     
v=v_0+v_1. 
\end{equation}
Here, $v_0=\e u^0$ which is a linear part of (\ref{IE_u_depend}) and 
$v_1$ is a solution to the inhomogeneous wave equation
\begin{equation}
\label{v_1}
\left\{
\begin{array}{ll}
\d(v_1)_{tt}-\Delta v_1=\frac{2(m-1)}{2m-1}M(F(\e f)|x,t)
& \mbox{in}\ \R^n\times[0,\infty),\\
v_1(x,0)=(v_1)_t(x,0)=0,
& x\in \R^n
\end{array}
\right.
\end{equation}
where $M$ is defined in (\ref{M}). 
Then we have the following lemma.
\begin{lem}[Agemi, Kubota and Takamura \cite{AKT94}]
\label{lem:decay_est_v}
Under the same assumption as in Theorem \ref{thm:main1},
there exists a positive constant $C_{n,k,f,g}$
depending only on $n$, $k$, $f$ and $g$ such that
$v_0$ and $v_1$ satisfies
\begin{equation}
\label{decay_est_v_0odd}
\sum_{|\alpha|\le1}|\nabla_x^{\alpha}v_0(x,t)|
\le\frac{C_{n,k,f,g}\e}{(t+|x|+2k)^{(n-1)/2}}
\end{equation}
and
\begin{equation}
\label{decay_est_v_1odd}
\sum_{|\alpha|\le1}|\nabla_x^{\alpha}v_1(x,t)|
\le\frac{C_{n,k,f^p}\e^p}{(t+|x|+2k)^{(n-1)/2}}
\end{equation}
when $n=5,7,9,\cdots$, or
\begin{equation}
\label{decay_est_v_0even}
\sum_{|\alpha|\le2}|\nabla_x^{\alpha}v_0(x,t)| 
\le\frac{C_{n,k,f,g}\e}{(t+|x|+2k)^{(n-1)/2}(t-|x|+2k)^{(n-1)/2}}
\end{equation}
and
\begin{equation}
\label{decay_est_v_1even}
\sum_{|\alpha|\le2}|\nabla_x^{\alpha}v_1(x,t)| 
\le\frac{C_{n,k,f^p}\e^p}{(t+|x|+2k)^{(n-1)/2}(t-|x|+2k)^{(n-3)/2}}
\end{equation}
when $n=4,6,8,\cdots$, where
$C_{n,k,f^p}$ is a non-negative constant depending on $n,k,f^p$
with $C_{n,k,f^p}=0$ if and only if $f\equiv0$.
\end{lem}
This lemma directly follows from Lemma 3.2 and Lemma 3.4
in \cite{AKT94}.
We omit the proof here.
\begin{rem}
\label{rem:F(f)}
It is trivial that
Lemma \ref{lem:huygens} and Lemma \ref{lem:decay_est_v}
are available with $v=\e u^0$ if $f\equiv0$.
Therefore we have to prove Theorem \ref{thm:main1} only for {\rm (\ref{IE})}
because all the estimates for {\rm (\ref{IE0})} follow from setting $f\equiv0$
which is a special case of {\rm (\ref{IE})}.
\end{rem}
\par
Taking into account of these lemmas,
we shall introduce {\it a priori} estimate in the weighted $L^\infty$ space as follows.
It is obvious that $v$ is a global classical solution of (\ref{v}).
Therefore our unknown function shall be $U=u-v$.
Then, (\ref{IE}) can be rewritten into
\begin{equation}
\label{IE'}
U=L(F(U+v)).
\end{equation}
In what follows, we shall construct a solution of (\ref{IE'})
in a weighted $L^\infty$ space.
For this purpose, define the sequence of functions $\{U_l\}_{l\in\N}$ by
\begin{equation}
\label{sequence1}
U_l=L(F(U_{l-1}+U_0)),\quad U_0=v\quad \mbox{and}\quad U_{00}=v_0,\quad U_{01}=v_1,
\end{equation}
where $L$ is the one in (\ref{L}).
We denote a weighted $L^\infty$ norm of $U$ by 
\begin{equation}
\label{norm}
\|U\|=\sup_{(x,t)\in\R^n\times[0,T]}\{w(|x|,t)|U(x,t)|\}
\end{equation}
with a weighted function
\begin{equation}
\label{weight}
w(r,t)=
\left\{
\begin{array}{lll}
\d\tau_+(r,t)^{(n-1)/2}\tau_-(r,t)^q
&\mbox{if}&\d p>\frac{n+1}{n-1},\\
\d\tau_+(r,t)^{(n-1)/2}\left(\log4\frac{\tau_+(r,t)}{\tau_-(r,t)}\right)^{-1}
&\mbox{if}&\d p=\frac{n+1}{n-1},\\
\d\tau_+(r,t)^{(n-1)/2+q}
&\mbox{if}&\d 1<p<\frac{n+1}{n-1},
\end{array}
\right.
\end{equation}
where we set
\[
\tau_+(r,t)=\frac{t+r+2k}{k}, \quad\tau_-(r,t)=\frac{t-r+2k}{k}
\]
and $q$ is defined by
\begin{equation}
\label{q}
q=\frac{n-1}{2}p-\frac{n+1}{2}.
\end{equation}
In order to get a $C^1$ solution of (\ref{IE'}),
we shall show the convergence of $\{U_l\}_{l\in\N}$
in a function space $X$ defined by 
\[
\begin{array}{ll}
X=\d \left\{U\in C^1(\R^n\times[0,T])\ :\right.
&\d\|U\|_{X}<\infty,\\
&\d\left.\mbox{supp}\ U(x,t)\subset\{|x|\le t+k\}\right\}
\end{array}
\]
which equips a norm
\[
\|U\|_{X}=\sum_{|\alpha|\le 1}\|\nabla_x^\alpha U\|.
\]
In view of (\ref{IE}), we note that $\partial U/\partial t$ can be expressed
in terms of $\nabla_x U$.
Hence we consider spatial derivatives of $U$ only.
Moreover, we see that $X$ is a Banach space for any fixed $T>0$.
Because it follows from the definition of the norm (\ref{norm})
that there exists a positive constant $C_T$ depending on $T$ such that 
\[
\|U\|\ge C_T|U(x,t)|,\quad t\in[0,T].
\]
\par
Later, we shall make use of H\"older's inequality
\begin{equation}
\label{Holder}
\||U_1|^a|U_2|^b\|\le\|U_1\|^a\|U_2\|^b,
\quad a+b=1,\quad a,b\in [0,1]. 
\end{equation}
Further more,
we denote $\partial/\partial x_i$ by $\partial_i$ for $i=1,2,\cdots,n$, 
and set
\[
\begin{array}{l}
\d\p^jW_l=\max_{|\alpha|\le j}\{|\nabla_x^{\alpha}U_l|,|\nabla_x^{\alpha}U_{l-1}|\},\\
\d\p^jW_0=\max_{|\alpha|\le j}\{|\nabla_x^{\alpha}U_0|\},
\quad\p^jW_{0a}=\max_{|\alpha|\le j}\{|\nabla_x^{\alpha}U_{0a}|\}
\end{array}
\]
for $j,a=0,1$. One may omit $\partial^0$ and denote $\partial=\partial^1$.
\section{A priori estimate}
In this section, we show {\it a priori} estimate
which plays a key role in the contraction mapping argument. 
The following lemma is one of the most essential estimates.
\begin{lem}
\label{lm:apriori}
Let $L$ be the linear integral operator defined by {\rm (\ref{L})}.
Assume that $U\in C^0(\R^n\times[0,T])$ with
{\rm  supp} $U\subset\{(x,t)\in \R^n\times[0,T] : |x|\le t+k\}$ and $\|U\|<\infty$.
Then, there exists a positive constant C independent of $k$ and $T$ such that 
\begin{equation}
\label{apriori}
\|L(|U|^p)\|\le Ck^2\|U\|^pD(T),
\end{equation}
where $D(T)$ is defined by 
\begin{equation}
\label{D_1}
D(T)=
\left\{
\begin{array}{ll}
\d1 & \mbox{if}\ p>p_0(n),\\
\d\log\frac{2T+3k}{k} & \mbox{if}\ p=p_0(n),\\
\d\left(\frac{2T+3k}{k}\right)^{\gamma(p,n)/2} & 
\mbox{if}\ 1<p<p_0(n)
\end{array}
\right.
\end{equation}
and $\gamma(p,n)$ is the one in {\rm (\ref{gamma})}.
\end{lem}
\par
We note that Lemma \ref{lm:apriori} is not sufficient to cover
all the cases on the exponent. 
In fact, we need the following variants of {\it a priori} estimate up to space dimensions.
\begin{lem}
\label{lm:apriori_odd_u0}
Let $n=5,7,9,\cdots$ and $L$ be the linear integral operator defined by {\rm (\ref{L})}.
Assume that $U,U_0 \in C^0(\R^n\times[0,T])$
with {\rm  supp} $U\subset\{(x,t)\in \R^n\times[0,T]:|x|\le t+k\}$, 
{\rm  supp} $U_0\subset\{(x,t)\in \R^n\times[0,T]:t-k\le |x|\le t+k\}$ 
and $\|U\|,\|\tau_{+}^{(n-1)/2}U_{0}w^{-1}\|<\infty$. 
Then, there exists a positive constant $C_{n,\nu,p}$
depending on $n$, $\nu$ and $p$ such that 
\begin{equation}
\label{apriori_odd_u0}
\|L(|U_0|^{p-\nu}|U|^\nu)\|
\le C_{n,\nu,p}k^2\left\|\frac{\tau_{+}^{(n-1)/2}}{w}U_0\right\|^{p-\nu}
\|U\|^{\nu}E_{\nu}(T),
\end{equation}
where $0\le \nu\le p$. $E_{\nu}(T)$ is defined by 
\begin{equation}
\label{E(nu)_1}
E_{\nu}(T)=
\left\{
\begin{array}{lll}
\d 1 & 
\mbox{if}\ \d p>\frac{n+1}{n-1},\\
\d  \left(\frac{2T+3k}{k}\right)^{\nu\delta} &
\mbox{if}\ \d p=\frac{n+1}{n-1},\\
\d \left(\frac{2T+3k}{k}\right)^{-\nu q} &
\mbox{if}\ \d p<\frac{n+1}{n-1}
\end{array}
\right.
\mbox{for}\ 0\le\nu<p,
\end{equation}
where $q$ is the one in (\ref{q}) and $\delta$ is a small positive constant. 
When $\nu=p$, (\ref{apriori_odd_u0}) coincides with (\ref{apriori})
as $E_p(T)=D(T)$ and $C_{n,p,p}=C$. 
\end{lem}

\begin{lem}
\label{lm:apriori_even_2}
Let $n=4,6,8,\cdots$ and $L$ be the linear integral operator defined by (\ref{L}). 
Assume that $U,U_{00}, U_{01}\in C^0(\R^n\times[0,T])$ 
with {\rm  supp}$(U,U_{00},U_{01})$$\subset\{(x,t)\in \R^n\times[0,T] : |x|\le t+k\}$
and $\|U\|,\|(\tau_+\tau_-)^{(n-1)/2}U_{0a}(w\tau_-^a)^{-1}\|<\infty\ (a=0,1)$.
Then, there exists a positive constant $C_{n,\nu,p}$ depending on 
$n$, $\nu$, and $p$ such that 
\begin{equation}
\begin{array}{ll}
\label{apriori_even_2}
\|L(|U_{0a}|^{p-\nu}|U|^{\nu})\|
\d \le C_{n,\nu,p}k^2\left\|(\tau_{+}\tau_{-})^{(n-1)/2}
\frac{U_{0a}}{w \tau_{-}^a}\right\|^{p-\nu}\|U\|^{\nu}E_{\nu,a}(T),
\end{array}
\end{equation}
where $0\le \nu\le p$ and $a=0,1$. When $0\le \nu<p$, 
$E_{\nu,a}(T)$ is defined by
\begin{equation}
\label{E(nu,a)_1}
E_{\nu,a}(T)=
\left\{
\begin{array}{lll}
\d 1 & 
\mbox{if}\ \mu<-1,\\
\d \log\frac{2T+3k}{k} &
\mbox{if}\ \mu=-1,\\
\d \left(\frac{2T+3k}{k}\right)^{1+\mu} &
\mbox{if}\ \mu>-1
\end{array}
\right.
\mbox{for}\ p>\frac{n+1}{n-1},
\end{equation}
where $\d\mu=(p-\nu)\left(a-\frac{n-1}{2}\right)-\nu q$
and $q$ is the one in (\ref{q}), and
\begin{equation}
\label{E(nu,t)_2}
E_{\nu,a}(T)=
\left\{
\begin{array}{lll}
\d \log\frac{2T+3k}{k} &
\mbox{if}\ \sigma=-1,\ \nu=0,\\
\d \left(\frac{2T+3k}{k}\right)^{1+\sigma} &
\mbox{if}\ \sigma>-1,\\
\d \left(\frac{2T+3k}{k}\right)^{\nu\delta} &
\mbox{otherwise}
\end{array}
\right.
\mbox{for}\ p=\frac{n+1}{n-1},
\end{equation}
where $\d\sigma=(p-\nu)\left(a-\frac{n-1}{2}\right)$ and 
$\delta$ stands for any positive constant, and
\begin{equation}
\label{E(nu,t)_3}
E_{\nu,a}(T)=
\left\{
\begin{array}{lll}
\d \left(\frac{2T+3k}{k}\right)^{-\nu q} & 
\mbox{if}\ \sigma<-1,\\
\d \log\frac{2T+3k}{k}\left(\frac{2T+3k}{k}\right)^{-\nu q} &
\mbox{if}\ \sigma=-1,\\
\d \left(\frac{2T+3k}{k}\right)^{1+\mu} &
\mbox{if}\ \sigma>-1
\end{array}
\right.
\mbox{for}\ p<\frac{n+1}{n-1}.
\end{equation}
When $\nu=p$, (\ref{apriori_even_2}) coincides with (\ref{apriori})
as $E_{p,a}(T)=D(T)$ for $a=0,1$ and $C_{n,p,p}=C$. 
\end{lem}
\begin{lem}
\label{lem:apriori_even_3}
Suppose that the same assumption as in Lemma \ref{lm:apriori_even_2} is fulfilled. 
Then, there exists a positive constant $C_{n,\nu,p}$ depending on 
$n$, $\nu$, and $p$ such that
\begin{equation}
\begin{array}{ll}
\label{apriori_even_3}
\|L(|U_{00}|^{p-\nu}|U_{01}|^{\nu})\|\\
\quad\d \le C_{n,\nu,p}k^2\left\|(\tau_{+}\tau_{-})^{(n-1)/2}
\frac{U_{00}}{w}\right\|^{p-\nu}
\left\|(\tau_{+}\tau_{-})^{(n-1)/2}\frac{U_{01}}{w\tau_{-}}\right\|^{\nu}F_{\nu}(T),
\end{array}
\end{equation}
where $0\le \nu \le p$. When $0<\nu<p$, $F_{\nu}(T)$ is defined by 
\begin{equation}
\label{F(nu,t)_4}
F_{\nu}(T)=
\left\{
\begin{array}{lll}
\d 1 & 
\mbox{if}\ \kappa<-1,\\
\d \log\frac{2T+3k}{k} &
\mbox{if}\ \kappa=-1,\\
\d \left(\frac{2T+3k}{k}\right)^{1+\kappa} &
\mbox{if}\ \kappa>-1,
\end{array}
\right.
\end{equation}
where $\d\kappa=\nu-\frac{n-1}{2}p$. 
When $\nu=0$ or $\nu=p$, (\ref{apriori_even_3}) coincides with (\ref{apriori_even_2})
as $F_{0}(T)=E_{0,0}(T)$ for $a=\nu=0$ and $F_{p}(T)=E_{0,1}(T)$ for $a=1$, $\nu=p$. 
\end{lem}
\par
Four lemmas above follows from the following basic estimate. 
\begin{lem}{\rm\bf (Basic estimate)}
\label{lem:basic_est}
Let $L$ be the linear integral operator defined by {\rm (\ref{L})} and 
$a_1\ge0$, $a_2 \in\R$ and $a_3\ge0$. Then, there exists a positive constant 
$C_{n,p,a_1,a_2,a_3}$ such that 
\begin{equation}
\label{basic_est}
\begin{array}{lll}
L\left\{\tau_{+}^{-(n-1)p/2+a_1}\tau_{-}^{a_2}
\left(\log(4\tau_{+}/\tau_{-})\right)^{a_3}\right\}(x,t)\\
\d \le C_{n,p,a_1,a_2,a_3}k^2w(r,t)^{-1}
\left(\frac{2T+3k}{k}\right)^{a_1}E_{a_1,a_2,a_3}(T)
\end{array}
\end{equation}
for $|x|\le t+k$, $t\in[0,T]$, where $E_{a_1,a_2,a_3}(T)$ 
is defined by 
\begin{equation}
\label{E_gen}
E_{a_1,a_2,a_3}(T)=
\left\{
\begin{array}{lll}
\d 1 & 
\mbox{if}\ \d a_2<-1\ \mbox{and}\ a_3=0,\\
\d \log\frac{2T+3k}{k} &
\mbox{if}\ \d a_2=-1\ \mbox{and}\ a_3=0,\\
\d \left(\frac{2T+3k}{k}\right)^{\delta a_3} &
\mbox{if}\ a_2\le-1\ \mbox{and}\ a_3>0,\\
\d \left(\frac{2T+3k}{k}\right)^{1+a_2} &
\mbox{if}\ \d a_2>-1,
\end{array}
\right.
\end{equation}
where $\delta$ stands for any positive constant.
\end{lem}
\par\noindent
{\bf Proof.}
First we employ the following fundamental identity for spherical means.
\begin{lem}[John \cite{J55}]
\label{lm:Planewave}
Let $b\in C([0,\infty))$.
Then, the identity    
\begin{equation}
\label{Planewave}
\begin{array}{ll}
\d \int_{|\omega|=1}b(|x+\rho \omega|)dS_\omega
\d = 2^{3-n}\omega_{n-1}(r\rho)^{2-n}\int_{|\rho-r|}^{\rho+r}\lambda h(\lambda, \rho ,r)
b(\lambda)d\lambda
\end{array} 
\end{equation}
holds for $x\in\R^n,\ r=|x|$ and $\rho>0$, where 
\begin{equation}
\label{h}
h(\lambda, \rho ,r)
=\{\lambda^2-(\rho-r)^2\}^{(n-3)/2}\{(\rho+r)^2-\lambda^2\}^{(n-3)/2}.
\end{equation}
\end{lem}
See \cite{J55} for the proof of this lemma.
\par
In order to continue the proof of Lemma \ref{lem:basic_est},
we need radially symmetric versions of $L$ which follows from Lemma \ref{lm:Planewave}. 
From now on, a positive constant $C$ depending only on $n$ and $p$
may change from line to line.
\begin{lem}
\label{lem:int_new_L}
Let $L$ be the linear integral operator defined by {\rm (\ref{L})} 
and $\Psi=\Psi(|x|,t)\in C([0,\infty)^2),\ x\in\R^n$.
Then,
\begin{equation}
\label{int_w_odd}
L\left(\Psi\right)(x,t)=L_{odd}\left(\Psi\right)(r,t),\ r=|x|
\end{equation}
holds for $n=5,7,9\cdots$ and
\begin{equation}
\label{int_w_even}
L\left(\Psi\right)(x,t)=L_{even,1}\left(\Psi\right)(r,t)+L_{even,2}\left(\Psi\right)(r,t),\ r=|x|
\end{equation}
holds for $n=4,6,8\cdots$, where $L_{odd}\left(\Psi\right)$ is defined by
\begin{equation}
\label{I_odd}
\begin{array}{ll}
L_{odd}\left(\Psi\right)(r,t)\\
\d =Cr^{2-n}\int_{0}^{t}(t-\tau)^{3-n}d\tau
\int_{|t-\tau-r|}^{t-\tau+r}\lambda h(\lambda,t-\tau,r)\Psi(\lambda,\tau)d\lambda
\end{array}
\end{equation}
and each $L_{even,i}\left(\Psi\right)\ (i=1,2)$ is defined by 
\begin{equation}
\label{even_int_domain1}
\begin{array}{ll}
L_{even,1}\left(\Psi\right)(r,t)=
&\d Cr^{2-n}\int_{0}^{t}(t-\tau)^{2-n}d\tau\int_{|t-r-\tau|}^{t+r-\tau}
\lambda \Psi(\lambda,\tau)d\lambda\times\\
&\d\times \int_{|\lambda-r|}^{t-\tau}
\frac{\rho h(\lambda,\rho,r)}{\sqrt{(t-\tau)^2-\rho^2}}d\rho,
\end{array}
\end{equation}
\begin{equation}
\label{even_int_domain2}
\begin{array}{ll}
L_{even,2}\left(\Psi\right)(r,t)=
&\d Cr^{2-n}\int_{0}^{(t-r)_{+}}\!\!\!(t-\tau)^{2-n}d\tau
\int_{0}^{t-r-\tau}\!\!\!\lambda \Psi(\lambda,\tau)d\lambda\times\\
&\d\quad \times \int_{|\lambda-r|}^{\lambda+r}
\frac{\rho h(\lambda,\rho,r)}{\sqrt{(t-\tau)^2-\rho^2}}d\rho.
\end{array}
\end{equation}
Here the usual notation $a_+=\max\{a,0\}$ is used. 
\end{lem}
{\bf Proof.}
(\ref{int_w_odd}) immediately follows from Lemma \ref{lm:Planewave}. 
For (\ref{int_w_even}), 
we make use of changing variables by $y-x=(t-\tau)\xi$ in (\ref{L}).
Then, we obtain
\[
L\left(\Psi\right)(x,t)=C\int_{0}^{t}(t-\tau)^{2-n}d\tau\int_{|y-x|\le t-\tau}
\frac{\Psi(|y|,\tau)}{\sqrt{(t-\tau)^2-|y-x|^2}}dy.
\]
Introducing polar coordinates, we have 
\[
\begin{array}{lll}
\d L\left(\Psi\right)(x,t)&=\d C\int_{0}^{t}(t-\tau)^{2-n}d\tau\int_{0}^{t-\tau}
\frac{\rho^{n-1}d\rho}{\sqrt{(t-\tau)^2-\rho^2}}&\times\\
&\quad \times \d \int_{|\omega|=1}
\Psi(|x+\rho\omega|,\tau)dS_{\omega}.&
\end{array}
\]
Thus Lemma \ref{lm:Planewave} yields
\begin{equation}
\label{int_new_L_last_even}
\begin{array}{llll}
L\left(\Psi\right)(x,t)& \d = Cr^{2-n}\int_{0}^{t}(t-\tau)^{2-n}d\tau\int_{0}^{t-\tau}
\frac{\rho d\rho}{\sqrt{(t-\tau)^2-\rho^2}}\times&\\
&\quad \times \d \int_{|\rho-r|}^{\rho+r}\lambda 
\Psi(\lambda,\tau)h(\lambda,\rho,r)d\lambda.
\end{array}
\end{equation}
Therefore, (\ref{int_w_even}) follows from inverting the order of $(\rho,\lambda)$-integral
in (\ref{int_new_L_last_even}).
\hfill$\Box$
\par
In order to estimate the kernel $h(\lambda,\rho,r)$,
we need the following lemma.
\begin{lem}[Agemi, Kubota and Takamura \cite{AKT94}]
\label{lm:h}
Let $h(\lambda, \rho, r)$ be the one in (\ref{h}).
Suppose that $|\rho-r|\le \lambda \le \rho+r$ and $\rho\ge 0$.
Then, the inequality
\begin{equation}
\label{diff_est_int_1}
|\lambda-r|\le \rho \le \lambda+r
\end{equation}
holds. Moreover, the following three estimates are available.
\begin{eqnarray}
\label{h_1}
h(\lambda,\rho,r)&\le& C r^{n-3}\lambda^{n-3},\\
\label{h_2}
h(\lambda,\rho,r)&\le& C \rho^{n-3}r^{(n-3)/2}\lambda^{(n-3)/2},\\
\label{h_3}
h(\lambda,\rho,r)&\le& C r^{n-3}\rho^{n-3}.
\end{eqnarray}
\end{lem}
See pp.257-258 in \cite{AKT94} for the proof of this lemma.
\par
Let us continue the proof of Lemma \ref{lem:basic_est}.
For simplicity, we set 
\[
\begin{array}{ll}
I_{odd}(r,t)
&\d=L_{odd}\left\{\tau_{+}^{-(n-1)p/2+a_1}\tau_{-}^{a_2}
\left(\log(4\tau_{+}/\tau_{-})\right)^{a_3}\right\}(r,t),\\
I_{even,i}(r,t)
&\d=L_{even,i}\left\{\tau_{+}^{-(n-1)p/2+a_1}\tau_{-}^{a_2}
\left(\log(4\tau_{+}/\tau_{-})\right)^{a_3}\right\}(r,t)\ (i=1,2).
\end{array}
\]
\par\noindent
{\bf Estimates for $\v{I_{odd}}$ and $\v{I_{even,1}}$.}
We shall estimate $I_{odd}$ and $I_{even,1}$ on the following three domains.
\[
\begin{array}{l}
D_1=\{(r,t)\ |\ r\ge t-r>-k\ \mbox{and}\ r\ge 2k\},\\
D_2=\{(r,t)\ |\ r\ge t-r>-k\ \mbox{and}\ r\le 2k\},\\
D_3=\{(r,t)\ |\ t-r\ge r\}.
\end{array}
\]
\par\noindent
(i) Estimate in $D_1$,
\par
Making use of (\ref{h_2}), we get 
\begin{equation}
\label{I_odd_est_1}
\begin{array}{lll}
I_{odd}(r,t)&\d\le \d Cr^{-(n-1)/2}\int_{0}^{t}d\tau
\int_{|t-\tau-r|}^{t+r-\tau}\lambda^{(n-1)/2}\tau_{-}(\lambda,\tau)^{a_2}\times&\\
&\quad\d \times \tau_{+}(\lambda,\tau)^{-(n-1)p/2+a_1}
\left(\log{4\frac{\tau_{+}(\lambda,\tau)}{\tau_{-}(\lambda,\tau)}}\right)^{a_3}d\lambda&
\end{array}
\end{equation}
and
\begin{equation}
\label{I_odd_est_2}
\begin{array}{lll}
\d I_{even,1}(r,t)&\d\le \d Cr^{-(n-1)/2}\int_{0}^{t}(t-\tau)^{2-n}d\tau
\int_{|t-\tau-r|}^{t+r-\tau}\lambda^{(n-1)/2}\times&\\
&\quad\d\times\tau_{+}(\lambda,\tau)^{-(n-1)p/2+a_1}
\tau_{-}(\lambda,\tau)^{a_2}\times&\\
&\quad\d\times\left(\log{4\frac{\tau_{+}(\lambda,\tau)}
{\tau_{-}(\lambda,\tau)}}\right)^{a_3}d\lambda
\d\int_{|\lambda-r|}^{t-\tau}\frac{\rho^{n-2}}{\sqrt{(t-\tau)^2-\rho^2}}d\rho.&
\end{array}
\end{equation}
If one apply the simple inequality
\begin{equation}
\label{h_even_est_1}
\int_{|\lambda-r|}^{t-\tau}\frac{\rho^{n-2}}{\sqrt{(t-\tau)^2-\rho^2}}
d\rho\le (t-\tau)^{n-2} \quad \mbox{for}\ 0\le\tau\le t
\end{equation}
to the right-hand side of (\ref{I_odd_est_2}),
the same quantity as the right-hand side of (\ref{I_odd_est_1}) appears.
Hence, we shall estimate for $I_{odd}$ only from now on. 
\par
Changing variables in (\ref{I_odd_est_1}) by
\begin{equation}
\label{alpha_beta}
\alpha=\tau+\lambda,\ \beta=\tau-\lambda,
\end{equation}
we get
\[
\begin{array}{lll}
I_{odd}(r,t)&\le\d Cr^{-(n-1)/2}\int_{-k}^{t-r}
\left(\frac{\beta+2k}{k}\right)^{a_2}d\beta
\int_{|t-r|}^{t+r}(\alpha-\beta)^{(n-1)/2}\times &\\
&\d \quad \times\d\left(\frac{\alpha+2k}{k}\right)^{-(n-1)p/2+a_1}
\left(\log{4\frac{\alpha+2k}{\beta+2k}}\right)^{a_3}d\alpha.
\end{array}
\]
It follows from
\[
\frac{r}{k}=\frac{2r+r+r}{4k}\ge \frac{\tau_{+}(r,t)}{4}
\]
that
\begin{equation}
\label{I_odd_est_3}
\begin{array}{lll}
&\d I_{odd}(r,t)&\\
&\d \le\d C\tau_{+}(r,t)^{-(n-1)/2}\left(\frac{t+r+2k}{k}\right)^{a_1}
\int_{-k}^{t-r}\left(\frac{\beta+2k}{k}\right)^{a_2}d\beta\times&\\
&\quad \times\d\int_{t-r}^{t+r}
\left(\frac{\alpha+2k}{k}\right)^{-1-q}
\left(\log{4\frac{\alpha+2k}{\beta+2k}}\right)^{a_3}d\alpha.&
\end{array}
\end{equation}
When $a_3=0$, $\alpha$-integral in (\ref{I_odd_est_3}) is dominated by 
\[
\left\{
\begin{array}{lll}
\d Ck\tau_{-}^{-q} & 
\mbox{if}\ \d p>\frac{n+1}{n-1},\\
\d k\log\frac{\tau_{+}}{\tau_{-}} &
\mbox{if}\ \d p=\frac{n+1}{n-1},\\
\d Ck\tau_{+}^{-q} &
\mbox{if}\ \d p<\frac{n+1}{n-1}
\end{array}
\right.
\]
and $\beta$-integral oin (\ref{I_odd_est_3}) is dominated by 
\[
\left\{
\begin{array}{lll}
\d \frac{k}{-(1+a_2)} & 
\mbox{if}\ \d a_2<-1,\\
\d k\log\frac{t-r+2k}{k} &
\mbox{if}\ \d a_2=-1,\\
\d \frac{k}{(1+a_2)}\left(\frac{t-r+2k}{k}\right)^{1+a_2} &
\mbox{if}\ \d a_2>-1.
\end{array}
\right.
\]
(\ref{basic_est}) is now established for $a_3=0$.
\par
When $a_3>0$, we employ the following simple lemma. 
\begin{lem}
\label{lm:log}
Let $\delta>0$ be any given constant. Then, we have 
\begin{equation}
\label{log}
\log X \le \frac{X^{\delta}}{\delta}\ for\ X\ge 1. 
\end{equation}
\end{lem}
The proof of this lemma follows from elementary computation.
We shall omit it.
Then, it follows from Lemma \ref{lm:log} that 
\[
\begin{array}{lll}
\d I_{odd}(r,t)&\le\d C(4\delta^{-1})^{a_3}\tau_{+}(r,t)^{-(n-1)/2}
\left(\frac{t+r+2k}{k}\right)^{a_1+\delta a_3}\times&\\
&\d\quad\times\int_{-k}^{t-r}\left(\frac{\beta+2k}{k}\right)^{a_2-\delta a_3}d\beta
\d\int_{t-r}^{t+r}\left(\frac{\alpha+2k}{k}\right)^{-1-q}d\alpha.&
\end{array}
\]
The $\alpha$-integral above can be estimated by the same manner in the case of $a_3=0$.
The $\beta$-integral is dominated by 
\begin{equation}
\label{beta_int_est_2}
\left\{
\begin{array}{lll}
\d \frac{-k}{1+a_2-\delta a_3} & 
\mbox{if}\ \d a_2\le-1,\\
\d \frac{k}{1+a_2-\delta a_3}\left(\frac{t-r+2k}{k}\right)^{1+a_2-\delta a_3} &
\mbox{if}\ \d a_2>-1
\end{array}
\right.
\end{equation}
with $\delta>0$ satisfying $1+a_2-\delta a_3>0$.
Therefore $I_{odd}$ and $I_{even,1}$ are bounded in $D_1$
by the quantity in the right-hand side of (\ref{basic_est}) as desired.
\par\noindent
(ii) Estimate in $D_2$.
\par
By making use of (\ref{h_3}), we have
\begin{equation}
\label{I_odd_est_5}
\begin{array}{ll}
I_{odd}(r,t)&
\le \d Cr^{-1}\int_{0}^{t}d\tau
\int_{|t-\tau-r|}^{t+r-\tau}\lambda \tau_{+}(\lambda,\tau)^{-(n-1)p/2+a_1}\times\\
&\quad\d \times\tau_{-}(\lambda,\tau)^{a_2}
\left(\log{4\frac{\tau_{+}(\lambda,\tau)}{\tau_{-}(\lambda,\tau)}}\right)^{a_3}d\lambda
\end{array}
\end{equation}
and
\[
\begin{array}{ll}
\d I_{even,1}(r,t)&
\le \d Cr^{-1}\int_{0}^{t}(t-\tau)^{-1}d\tau
\int_{|t-\tau-r|}^{t+r-\tau}\lambda \tau_{+}(\lambda,\tau)^{-(n-1)p/2+a_1}\times\\
&\quad\d \times\tau_{-}(\lambda,\tau)^{a_2}
\left(\log{4\frac{\tau_{+}(\lambda,\tau)}{\tau_{-}(\lambda,\tau)}}\right)^{a_3}d\lambda
\int_{|\lambda-r|}^{t-\tau}\frac{\rho}{\sqrt{(t-\tau)^2-\rho^2}}d\rho.
\end{array}
\]
Similarly to the estimate in $D_1$, the simple inequality
\[
\int_{|\lambda-r|}^{t-\tau}\frac{\rho}{\sqrt{(t-\tau)^2-\rho^2}}d\rho
\le t-\tau \quad \mbox{for}\ 0\le\tau\le t
\]
helps us to estimate $I_{odd}$ only.
Changing variables by (\ref{alpha_beta}), we get 
\begin{equation}
\label{I_odd_est_6}
\begin{array}{lll}
I_{odd}(r,t)&\d \le\d Ckr^{-1}\left(\frac{t+r+2k}{k}\right)^{a_1}
\int_{-k}^{t-r}\left(\frac{\beta+2k}{k}\right)^{a_2}d\beta\times&\\
&\quad \times\d\int_{t-r}^{t+r}
\left(\frac{\alpha+2k}{k}\right)^{1-(n-1)p/2}
\left(\log{4\frac{\alpha+2k}{\beta+2k}}\right)^{a_3}d\alpha.&
\end{array}
\end{equation}
\par
Note that both $\tau_{+}$ and $\tau_{-}$ are numerical constants in this domain,
and that the integrand of both $\alpha$-integral and $\beta$-integral
in (\ref{I_odd_est_6}) 
is numerical constant $C_{a_1,a_2,a_3}$ depending on $a_1$, $a_2$ and $a_3$.
Hence we have
\begin{equation}
\label{I_odd_est_8}
I_{odd}(r,t)\le CC_{a_1,a_2,a_3}kr^{-1}\int_{-k}^{t-r}d\beta\d\int_{t-r}^{t+r}d\alpha
\le\d CC_{a_1,a_2,a_3}k^2.
\end{equation}
This is the desired estimate in $D_2$.
\par\noindent
(iii) Estimate in $D_3$.
\par
By the same reason, we have to estimate $I_{odd}$ in (\ref{I_odd_est_6}) only.
In $D_3$, since $1-(n-1)p/2<0$ is trivial, we get
\[
\left(\frac{\alpha+2k}{k}\right)^{1-(n-1)p/2}
\le\left(\frac{t-r+2k}{k}\right)^{1-(n-1)p/2}\le Cw(r,t)^{-1}
\]
because $t-r\ge r$ is equivalent to $3(t-r)\ge t+r$.
Hence, when $a_3=0$, we obtain 
\[
I_{odd}(r,t)\le\d Ckw(r,t)^{-1}\left(\frac{t+r+2k}{k}\right)^{a_1}
\int_{-k}^{t-r}\left(\frac{\beta+2k}{k}\right)^{a_2}d\beta.
\]
When $a_3>0$, due to Lemma \ref{lm:log}, we have
\[
\begin{array}{ll}
\d I_{odd}(r,t)\le
&\d C(4\delta^{-1})^{a_3}kw(r,t)^{-1}\times\\
&\d\times\left(\frac{t+r+2k}{k}\right)^{a_1+\delta a_3}
\int_{-k}^{t-r}\left(\frac{\beta+2k}{k}\right)^{a_2-\delta a_3}d\beta.
\end{array}
\]
Therefore, in view of (\ref{beta_int_est_2}),
$I_{odd}$ and $I_{even,1}$ are bounded in $D_3$ by the quantity
in the right-hand side of (\ref{basic_est}).
\vskip10pt
\par\noindent
{\bf Estimates for $\v{I_{even,2}}$.}
We shall estimate $I_{even,2}$ on the following three domains.
\[
\begin{array}{l}
D_4=\{(r,t)\ |\ 0<t-r\le k\ \mbox{and}\ t\le 2k\},\\
D_5=\{(r,t)\ |\ 0<t-r\le k\ \mbox{and}\ t\ge 2k\},\\
D_6=\{(r,t)\ |\ t-r\ge k\}.
\end{array}
\]
\par\noindent
(iv) Estimate in $D_4$,
\par
Note that $\tau_{+}$ and $\tau_{-}$ are numerical constants in this case. 
By virtue of (\ref{diff_est_int_1}) and (\ref{h_1}), we get 
\[
\begin{array}{lll}
I_{even,2}(r,t)
&\d\le \d CC_{a_1,a_2,a_3}r^{-1}\int_{0}^{t-r}(t-\tau)^{2-n}d\tau
\int_{0}^{t-r-\tau}\lambda^{n-2}d\lambda\\
&\quad\d\times\int_{|\lambda-r|}^{\lambda+r}\frac{\rho}{\sqrt{(t-\tau)^2-\rho^2}}d\rho.
\end{array}
\]
It follows from 
\begin{equation}
\label{rho-est}
\int_{|\lambda-r|}^{\lambda+r}\frac{\rho d\rho}{\sqrt{(t-\tau)^2-\rho^2}}
\le \frac{2r\lambda}{\sqrt{t-\tau+\lambda+r}\sqrt{t-\tau-\lambda-r}}
\end{equation}
that 
\[
\begin{array}{lll}
&\d I_{even,2}(r,t)&\\
&\le \d CC_{a_1,a_2,a_3}\int_{0}^{t-r}(t-\tau)^{3/2-n}d\tau
\int_{0}^{t-r-\tau}\frac{\lambda^{n-1}}{\sqrt{t-\tau-\lambda-r}}d\lambda.&
\end{array}
\]
Noticing that
\[
\lambda\le t-r-\tau\le t-\tau \quad
\mbox{for}\ \tau\ge 0,
\]
we obtain
\[
\begin{array}{llll}
\d I_{even,2}(r,t)&=
\d CC_{a_1,a_2,a_3}\int_{0}^{t-r}(t-\tau)^{1/2}d\tau
\int_{0}^{t-r-\tau}\frac{d\lambda}{\sqrt{t-\tau-\lambda-r}}.&
\end{array}
\]
Making use of (\ref{alpha_beta}), we have
\[
\begin{array}{lll}
I_{even,2}(r,t)&\d \le CC_{a_1,a_2,a_3}k^{1/2}\int_{-k}^{t-r}d\beta
\int_{\beta}^{t-r}\frac{d\alpha}{\sqrt{t-r-\alpha}}&\\
&\le CC_{a_1,a_2,a_3}k^2.&
\end{array}
\]
This is the desired estimate in $D_4$. 
\par\noindent
(v) Estimate in $D_5$.
\par
In this domain, (\ref{basic_est}) follows from
\[
I_{even,2}(r,t)\le Ck^2C_{a_1,1_2,a_3}\tau_{+}(r,t)^{-(n-1)/2}.
\] 
To see this, we shall employ (\ref{diff_est_int_1}) and (\ref{h_2}).
Then we have 
\[
\begin{array}{lll}
\d I_{even,2}(r,t)&\d\le \d Cr^{-(n-1)/2}\int_{0}^{t-r}(t-\tau)^{2-n}d\tau
\int_{0}^{t-r-\tau}\lambda^{(n-1)/2}\times&\\
&\quad\d\times\tau_{+}(\lambda,\tau)^{-(n-1)p/2+a_1}
\tau_{-}(\lambda,\tau)^{a_2}\times&\\
&\quad\d\times\left(\log{4\frac{\tau_{+}(\lambda,\tau)}
{\tau_{-}(\lambda,\tau)}}\right)^{a_3}d\lambda
\d\int_{|\lambda-r|}^{\lambda+r}\frac{\rho^{n-2}}{\sqrt{(t-\tau)^2-\rho^2}}d\rho.&
\end{array}
\]
Similarly to (\ref{I_odd_est_8}),
it follows from (\ref{h_even_est_1}) that
\[
I_{even,2}(r,t)\le CC_{a_1,a_2,a_3}r^{-(n-1)/2}k^2.
\]
In $D_5$, we have $r\ge k$ which implies $r/k\ge C\tau_{+}(r,t)$.
Hence, we obtain the desired estimate. 
\par\noindent
(vi) Estimate in $D_6$.
\par
By virtue of (\ref{h_1}) and (\ref{rho-est}), we get  
\[
\begin{array}{lll}
&\d I_{even,2}(r,t)&\\
&\le \d C\int_{0}^{t-r}(t-\tau)^{3/2-n}d\tau
\int_{0}^{t-r-\tau}\frac{\lambda^{n-1}}{\sqrt{t-\tau-\lambda-r}}\times&\\
&\d\quad \times \tau_{+}(\lambda,\tau)^{-(n-1)p/2+a_1}
\tau_{-}(\lambda,\tau)^{a_2}\left(\log{4\frac{\tau_{+}(\lambda,\tau)}
{\tau_{-}(\lambda,\tau)}}\right)^{a_3}d\lambda.&\\
\end{array}
\]
Then, we divide the integral of the right-hand side above as
\[
I_{even,3}(r,t)+I_{even,4}(r,t),
\]
where 
\[
\begin{array}{llll}
\d I_{even,3}(r,t)&=
\d C\int_{0}^{(t-r)/2}(t-\tau)^{3/2-n}d\tau
\int_{0}^{t-r-\tau}\frac{\lambda^{n-1}}{\sqrt{t-\tau-\lambda-r}}\times&\\
&\d\quad \times \tau_{+}(\lambda,\tau)^{-(n-1)p/2+a_1}
\tau_{-}(\lambda,\tau)^{a_2}\left(\log{4\frac{\tau_{+}(\lambda,\tau)}{\tau_{-}(\lambda,\tau)}}\right)^{a_3}d\lambda&
\end{array}
\]
and
\[
\begin{array}{llll}
\d I_{even,4}(r,t)&=
\d C\int_{(t-r)/2}^{t-r}(t-\tau)^{3/2-n}d\tau
\int_{0}^{t-r-\tau}\frac{\lambda^{n-1}}{\sqrt{t-\tau-\lambda-r}}\times&\\
&\d\quad \times \tau_{+}(\lambda,\tau)^{-(n-1)p/2+a_1}
\tau_{-}(\lambda,\tau)^{a_2}\left(\log{4\frac{\tau_{+}(\lambda,\tau)}
{\tau_{-}(\lambda,\tau)}}\right)^{a_3}d\lambda.&
\end{array}
\]
\par
First, we shall estimate $I_{even,3}$. 
It follows from (\ref{alpha_beta}) that 
\[
\begin{array}{l}
I_{even,3}(r,t)\\
\d \le Ck^{1/2}\left(\frac{t+r+2k}{k}\right)^{3/2-n}
\left(\frac{t-r+2k}{k}\right)^{n-1-(n-1)p/2+a_1}\times\\
\d\quad\times\int_{-k}^{t-r}\left(\frac{\beta+2k}{k}\right)^{a_2}d\beta
\int_{\beta}^{t-r}\left(\log{4\frac{\alpha+2k}{\beta+2k}}\right)^{a_3}
\frac{d\alpha}{\sqrt{t-r-\alpha}}
\end{array}
\]
because of $n-1-(n-1)p/2\ge0$ for $p\le2$.
When $a_3=0$, we get 
\[
\begin{array}{ll}
\d I_{even,3}(r,t)\le
&\d Ck\left(\frac{t+r+2k}{k}\right)^{3/2-n}\times\\
&\d \times\left(\frac{t-r+2k}{k}\right)^{-q+(n-2)/2+a_1}
\int_{-k}^{t-r}\left(\frac{\beta+2k}{k}\right)^{a_2}d\beta.\\
\end{array}
\]
Hence the desired estimate follows from
\begin{equation}
\label{case1_est4_even_1}
\frac{\tau_{-}(r,t)^{-q+(n-2)/2}}{\tau_{+}(r,t)^{n-3/2}}\le Cw(r,t)^{-1}
\end{equation}
in this case.
When $a_3>0$, Lemma \ref{lm:log} yields that 
\[
\begin{array}{ll}
\d I_{even,3}(r,t)\le
&\d Ck(4\delta^{-1})^{a_3}\left(\frac{t+r+2k}{k}\right)^{3/2-n}\times\\
&\d \left(\frac{t-r+2k}{k}\right)^{-q+(n-2)/2+a_1+\delta a_3}
\int_{-k}^{t-r}\left(\frac{\beta+2k}{k}\right)^{a_2-\delta a_3}d\beta.
\end{array}
\]
Hence the desired estimate follows
from (\ref{beta_int_est_2}) and (\ref{case1_est4_even_1}). 
\par
Next, we shall estimate $I_{even,4}$.
If $r\ge t-r\ge k$, (\ref{alpha_beta}) yields that
\[
\begin{array}{ll}
\d I_{even,4}(r,t)\le
&\d \frac{Ck^{n-1}}{r^{n-3/2}}\left(\frac{t-r+2k}{k}\right)^{n-1-(n-1)p/2+a_1}\times\\
&\d \times\int_{-k}^{t-r}\left(\frac{\beta+2k}{k}\right)^{a_2}d\beta
\int_{\beta}^{t-r}\!\!\!\frac{d\alpha}{\sqrt{t-r-\alpha}}
\left(\log{4\frac{\alpha+2k}{\beta+2k}}\right)^{a_3}.
\end{array}
\]
In this case, we have $r/k\ge C\tau_{+}(r,t)$,
so that (\ref{basic_est}) follows from the same argument as for $I_{even,3}$.
On the other hand, if $t-r\ge r$ and $t-r\ge k$, we have
\[
\tau+\lambda+2k\ge \frac{t-r}{2}+2k\ge \frac{t+r+2k}{6}\ 
\mbox{for}\ \tau\ge \frac{t-r}{2},\ \lambda\ge0.
\]
Hence (\ref{rho-est}) yields that
\[
\begin{array}{ll}
I_{even,4}(r,t)\le
&\d Ck^{1/2}\left(\frac{t+r+2k}{k}\right)^{1/2-(n-1)p/2+a_1}\times\\
&\d\times\int_{0}^{t-r}\!\!\!d\tau\int_{0}^{t-r-\tau}\!\!\!
\frac{\tau_{-}(\lambda,\tau)^{a_2}}{\sqrt{t-\tau-\lambda-r}}
\left(\log{4\frac{\tau_{+}(\lambda,\tau)}{\tau_{-}(\lambda,\tau)}}\right)^{a_3}
d\lambda.
\end{array}
\]
Changing variables by (\ref{alpha_beta}), we have 
\[
\begin{array}{ll}
I_{even,4}(r,t)\le
&\d Ck^{1/2}\left(\frac{t+r+2k}{k}\right)^{1/2-(n-1)p/2+a_1}\times\\
&\d\times\int_{-k}^{t-r}\left(\frac{\beta+2k}{k}\right)^{a_2}d\beta
\int_{\beta}^{t-r}\left(\log{4\frac{\alpha+2k}{\beta+2k}}\right)^{a_3}
\frac{d\alpha}{\sqrt{t-r-\alpha}}.
\end{array}
\]
Therefore, applying the simple inequality
\[
\tau_{+}(r,t)^{1-(n-1)p/2}\le Cw(r,t)^{-1},
\]
we obtain the desired estimates by the same argument as for $I_{even,3}$.
The proof of Lemma \ref{lem:basic_est} is now completed. 
\hfill$\Box$
\vskip10pt
\par\noindent
{\bf Proof of Lemma \ref{lm:apriori}.} 
Since (\ref{L}) and (\ref{norm}) yield that 
\[
L(|U|^p)(x,t)\le \|U\|^pL(w^{-p})(x,t), 
\]
it is enough to show the inequality
\[
w(r,t)L(w^{-p})(x,t)\le Ck^2D(T).
\]
This is established by (\ref{basic_est}) with setting
\[
\left\{
\begin{array}{lll}
\d a_1=a_3=0,\ a_2=-pq & 
\mbox{if}\ \d p>\frac{n+1}{n-1},\\
\d a_1=a_2=0,\ a_3=p &
\mbox{if}\ \d p=\frac{n+1}{n-1},\\
\d a_1=-pq,\ a_2=a_3=0 &
\mbox{if}\ \d p<\frac{n+1}{n-1}.
\end{array}
\right.
\]
\hfill$\Box$
\vskip10pt
\par\noindent
{\bf Proof of Lemma \ref{lm:apriori_odd_u0}.}
Due to Huygens' principle in Lemma \ref{lem:huygens},
one can replace $\tau_{-}$ by $\d \tau_{-}\chi_{\{-k\le t-r\le k\}}$ in (\ref{basic_est}).
Then, the integral with respect to the variable of $\beta=\tau-\lambda$ is bounded.  
In order to establish (\ref{apriori_odd_u0}),
it is enough to show the inequality
\[
w(r,t)L\left(\tau_{+}^{-(n-1)(p-\nu)/2}w^{-\nu}
\chi_{\{-k\le t-r\le k\}}\right)\le C_{n,\nu,p}k^2E_{\nu}(T).
\]
This is established by (\ref{basic_est}) with setting
\[
\left\{
\begin{array}{lll}
\d a_1=a_3=0,\ a_2=-\nu q & 
\mbox{if}\ \d p>\frac{n+1}{n-1},\\
\d a_1=a_2=0,\ a_3=\nu &
\mbox{if}\ \d p=\frac{n+1}{n-1},\\
\d a_1=-\nu q,\ a_2=a_3=0 &
\mbox{if}\ \d p<\frac{n+1}{n-1}.
\end{array}
\right.
\]
\hfill$\Box$
\vskip10pt
\par\noindent
{\bf Proof of Lemma \ref{lm:apriori_even_2} and Lemma \ref{lem:apriori_even_3}.}
In order to prove (\ref{apriori_even_2}) and (\ref{apriori_even_3}),
it is enough to show inequalities
\[
w(r,t)L\left(\tau_{+}^{-(n-1)(p-\nu)/2}\tau_{-}^{\sigma}w^{-\nu}\right)
\le C_{n,\nu,p}k^2E_{\nu,a}(T)
\]
and
\[
w(r,t)L\left(\tau_{+}^{-(n-1)p/2}\tau_{-}^{\kappa}\right)\le C_{n,\nu,p}k^2F_{\nu}(T).
\]
If we set
\[
\left\{
\begin{array}{lll}
\d a_1=a_3=0,\ a_2=\mu & 
\mbox{if}\ \d p>\frac{n+1}{n-1},\\
\d a_1=0,\ a_2=\sigma,\ a_3=\nu &
\mbox{if}\ \d p=\frac{n+1}{n-1},\\
\d a_1=-\nu q,\ a_2=\sigma,\ a_3=0 &
\mbox{if}\ \d p<\frac{n+1}{n-1}
\end{array}
\right.
\]
in (\ref{basic_est}), we have (\ref{apriori_even_2}). 
If we set $a_1=a_3=0$, $a_2=\kappa$ in (\ref{basic_est}), we have (\ref{apriori_even_3}).
\hfill$\Box$
\section{Lower bound in odd space dimensions}
In this section, we prove Theorem \ref{thm:main1} in odd space dimensions.
It is obviously enough for this to show the following proposition. 

\begin{prop}
\label{prop:lifespan_odd}
Let $n=5,7,9,\cdots$. Assume {\rm (\ref{hypo_F})} and {\rm (\ref{hypo_data})}.
Then, there exists a positive constant $\e_0=\e_0(f,g,n,p,k)$ 
such that each of {\rm (\ref{IE0})} and {\rm (\ref{IE})}
admits a unique solution $u\in C^1(\R^n\times[0,T])$
as far as T satisfies
\begin{equation}
\label{lower_lifespan_odd}
T\le
\left\{
\begin{array}{llll}
\d c\e^{-2p(p-1)/\gamma(p,n)}\ \mbox{if}\ 1<p<p_0(n) \\
\d \exp\left(c\e^{-p(p-1)}\right)\ \mbox{if}\ p=p_0(n)
\end{array}
\right.
\end{equation}
for $0<\e \le \e_0$, where $c$ is a positive constant independent of $\e$.
\end{prop}
\par
Our purpose is to construct a solution of the integral equation (\ref{IE'})
as a limit of $\{U_l\}_{l\in\N}$ in $X$.
To end this, we define a closed subspace $Y_0$ in $X$ by 
\[
Y_0=\{U \in X\ :\ \|\nabla_{x}^\alpha U\|\le2M_0\e^p\ (|\alpha|\le1)\},
\]
where we set
\[
M_0=2^ppAC_{n,0,p}k^2C_0^p>0.
\]
Recall that $A$ is the one in (\ref{hypo_F}) and
$C_{n,0,p}$ is the one in (\ref{apriori_odd_u0}). 
We note that there exists a positive constant $C_0$ independent of $\e$
which satisfies that
\begin{equation}
\label{W_0_odd}
\left\|\frac{\tau_{+}^{(n-1)/2}}{w}\p W_0\right\| \le C_0\e
\quad\mbox{for}\ 0<\e\le1.
\end{equation}
It is easy to check that
this fact follows from the definition of the weight function in (\ref{weight})
and the decay estimate for $U_0=v$ in (\ref{decay_est_v_0odd}) and 
(\ref{decay_est_v_1odd}).
\vskip10pt
\par\noindent
{\bf Proof of Proposition \ref{prop:lifespan_odd}.}
First of all, we assume that
\begin{equation}
\label{ep}
0<\e\le1
\end{equation}
to make use of (\ref{W_0_odd}).
We shall show the convergence of $\{U_l\}_{l\in\N}$.
The boundedness of $\{U_l\}_{l\in\N}$,
\begin{equation}
\label{induction_odd}
\|U_l\|\le2M_0\e^p\quad (l\in \N),
\end{equation}
can be obtained by induction with respect to $l$ as follows.
Recall that $L$ is a positive and linear operator by its definition, (\ref{L}).
It follows from (\ref{hypo_F}), (\ref{sequence1}), (\ref{apriori_odd_u0}) with $\nu=0$
and (\ref{W_0_odd}) that
\begin{equation}
\label{U_1_odd}
\begin{array}{ll}
\|U_1\|
&\le\|L\{|F(2U_0)|\}\|\le2^pA\|L(|U_0|^p)\|\\
&\le2^pAC_{n,0,p}k^2\|\tau_{+}^{(n-1)/2}w^{-1}U_{0}\|^pE_{0}(T)\le M_0\e^p,
\end{array}
\end{equation}
where we have used $E_{0}(T)=1$ for $p>1$. 
Assume that $\|U_{l-1}\|\le 2M_0\e^p\ (l\ge2)$.
It follows from the simple estimate
\begin{equation}
\label{U_l_odd}
\begin{array}{ll}
|U_l|
&\le AL(|U_{l-1}+U_0|^p)\\
&\le 2^pA\{L(|U_{l-1}|^p)+L(|U_0|^p)\},
\end{array}
\end{equation}
(\ref{apriori}) and (\ref{apriori_odd_u0}) with $\nu=0$ that
\[
\|U_l\|\le 2^pAk^2\{C\|U_{l-1}\|^pD(T)
+C_{n,0,p}\|\tau_{+}^{(n-1)/2}w^{-1}U_0\|^pE_{0}(T)\}.
\]
Hence (\ref{W_0_odd}) and the assumption of the induction yield that
\[
\|U_l\|\le 2^pACk^2(2M_0\e^p)^pD(T)+M_0\e^p.
\]
This inequality shows (\ref{induction_odd}) provided 
\begin{equation}
\label{det_lifespan_1}
2^pACk^2(2M_0)^p\e^{p^2}D(T)\le M_0\e^p.
\end{equation}
\par
Next we shall estimate the differences of $\{U_l\}_{l\in\N}$ 
under the conditions, (\ref{ep}) and (\ref{det_lifespan_1})
ensuring the boundedness (\ref{induction_odd}). 
By virture of (\ref{hypo_F}),
there exists a $\theta\in(0,1)$ such that
\begin{equation}
\label{diff_u_odd_0}
\begin{array}{ll}
|U_{l+1}-U_l|&
=|L\{F(U_l+U_0)-F(U_{l-1}+U_0)\}|\\
&=|L\{F'(U_{l-1}+U_0+\theta(U_l-U_{l-1}))(U_l-U_{l-1})\}|\\
&\le pAL\{|U_{l-1}+U_0+\theta(U_l-U_{l-1})|^{p-1}|U_l-U_{l-1}|\}.
\end{array}
\end{equation}
Hence we have
\begin{equation}
\label{diff_u_odd}
|U_{l+1}-U_l|\le2^{p-1}pAL\{(|3W_l|^{p-1}+|U_0|^{p-1})|U_l-U_{l-1}|\}.
\end{equation}
H\"older's inequality (\ref{Holder})
and {\it a priori} estimate (\ref{apriori}) yield that 
\begin{equation}
\label{Holder_apriori_odd}
\begin{array}{l}
\|L(|3W_l|^{p-1}|U_l-U_{l-1}|)\|\\
=\|L\{(|3W_l|^{(p-1)/p}|U_l-U_{l-1}|^{1/p})^p\}\|\\
\le Ck^2\||3W_l|^{(p-1)/p}|U_l-U_{l-1}|^{1/p}\|^pD(T)\\
\le Ck^2\|3W_l\|^{p-1}D(T)\|U_l-U_{l-1}\|.
\end{array}
\end{equation}
We note that (\ref{induction_odd}) implies $\|W_l\|\le2M_0\e^p\ (l\in\N)$.
Moreover, (\ref{apriori_odd_u0}) with $\nu=1$ implies that
\begin{equation}
\label{support_odd}
\begin{array}{l}
\|L(|U_0|^{p-1}|U_l-U_{l-1}|)\|\\
\le C_{n,1,p}k^2\|\tau_{+}^{(n-1)/2}w^{-1}U_0\|^{p-1}\|U_l-U_{l-1}\|E_{1}(T).
\end{array}
\end{equation}
Since (\ref{induction_odd}) implies that $\|W_l\|\le 2M_0\e^p$ for $l\ge2$,
the convergence of $\{U_l\}_{l\in\N}$ follows from 
\[
\|U_{l+1}-U_l\| \le \frac{1}{2} \|U_l-U_{l-1}\|\quad\mbox{for}\ l\ge2
\]
provided 
\begin{equation}
\label{det_lifespan_2}
2^{p-1}pAk^2\{C(6M_0\e^p)^{p-1}D(T)+C_{n,1,p}(C_0\e)^{p-1}E_1(T)\}\le\frac{1}{2},
\end{equation}
In fact, we obtain 
\begin{equation}
\label{conver_seq_odd}
\|U_{l+1}-U_l\|\le \frac{1}{2^{l-1}}\|U_2-U_1\|\quad\mbox{for}\ l\ge2
\end{equation}
which implies the convergence of $\{U_l\}_{l\in\N}$.
\par
Now we shall show the convergence of $\{\partial_iU_l\}_{l\in\N}$ for $i=1,2,\cdots,n$
under the conditions, (\ref{ep}), (\ref{det_lifespan_1}) and (\ref{det_lifespan_2})
which ensure the convergence of $\{U_l\}_{l\in\N}$.
As before, the boundedness of $\{\partial_iU_l\}_{l\in\N}$ for $i=1,2,\cdots,n$, 
\begin{equation}
\label{derivative_est_odd}
\|\partial_iU_l\|\le 2M_0\e^p\quad(l\in\N,\ i=1,2,\cdots,n),
\end{equation}
can be obtained by induction as follows.
Similarly to (\ref{U_1_odd}), we have
\[
\begin{array}{ll}
\|\p_iU_1\|
&\le\|L(|F'(2U_0)2\p_iU_0|)\|\le2^ppA\|L(|U_0|^{p-1}|\p_iU_0|)\|\\
&\le2^ppA\|L(|\p W_0|^p)\|\le2^ppAC_{n,0,p}k^2\|\tau_{+}^{(n-1)/2}w^{-1}\p W_0\|^pE_{0}(T).
\end{array}
\]
Hence (\ref{W_0_odd}) and $E_{0}(T)=1$ for $p>1$ implies $\|\p_iU_1\|\le M_0\e^p$.
Assume that $\|\p_iU_{l-1}\|\le2M_0\e^p\ (l\ge2)$.
We note that this means $\|\p W_{l-1}\|\le2M_0\e^p\ (l\ge2)$.
It follows from (\ref{hypo_F}) and (\ref{sequence1}) that 
\[
\begin{array}{ll}
|\partial_i U_l|
&\le |L(|F'(U_{l-1}+U_0)||\partial_i(U_{l-1}+U_0)|)|\\
&\le pA|L(|U_{l-1}+U_0|^{p-1}|\partial_i(U_{l-1}+U_0)|)|\\
&\le 2^{p-1}pAL\{(|U_{l-1}|^{p-1}+|U_0|^{p-1})(|\partial_iU_{l-1}|+|\partial_iU_0|)\}.
\end{array}
\]
Similarly to (\ref{Holder_apriori_odd}) and (\ref{support_odd}),
we obtain that
\[
\begin{array}{l}
\|L(|U_{l-1}|^{p-1}|\partial_iU_{l-1}|)\|
\le Ck^2\|U_{l-1}\|^{p-1}\|\partial_iU_{l-1}\|D(T),\\
\|L(|U_{l-1}|^{p-1}|\partial_iU_0|)\|
\le C_{n,p-1,p}k^2\|U_{l-1}\|^{p-1}\|\tau_{+}^{(n-1)/2}w^{-1}\partial_iU_0\|E_{p-1}(T),\\
\|L(|U_0|^{p-1}|\partial_iU_{l-1}|)\|
\le C_{n,1,p}k^2\|\tau_{+}^{(n-1)/2}U_{0}w^{-1}\|^{p-1}\|\partial_iU_{l-1}\|E_{1}(T),\\
\|L(|U_0|^{p-1}||\partial_iU_0|)\|
\le C_{n,0,p}k^2\|\tau_{+}^{(n-1)/2}\p W_{0}w^{-1}\|^{p}E_{0}(T).
\end{array}
\]
Hence, we get  
\[
\begin{array}{rl}
\|\partial_i U_l\|\le
&2^{p-1}pAk^2\{C\|\p W_{l-1}\|^pD(T)\\
&\quad+C_{n,p-1,p}\|\p W_{l-1}\|^{p-1}\|\tau_{+}^{(n-1)/2}w^{-1}\p W_0\|E_{p-1}(T)\\
&+C_{n,1,p}\|\p W_{l-1}\|\|\tau_{+}^{(n-1)/2}w^{-1}\p W_0\|^{p-1}E_{1}(T)\\
&\quad+C_{n,0,p}\|\tau_{+}^{(n-1)/2}w^{-1}\p W_0\|^p\}\\
\le&2^{p-1}pAk^2\{C(2M_0\e^p)^pD(T)+C_{n,p-1,p}(2M_0\e^p)^{p-1}C_0\e E_{p-1}(T)\\
&+C_{n,1,p}(2M_0\e^p)(C_0\e)^{p-1}E_{1}(T)+C_{n,0,p}(C_0\e)^p\}.
\end{array}
\]
This inequality shows (\ref{derivative_est_odd}) provided 
\begin{equation}
\label{det_lifespan_3}
\begin{array}{ll}
(3/2)M_0\e^p\ge
&2^{p-1}pAk^2\{C(2M_0)^p\e^{p^2}D(T)\\
&+C_{n,p-1,p}(2M_0)^{p-1}C_0\e^{p^2-p+1}E_{p-1}(T)\\
&+C_{n,1,p}2M_0C_0^{p-1}\e^{2p-1}E_{1}(T)\}.
\end{array}
\end{equation} 
\par
Next we shall estimate the differences of $\{\partial_iU_l\}_{l\in\N}$
for $i=1,2,\cdots,n$, 
under the conditions, (\ref{ep}), (\ref{det_lifespan_1}), (\ref{det_lifespan_2})
and (\ref{det_lifespan_3})
which ensure the convergence of $\{U_l\}_{l\in\N}$
and the boundedness of $\{\p W_l\}_{l\in\N}$. 
(\ref{sequence1}) implies that
\[
\begin{array}{l}
|\p_iU_{l+1}-\p_iU_l|\\
=|L\{F'(U_l+U_0)(\p_iU_l+\p_iU_0)-F'(U_{l-1}+U_0)(\p_iU_{l-1}+\p_iU_0)\}|\\
\le L\{|F'(U_l+U_0)||\p_iU_l-\p_iU_{l-1}|\}\\
\quad+L\{|F'(U_l+U_0)-F'(U_{l-1}+U_0)||\p_iU_{l-1}+\p_iU_0|\}.
\end{array}
\]
Hence it follows from (\ref{hypo_F}) that
\[
\begin{array}{rl}
|\p_iU_{l+1}-\p_iU_l|\le
& pAL\{|U_l+U_0|^{p-1}|\p_iU_l-\p_iU_{l-1}|\}\\
& +pAL\{|U_l-U_{l-1}|^{p-1}|\p_iU_{l-1}+\p_iU_0|\}\\
\le & 2^{p-1}pAL\{(|U_l|^{p-1}+|U_0|^{p-1})|\p_iU_l-\p_iU_{l-1}|\}\\
& +pAL\{|U_l-U_{l-1}|^{p-1}(|\p_iU_{l-1}|+|\p_iU_0|)\}.
\end{array}
\]
Similarly to the proof of the convergence of $\{U_l\}_{l\in\N}$, we obtain that
\[
\begin{array}{l}
\|L\{|U_l|^{p-1}|\p_iU_l-\p_iU_{l-1}|\}\|
\le Ck^2\|U_l\|^{p-1}\|\p_iU_l-\p_iU_{l-1}\|D(T),\\
\|L\{|U_0|^{p-1}|\p_iU_l-\p_iU_{l-1}|\}\|\\
\quad\le C_{n,1,p}k^2\|\tau_{+}^{(n-1)/2}w^{-1}U_0\|^{p-1}
\|\p_iU_l-\p_iU_{l-1}\|E_{1}(T),\\
\|L\{|U_l-U_{l-1}|^{p-1}|\p_iU_{l-1}|\}\|
\le Ck^2\|U_l-U_{l-1}\|^{p-1}\|\p_iU_{l-1}\|D(T),\\
\|L\{|U_l-U_{l-1}|^{p-1}|\p_iU_0|\}\|\\
\quad\le C_{n,p-1,p}k^2\|U_l-U_{l-1}\|^{p-1}\|\tau_{+}^{(n-1)/2}w^{-1}\p_iU_0\|E_{p-1}(T).
\end{array}
\]
Therefore, due to (\ref{conver_seq_odd}),
all the assumptions imply that
\[
\begin{array}{l}
\|\p_iU_{l+1}-\p_iU_l\|\\
\le2^{p-1}pAk^2\{C\|W_l\|^{p-1}D(T)
+C_{n,1,p}\|\tau_{+}^{(n-1)/2}w^{-1}W_0\|^{p-1}E_{1}(T)\}\times\\
\quad\times\|\p_iU_l-\p_iU_{l-1}\|\\
\quad+pAk^2\{C\|\p W_{l-1}\|D(T)
+C_{n,p-1,p}\|\tau_{+}^{(n-1)/2}w^{-1}\p W_0\|E_{p-1}(T)\}\times\\
\quad\times\|U_l-U_{l-1}\|^{p-1}\\
\le2^{p-1}pAk^2\{C(2M_0\e^p)^{p-1}D(T)
+C_{n,1,p}(C_0\e)^{p-1}E_{1}(T)\}\|\p_iU_l-\p_iU_{l-1}\|\\
\quad+pAk^2\{C(2M_0\e^p)D(T)+C_{n,p-1,p}C_0\e E_{p-1}(T)\}(\|U_2-U_1\|2^{-(l-1)})^{p-1}.
\end{array}
\]
This inequality yields
\begin{equation}
\label{conver_pseq_odd}
\|\p_iU_{l+1}-\p_iU_l\|
\le\frac{1}{2}\|\p_iU_l-\p_iU_{l-1}\|
+\frac{N_0(\e,T)}{2^{(l-1)(p-1)}},
\end{equation}
where we set
\[
N_0(\e,T)=pAk^2\|U_2-U_1\|^{p-1}\{C(2M_0\e^p)D(T)+C_{n,p-1,p}C_0\e E_{p-1}(T)\}
\]
provided
\begin{equation}
\label{det_lifespan_4}
2^{p-1}pAk^2\{C(2M_0)^{p-1}\e^{p(p-1)}D(T)
+C_{n,1,p}(C_0)^{p-1}\e^{p-1}E_{1}(T)\}\le\frac{1}{2}.
\end{equation}
We note that (\ref{conver_pseq_odd}) implies
\[
\|\p_iU_{l+1}-\p_iU_l\|
\le\frac{1}{2^{l-1}}\|\p_iU_2-\p_iU_1\|
+\frac{N_0(\e,T)}{2^{(l-1)(p-1)}}\sum_{\nu=0}^{l-2}\frac{1}{2^\nu}
\quad\mbox{for}\ l\ge2.
\]
The convergence of $\{\p_iU_l\}_{l\in\N}\ (i=1,2,3,\cdots,n)$ follows from this estimate.
\par
In this way, the convergence of $\{U_l\}_{l\in\N}$ in $Y_0\subset X$
can be established if all the five conditions,
(\ref{ep}), (\ref{det_lifespan_1}), (\ref{det_lifespan_2}), (\ref{det_lifespan_3}),
(\ref{det_lifespan_4}), are satisfied.
In order to complete the proof of Proposition \ref{prop:lifespan_odd},
we shall fix $\e_0=\e_0(f,g,n,p,k)$ and $c$ 
in the statement of Proposition \ref{prop:lifespan_odd}. 
First, we propose a sufficient condition to (\ref{ep}), (\ref{det_lifespan_1}) 
as well as related factors in (\ref{det_lifespan_2}), (\ref{det_lifespan_3}), 
(\ref{det_lifespan_4}) to $D(T)$ by 
\begin{equation}
\label{det_c_odd_1}
2^{2p}3^{p-1}pACk^2M_0^{p-1}\e^{p(p-1)}D(T)\le1.
\end{equation}
\par\noindent
Next, we propose sufficient conditions to related factors in (\ref{det_lifespan_2}), 
(\ref{det_lifespan_3}), (\ref{det_lifespan_4}) to $E_1(T)$ and $E_{p-1}(T)$ 
according to $p$ by the following. 
\par
In the case of $\d 1<p<\frac{n+1}{n-1}$, such conditions are 
\begin{equation}
\label{det_c_odd_2}
2^{p+1}pAk^2C_{n,1,p}C_0^{p-1}\e^{p-1}\left(\frac{2T+3k}{k}\right)^{-q}\le1
\end{equation}
and
\begin{equation}
\label{det_c_odd_3}
1\ge 2^{2p-1}pAk^2C_{n,p-1,p}M_0^{p-2}C_0\e^{(p-1)^2}
\left(\frac{2T+3k}{k}\right)^{-(p-1)q}.
\end{equation}
If we put
\[
\begin{array}{ll}
c=\min
&\{(2^{2p}3^{p-1}pACk^2(M_0)^{p-1})^{-1},(2^{p+1}pAk^2C_{n,1,p}C_0^{p-1})^{-p},\\
&(2^{2p-1}pAk^2C_{n,p-1,p}M_0^{p-2}C_0)^{-p/(p-1)}\}>0,
\end{array}
\]
the inequality $\e^{p(p-1)}D(T)\le c$ implies (\ref{det_c_odd_1}), 
(\ref{det_c_odd_2}) and (\ref{det_c_odd_3}), because of $-pq<\gamma(n,p)/2$. 
Furthermore, one can readily check that $\e_0=1$ by (\ref{ep}). 
\par
In the case of $\d p=\frac{n+1}{n-1}$, let us fix $\delta$ with
\begin{equation}
\label{delta_set_1}
0<\delta<\frac{1}{p}.
\end{equation}
Then, similarly to the above, our conditions are 
\begin{equation}
\label{det_c_odd_4}
2^{p+1}pAC_{n,1,p}k^2C_0^{p-1}\e^{p-1}\left(\frac{2T+3k}{k}\right)^{\delta}\le1
\end{equation}
and
\begin{equation}
\label{det_c_odd_5}
1\ge 2^{2p-1}pAk^2C_{n,p-1,p}M_0^{p-2}C_0\e^{(p-1)^2}
\left(\frac{2T+3k}{k}\right)^{\delta(p-1)}.
\end{equation}
If we put 
\[
\begin{array}{ll}
c=\min
&\{(2^{2p}3^{p-1}pACk^2(M_0)^{p-1})^{-1},
(2^{p+1}pAk^2C_{n,1,p}C_0^{p-1})^{-p},\\
&(2^{2p-1}pAk^2C_{n,p-1,p}M_0^{p-2}C_0)^{-p/(p-1)}\}>0,
\end{array}
\]
the inequality $\e^{p(p-1)}D(T)\le c$ implies (\ref{det_c_odd_1}), 
(\ref{det_c_odd_4}) and (\ref{det_c_odd_5}), because of (\ref{delta_set_1}).
Furthermore, one can readily check that $\e_0=1$ by (\ref{ep}).
\par
Finally, in the case of $\d p>\frac{n+1}{n-1}$, our conditions are 
\begin{equation}
\label{e_0_1}
2^{p+1}pAk^2C_{n,1,p}C_0^{p-1}\e^{p-1}\le1
\end{equation}
and
\begin{equation}
\label{e_0_2}
2^{2p-1}pAk^2C_{n,p-1,p}M_0^{p-2}C_0\e^{(p-1)^2}\le1.
\end{equation}
If we put 
\[
c=(2^{2p}3^{p-1}pACk^2(M_0)^{p-1})^{-1},
\]
the inequality $\e^{p(p-1)}D(T)\le c$ implies (\ref{det_c_odd_1}). 
Furthermore, one can find that 
\[
\begin{array}{ll}
\e_0=\min
&\{1,(2^{p+1}pAk^2C_{n,1,p})^{-1/(p-1)}C_0^{-1},\\
&(2^{2p-1}pAk^2M_0^{p-2}C_{n,p-1,p}C_0)^{-1/(p-1)^2}\}>0,
\end{array}
\]
by (\ref{ep}), (\ref{e_0_1}) and (\ref{e_0_2}). 
Therefore, the proof of proposition \ref{prop:lifespan_odd} is completed.
\hfill$\Box$

\section{Lower bound in even space dimensions}
Similarly to the previous section,
we investigate the lower bound of the lifespan in the even dimensional case.
Our purpose is to show the following proposition.
\begin{prop}
\label{prop:lifespan_even}
Let $n=4,6,8,\cdots$. Suppose that the same assumptions 
in Theorem \ref{thm:main1}, {\rm (\ref{hypo_F})} and {\rm (\ref{hypo_data})},
are fulfilled. 
Then, there exists a positive constant $\e_0=\e_0(f,g,n,p,k)$ 
such that each of {\rm (\ref{NWIVP0})} and {\rm (\ref{NWIVP})} admits a unique 
classical solution $u\in C^2(\R^4\times[0,T])$ if $n=4$ and $p=p_0(4)=2$,
or each of {\rm (\ref{IE0})} and {\rm (\ref{IE})} admits a unique solution 
$u\in C^1(\R^n\times[0,T])$ otherwise, as far as T satisfies
\begin{equation}
\label{lower_lifespan_even}
T\le
\left\{
\begin{array}{lllll}
\d c\e^{-2p(p-1)/\gamma(p,n)}\ \mbox{if}\ 1<p<p_0(n) \\
\d \exp\left(c\e^{-p(p-1)}\right)\ \mbox{if}\ p=p_0(n)
\end{array}
\right.
\end{equation}
for $0<\e \le \e_0$ and $c$ is a positive constant independent of $\e$.
\end{prop}
\par
Employing the similar argument to odd dimensions,
we shall construct a $C^1$ solution of the integral equation (\ref{IE'})
as a limit of $\{U_l\}_{l\in \N}$ in $X$. 
We also remark that it is possible to construct a $C^2$ solution 
if and only if $(n,p)=(4,2)$ in our problem.
However, its construction is almost the same as for $C^1$ solution.  
Therefore we shall omit it.
Now, define a closed subspace $Y_1$ in $X$ by 
\[
Y_1=\{U \in X : \|\nabla_{x}^\alpha U\|\le M_1 \e^p\ (|\alpha|\le1)\},
\]
where $M_1$ is defined by 
\[
M_1=2^{2p}pAk^2(C_{n,0,p}C_{1}^p+C_{n,1,p}C_{1}^{p-1}C_2
+C_{n,p-1,p}C_{1}C_{2}^{p-1}+C_{n,0,p}C_{2}^p) >0.
\]
Recall that $A$ is the one in (\ref{hypo_F}) and
that $C_{n,0,p}$, $C_{n,1,p}$ and $C_{n,p-1,p}$ 
are the one in (\ref{apriori_even_2}) and (\ref{apriori_even_3}).
$C_1$ and $C_2$ are positive constants 
independent of $\e$ which satisfy that 
\begin{equation}
\label{W_0_even}
\left\|(\tau_{+}\tau_{-})^{(n-1)/2}\frac{\p W_{00}}{w}\right\| \le C_1\e,\quad
\left\|(\tau_{+}\tau_{-})^{(n-1)/2}\frac{\p W_{01}}{w\tau_{-}}\right\| \le C_2\e^p. 
\end{equation}
The existence of $C_1$ and $C_2$ is trivial by
the definition of the weight function in (\ref{weight})
and the decay estimate for $U_{00}=v_0$ and $U_{01}=v_1$ 
in (\ref{decay_est_v_0even}) and (\ref{decay_est_v_1even}).
\vskip10pt
\par\noindent
{\bf Proof of the Proposition \ref{prop:lifespan_even}.} 
First of all, we shall show the convergence of $\{U_l\}_{l\in\N}$. 
The boundedness of $\{U_l\}_{l\in\N}$ in $Y_1$, 
\begin{equation}
\label{induction_even_1}
\|U_l\|\le2M_1\e^p\ (l\in \N),
\end{equation}
can be obtained by induction with respect to $l$ as follows. 
Similarly to (\ref{U_1_odd}), it follows from (\ref{sequence1}) 
and (\ref{apriori_even_2}) of Lemma \ref{lm:apriori_even_2} with $\nu=0$ 
and $a=0,1$ that 
\begin{equation}
\label{U_1_even}
\begin{array}{ll}
\|U_1\|\le
&\d 2^{2p}Ak^2C_{n,0,p}
\left(\left\|(\tau_{+}\tau_{-})^{(n-1)/2}
\frac{U_{00}}{w}\right\|^pE_{0,0}(T)\right.+\\
&\d+\left.\left\|(\tau_{+}\tau_{-})^{(n-1)/2}
\frac{U_{01}}{w\tau_{-}}\right\|^pE_{0,1}(T)\right).
\end{array}
\end{equation}
By the definition of (\ref{E(nu,a)_1}), (\ref{E(nu,t)_2})
and (\ref{E(nu,t)_3}) with $\nu=a=0$, we get 
\begin{equation}
\label{det_U_1_E_F}
E_{0,0}(T)=1.
\end{equation}
It follows from (\ref{W_0_even}) that 
\[
\|U_1\|\le2^{2p}Ak^2C_{n,0,p}\e^p(C_1^p+C_2^p\e^{p(p-1)}E_{0,1}(T)).
\]
This inequality show $\|U_1\|\le M_1\e^p$ provided 
\begin{equation}
\label{det_lifespan_even_0}
\e^{p(p-1)}E_{0,1}(T)\le 1.
\end{equation}
Assume that $\|U_{l-1}\|\le 2M_1\e^p$ $(l\ge2)$.
Making use of (\ref{U_l_odd}),
the assumption of the induction and (\ref{apriori}) yield that
\[
\begin{array}{ll}
\|U_l\|
&\d\le 2^pAk^2\left\{C\|U_{l-1}\|^pD(T)+
C_{n,0,p}\left\|(\tau_{+}\tau_{-})^{(n-1)/2}\frac{U_{00}}{w}\right\|^pE_{0,0}(T)\right.\\
&\d\qquad\left.+C_{n,0,p}\left\|(\tau_{+}\tau_{-})^{(n-1)/2}
\frac{U_{01}}{w\tau_{-}}\right\|^pE_{0,1}(T)\right\}\\
&\le 2^pACk^2(2M_1\e^p)^pD(T)+M_1\e^p.
\end{array}
\]
This inequality shows (\ref{induction_even_1}) provided 
\begin{equation}
\label{det_lifespan_even1}
2^pACk^2(2M_1)^p\e^{p^2}D(T)\le M_1\e^p.
\end{equation} 
\par
Next we shall estimate the differences under (\ref{det_lifespan_even_0}) 
and (\ref{det_lifespan_even1}) ensuring the boundedness (\ref{induction_even_1}). 
We note that (\ref{induction_even_1}) implies $\|W_l\|\le2M_1\e^p\ (l\in\N)$.
Making use of the inequalities (\ref{diff_u_odd_0}) 
and (\ref{diff_u_odd}), we have 
\[
|U_{l+1}-U_l|\le2^{p-1}pAL\{\{|3W_l|^{p-1}
+2^{p-1}(|U_{00}|^{p-1}+|U_{01}|^{p-1})\}|U_l-U_{l-1}|\}.
\]
Applying (\ref{apriori_even_2}) of Lemma \ref{lm:apriori_even_2} with $\nu=1$ and $a=0,1$,
we get 
\begin{equation}
\label{diff_U00_1}
\begin{array}{ll}
\|L(|U_{0a}|^{p-1}|U_l-U_{l-1}|)\|\\
\le k^2C_{n,1,p}E_{1,a}(T)
\d \left\|(\tau_{-}\tau_{+})^{(n-1)/2}
\frac{U_{0a}}{w\tau_{-}^a}\right\|^{p-1}\|U_l-U_{l-1}\|
\end{array}
\end{equation}
for $a=0,1$.
Since (\ref{induction_even_1}) implies that $\|W_l\|\le 2M_1\e^p$ for $l\ge2$,
the convergence of $\{U_l\}_{l\in\N}$ follows from 
\[
\|U_{l+1}-U_l\| \le \frac{1}{2} \|U_l-U_{l-1}\|\quad\mbox{for}\ l\ge2
\]
provided 
\begin{equation}
\label{det_lifespan_even2}
\begin{array}{ll}
\d 2^{p-1}pAk^2\{C(6M_1\e^p)^{p-1}D(T)+2^{p-1}C_{n,1,p}(C_1\e)^{p-1}
E_{1,0}(T)\\
\d\qquad+2^{p-1}C_{n,1,p}C_2^{p-1}E_{1,1}(T)\e^{p(p-1)}\}\le\frac{1}{2},
\end{array}
\end{equation}
Thus, we obtain   
\begin{equation}
\label{conver_seq_even2}
\|U_{l+1}-U_l\|\le \frac{1}{2^{l-1}}\|U_2-U_1\|\quad\mbox{for}\ l\ge2
\end{equation}
which implies the convergence of $\{U_l\}_{l\in\N}$. 
\par
Now we shall show the convergence of $\{\partial_iU_l\}_{l\in\N}$ for $i=1,2,\cdots,n$
under the conditions, (\ref{det_lifespan_even_0}), 
(\ref{det_lifespan_even1}), (\ref{det_lifespan_even2}) 
which ensure the convergence of $\{U_l\}_{l\in\N}$.
As before, the boundedness
\begin{equation}
\label{derivative_est_even}
\|\partial_iU_l\|\le 2M_1\e^p\quad(l\in\N,\ i=1,2,\cdots,n)
\end{equation}
can be obtained by induction as follows.
Similarly to (\ref{U_1_even}), we have
\[
\begin{array}{ll}
\|\p_iU_1\|
&\le2^ppA\|L(|U_0|^{p-1}|\p_iU_0|)\|\\
&\le2^{2p-1}pA\|L((|U_{00}|^{p-1}+|U_{01}|^{p-1})(|\p_iU_{00}|+|\p_iU_{01}|))\|\\
&\le2^{2p-1}pA\|L(\p W_{00}^{p}+W_{00}^{p-1}\p W_{01}+
W_{01}^{p-1}\p W_{00}+\p W_{01}^p)\|.
\end{array}
\]
Applying (\ref{apriori_even_3}) of Lemma \ref{lm:apriori_even_2} with 
$\kappa=1$ and $\kappa=p-1$, we get 
\[
\begin{array}{ll}
\|L(W_{00}^{p-1}\p W_{01})\|\\
\d\quad\le C_{n,1,p}k^2\d \left\|(\tau_{-}\tau_{+})^{(n-1)/2}\frac{W_{00}}{w}\right\|^{p-1}
\left\|(\tau_{-}\tau_{+})^{(n-1)/2}\frac{\p W_{01}}{w\tau_{-}}\right\|F_{1}(T)
\end{array}
\]
and 
\[
\begin{array}{ll}
\|L(W_{01}^{p-1}\p W_{00})\|\\
\d\le C_{n,p-1,p}k^2\d \left\|(\tau_{-}\tau_{+})^{(n-1)/2}
\frac{W_{01}}{w\tau_{-}}\right\|^{p-1}
\left\|(\tau_{-}\tau_{+})^{(n-1)/2}\frac{\p W_{00}}{w}\right\|F_{p-1}(T).
\end{array}
\]
By virtue of (\ref{W_0_even}) and (\ref{det_U_1_E_F}), we obtain 
\[
\begin{array}{ll}
\|\p_iU_1\|
&\le 2^{2p-1}pAk^2\e^p\left\{C_{n,0,p}C_1^p+C_{n,1,p}C_1^{p-1}C_2\e^{p-1}F_1(T)\right.\\
&\quad+C_{n,p-1,p}C_2^{p-1}C_1\e^{(p-1)^2}F_{p-1}(T)\\
&\left.\quad+C_{n,0,p,}C_2^p\e^{p(p-1)}E_{0,1}(T)\right\}.
\end{array}
\]
This inequality shows $\|\p_iU_1\|\le M_1\e^p$
provided (\ref{det_lifespan_even_0}), 
\begin{equation}
\label{det_lifespan_even3}
\e^{p-1}F_1(T)\le1
\end{equation}
and
\begin{equation}
\label{det_lifespan_even4}
\e^{(p-1)^2}F_{p-1}(T)\le1
\end{equation}
hold.
\par
Assume that $\|\p_iU_{l-1}\|\le2M_1\e^p\ (l\ge2)$.
Then we get 
\[
\begin{array}{ll}
|\partial_i U_l|
&\le pA|L(|U_{l-1}+U_0|^{p-1}|\partial_i(U_{l-1}+U_0)|)|\\
&\le 2^{p-1}pAL\{(|U_{l-1}|^{p-1}+|U_0|^{p-1})(|\partial_iU_{l-1}|+|\partial_iU_0|)\}\\
&\le 2^{p-1}pAL\{(|U_{l-1}|^{p-1}+2^{p-1}(|U_{00}|^{p-1}+|U_{01}|^{p-1}))\\
&\quad\times(|\partial_iU_{l-1}|+|\partial_iU_{00}|+|\p_iU_{01}|)\}.
\end{array}
\]
Similarly to (\ref{diff_U00_1}), we obtain that
\[
\begin{array}{l}
\|L(|U_{l-1}|^{p-1}|\partial_iU_{0a}|)\|\\
\quad \le C_{n,p-1,p}k^2\|U_{l-1}\|^{p-1}
\d\left\|(\tau_{-}\tau_{+})^{(n-1)/2}\frac{\p_iU_{0a}}{w\tau_{-}^a}\right\|E_{p-1,a}(T),\\
\|L(|U_{0a}|^{p-1}|\partial_iU_{l-1}|)\|\\
\quad\d\le C_{n,1,p}k^2\|\p_iU_{l-1}\|
\left\|(\tau_{-}\tau_{+})^{(n-1)/2}\frac{U_{0a}}{w\tau_{-}^a}\right\|^{p-1}E_{1,a}(T)\\
\end{array}
\]
for $a=0,1$. 
Making use of (\ref{W_0_even}) and (\ref{det_U_1_E_F}), we get  
\[
\begin{array}{rl}
&\|\partial_i U_l\|\\
&\le2^{p-1}pAk^2\lbrack C(2M_1\e^p)^pD(T)+C_{n,p-1,p}(2M_1\e^p)^{p-1}C_1\e E_{p-1,0}(T)\\
&\quad+C_{n,p-1,p}(2M_1\e^p)^{p-1}C_2\e^pE_{p-1,1}(T)\\
&\quad+2^{p-1}\left\{C_{n,1,p}(C_1\e)^{p-1}
(2M_1\e^p)E_{1,0}(T)+C_{n,0,p}C_1^p\e^p\right.\\
&\quad+C_{n,1,p}(C_1\e)^{p-1}C_2\e^pF_1(T)
+C_{n,1,p}(C_2\e^p)^{p-1}(2M_1\e^p)E_{1,1}(T)\\
&\quad\left.+C_{n,p-1,p}(C_2\e^p)^{p-1}C_1\e F_{p-1}(T)
+C_{n,0,p}(C_2\e^p)^pE_{0,1}(T)\right\}\rbrack.
\end{array}
\]
Under the assumptions (\ref{det_lifespan_even_0}), (\ref{det_lifespan_even3}) 
and (\ref{det_lifespan_even4}), this inequality shows (\ref{derivative_est_even}) provided 
\begin{equation}
\label{det_lifespan_even5}
\begin{array}{ll}
&(7/4)M_1\e^p\\
&\ge2^{p-1}pAk^2\lbrack C(2M_1)^p\e^{p^2}D(T)\\
&\quad+C_{n,p-1,p}(2M_1)^{p-1}C_1\e^{p^2-p+1}E_{p-1,0}(T)\\
&\quad+C_{n,p-1,p}(2M_1)^{p-1}C_2\e^{p^2}E_{p-1,1}(T)
+2^{p-1}\left\{C_{n,1,p}C_1^{p-1}\e^{2p-1}\right.\\
&\quad\times 2M_1E_{1,0}(T)+C_{n,1,p}C_2^{p-1}2M_1\e^{p^2}E_{1,1}(T)\}\rbrack.
\end{array}
\end{equation} 
\par
Next we shall estimate the differences under the conditions, 
(\ref{det_lifespan_even_0}), (\ref{det_lifespan_even1}), 
(\ref{det_lifespan_even2}), (\ref{det_lifespan_even3}), 
(\ref{det_lifespan_even4}), (\ref{det_lifespan_even5})
which ensure the convergence of $\{U_l\}_{l\in\N}$
and the boundedness of $\{\p U_l\}_{l\in\N}$. 
Similarly to odd dimensions, we have
\[
\begin{array}{rl}
&|\p_iU_{l+1}-\p_iU_l|\\
&\le  2^{p-1}pAL\{\{|U_l|^{p-1}+2^{p-1}\{|U_{00}|^{p-1}+|U_{01}|^{p-1}\}\}
|\p_iU_l-\p_iU_{l-1}|\}\\
& \quad+pAL\{|U_l-U_{l-1}|^{p-1}(|\p_iU_{l-1}|+|\p_iU_{00}|+|\p_iU_{01}|)\}.
\end{array}
\]
Applying (\ref{apriori_even_2}) of Lemma \ref{lm:apriori_even_2} with $\nu=1$, $\nu=p-1$ 
and $a=0,1$, we get 
\[
\begin{array}{l}
\|L\{|U_{0a}|^{p-1}|\p_iU_l-\p_iU_{l-1}|\}\|\\
\d\quad\le C_{n,1,p}k^2\left\|(\tau_{-}\tau_{+})^{(n-1)/2}
\frac{U_{0a}}{w\tau_{-}^a}\right\|^{p-1}
\|\p_iU_l-\p_iU_{l-1}\|E_{1,a}(T),\\
\|L\{|U_l-U_{l-1}|^{p-1}|\p_iU_{0a}|\}\|\\
\quad\le C_{n,p-1,p}k^2\|U_l-U_{l-1}\|^{p-1}
\d \left\|(\tau_{-}\tau_{+})^{(n-1)/2}\frac{\p_iU_{0a}}{w\tau_{-}^a}\right\|
E_{p-1,a}(T)
\end{array}
\]
for $a=0,1$.
Then, all the assumptions imply that
\[
\begin{array}{rl}
&\|\p_iU_{l+1}-\p_iU_l\|\\
&\le2^{p-1}pAk^2\lbrack\{C(2M_1\e^p)^{p-1}D(T)
+2^{p-1}\{C_{n,1,p}(C_1\e)^{p-1}E_{1,0}(T)\\
&\quad+C_{n,1,p}(C_2\e^p)^{p-1}E_{1,1}(T)\}\rbrack\|\p_iU_l-\p_iU_{l-1}\|\\
&\quad+pAk^2\{C(2M_1\e^p)D(T)+C_{n,p-1,p}C_1\e E_{p-1,0}(T)\\
&\quad+C_{n,p-1,p}C_2\e^{p}E_{p-1,1}(T)\}(\|U_2-U_1\|2^{-(l-1)})^{p-1}.
\end{array}
\]
This inequality yields
\begin{equation}
\label{conver_pseq_even3}
\|\p_iU_{l+1}-\p_iU_l\|
\le\frac{1}{2}\|\p_iU_l-\p_iU_{l-1}\|
+\frac{N_1(\e,T)}{2^{(l-1)(p-1)}},
\end{equation}
where we set
\[
\begin{array}{ll}
N_1(\e,T)&=pAk^2\{C(2M_1\e^p)D(T)+C_{n,p-1,p}C_1\e E_{p-1,0}(T)\\
&\quad+C_{n,p-1,p}C_2\e^{p}E_{p-1,1}(T)\}
\end{array}
\]
provided
\begin{equation}
\label{det_lifespan_even6}
\begin{array}{l}
2^{p-1}pAk^2\lbrack\{C(2M_1)^{p-1}\e^{p(p-1)}D(T)\\
\d\qquad+2^{p-1}\{C_{n,1,p}(C_1)^{p-1}\e^{p-1}E_{1,0}(T)\\
\d\qquad+C_{n,1,p}C_{2}^{p-1}\e^{p(p-1)}E_{1,1}(T)\}\rbrack\le\frac{1}{2}.
\end{array}
\end{equation}
Hence, the convergence of $\{\p_iU_l\}_{l\in\N}\ (i=1,2,3,\cdots,n)$ follows 
from this estimate. 
\par
In this way, the convergence of $\{U_l\}_{l\in\N}$ in $Y_0\subset X$
can be established if all the seven conditions,
(\ref{det_lifespan_even_0}), (\ref{det_lifespan_even1}), 
(\ref{det_lifespan_even2}), (\ref{det_lifespan_even3}), (\ref{det_lifespan_even4}), 
(\ref{det_lifespan_even5}), (\ref{det_lifespan_even6}) are satisfied.
In order to complete the proof of Proposition \ref{prop:lifespan_even},
we shall fix $\e_0=\e_0(f,g,n,p,k)$ and $c$ 
in the statement of Proposition \ref{prop:lifespan_even}. 
First, we propose a sufficient condition to (\ref{det_lifespan_even1}) 
as well as related factors in (\ref{det_lifespan_even2}), (\ref{det_lifespan_even5}), 
(\ref{det_lifespan_even6}) to $D(T)$ as 
\begin{equation}
\label{det_c_even_1}
2^{2p-1}3^{p}pACk^2M_1^{p-1}\e^{p(p-1)}D(T)\le1.
\end{equation}
\par\noindent
Next we propose sufficient conditions to related factors in 
(\ref{det_lifespan_even_0}), (\ref{det_lifespan_even2}), 
(\ref{det_lifespan_even3}), (\ref{det_lifespan_even4}), 
(\ref{det_lifespan_even5}) and (\ref{det_lifespan_even6}) to
\[
E_{0,1}(T),\ E_{1,0}(T),\ E_{1,1}(T),\ F_1(T),\ 
F_{p-1}(T),\ E_{p-1,0}(T),\ E_{p-1,1}(T)
\]
up to $p$. 
\vskip10pt
\par\noindent
{\bf Conditions in the case of $\v{\d 1<p<\frac{n+1}{n-1}}$.}
\par\noindent
(i) Conditions from $E_{0,1}(T)$.
\par
In (\ref{E(nu,t)_3}), setting $\nu=0$ and $a=1$, we have that
$\sigma=-(n-3)p/2<-1$ when $n\ge6$ and $\sigma=\mu=-p/2>-1$ when $n=4$,
which imply 
\begin{equation}
\label{det_E_01}
E_{0,1}(T)=
\left\{
\begin{array}{lll}
\d 1 & 
\mbox{if}\ \d n\ge6,\\
\d \left(\frac{2T+3k}{k}\right)^{1-p/2} &
\mbox{if}\ \d n=4.
\end{array}
\right.
\end{equation}
Since $E_{0,1}(T)$ appears in (\ref{det_lifespan_even_0}), the conditions are 
\begin{eqnarray}
\label{det_c_even}
\e^{p(p-1)}\le1 & \mbox{if}\ n\ge6,\\
\label{det_c_even0}
\d \left(\frac{2T+3k}{k}\right)^{1-p/2}\e^{p(p-1)}\le 1 & \mbox{if}\ n=4. 
\end{eqnarray}
\par\noindent
(ii) Condition from $E_{1,0}(T)$.
\par
In (\ref{E(nu,t)_3}), setting $\nu=1$ and $a=0$, we have that 
$\sigma=-(n-1)(p-1)/2>-1$ and $\mu=-2q-1$, which imply
\[
E_{1,0}(T)=\left(\frac{2T+3k}{k}\right)^{-2q}. 
\]
Since $E_{1,0}(T)$ appears in (\ref{det_lifespan_even2}),
(\ref{det_lifespan_even5}) and (\ref{det_lifespan_even6}),
the condition is
\begin{equation}
\label{det_c_even1}
2^{2p-1}3pAk^2C_{n,1,p}C_1^{p-1}\e^{p-1}\left(\frac{2T+3k}{k}\right)^{-2q}\le1.
\end{equation}
(iii) Condition from $E_{1,1}(T)$.
\par
In (\ref{E(nu,t)_3}), setting $\nu=a=1$, we have that 
$\sigma=-(n-3)(p-1)/2>-1$ and $\mu=n-1-(n-2)p$, which imply
\[
E_{1,1}(T)=\left(\frac{2T+3k}{k}\right)^{n-(n-2)p}.
\]
Since $E_{1,1}(T)$ appears in (\ref{det_lifespan_even2}),
(\ref{det_lifespan_even5}) and (\ref{det_lifespan_even6}),
the condition is 
\begin{equation}
\label{det_c_even2}
2^{2p-1}3pAk^2C_{n,1,p}C_2^{p-1}\e^{p(p-1)}\left(\frac{2T+3k}{k}\right)^{n-(n-2)p}\le1.
\end{equation}
\par\noindent
(iv) Conditions from $F_1(T)$.
\par
In (\ref{F(nu,t)_4}), setting $\nu=1$, we have that
\[
\begin{array}{ll}
\d \kappa=1-\frac{n-1}{2}p<-1 & \mbox{if}\ n\ge6\ ,
\mbox{or}\ n=4\ \mbox{and}\ \d p>\frac{4}{3},\\
\kappa=-1 & \mbox{if}\ n=4\ \mbox{and}\ \d p=\frac{4}{3},\\
\kappa>-1 & \mbox{if}\ n=4\ \mbox{and}\ \d 1<p<\frac{4}{3},
\end{array}
\]
which imply
\begin{equation}
\label{det_derivU_1_F_2}
F_{1}(T)=
\left\{
\begin{array}{ll}
\d 1 & \mbox{if}\ n\ge6\ ,\mbox{or}\ n=4\ \mbox{and}\ \d p>\frac{4}{3},\\
\d \log\frac{2T+3k}{k} & \mbox{if}\ n=4\ \mbox{and}\ \d p=\frac{4}{3},\\
\d \left(\frac{2T+3k}{k}\right)^{2-3p/2} & \mbox{if}\ n=4\ \mbox{and}\ \d 1<p<\frac{4}{3}.
\end{array}
\right.
\end{equation}
Since $F_{1}(T)$ appears in (\ref{det_lifespan_even3}), the conditions are
\begin{eqnarray}
\label{det_c_even3}
\e^{p(p-1)}\le1 & \mbox{if} & n\ge6\ ,\mbox{or}\ n=4\ \mbox{and}\ \d p>\frac{4}{3},\\ 
\label{det_c_even4}
\d \log\frac{2T+3k}{k}\e^{p(p-1)}\le 1 & \mbox{if} & n=4\ \mbox{and}\ \d p=\frac{4}{3},\\
\label{det_c_even5}
\left(\frac{2T+3k}{k}\right)^{2-3p/2}\e^{p(p-1)}\le 1
& \mbox{if} & n=4\ \mbox{and}\ \d 1<p<\frac{4}{3}.
\end{eqnarray}
\par\noindent
(v) Condition from $F_{p-1}(T)$.
\par
In (\ref{F(nu,t)_4}), setting $\nu=p-1$,
we have $\kappa=-(n-3)p/2-1<-1$, which implies $F_{p-1}(T)=1$. 
Since $F_{p-1}(T)$ appears in (\ref{det_lifespan_even4}), the condition is  
\begin{equation}
\label{det_c_even6}
\e^{p(p-1)}\le1.
\end{equation}
\par\noindent
(vi) Condition from $E_{p-1,0}(T)$.
\par
In (\ref{E(nu,t)_3}), setting $\nu=p-1$ and $a=0$,
we have $\sigma=-(n-1)/2<-1$, which implies  
\[
E_{p-1,0}(T)=\left(\frac{2T+3k}{k}\right)^{-(p-1)q}.
\]
Since $E_{p-1,0}(T)$ appears in (\ref{det_lifespan_even5}), the condition is 
\begin{equation}
\label{det_c_even7}
2^{2p}5\cdot7^{-1}pAk^2C_{n,p-1,p}C_1M_1^{p-2}\e^{(p-1)^2}
\left(\frac{2T+3k}{k}\right)^{-(p-1)q}\le1.
\end{equation}
\par\noindent
(vii) Conditions from $E_{p-1,1}(T)$.
\par
In (\ref{E(nu,t)_3}), setting $\nu=p-1$ and $a=1$, we have that 
$\sigma=-(n-3)/2<-1$ when $n\ge 6$ and $\sigma>-1,\mu=-1/2-(p-1)q$ when $n=4$,
which imply
\begin{equation}
\label{det_E_p-11}
E_{p-1,1}(T)=
\left\{
\begin{array}{ll}
\d \left(\frac{2T+3k}{k}\right)^{-(p-1)q} &
\mbox{if}\ n\ge6,\\
\d \left(\frac{2T+3k}{k}\right)^{1/2-(p-1)q} &
\mbox{if}\ n=4.
\end{array}
\right.
\end{equation}
Since $E_{p-1,1}(T)$ appear in (\ref{det_lifespan_even5}), the conditions are
\begin{equation}
\begin{array}{lll}
\label{det_c_even8}
&\d 2^{2p}5\cdot7^{-1}pAk^2C_{n,p-1,p}C_2M_1^{p-2}\e^{p(p-1)}\\
\qquad&\d\times\left(\frac{2T+3k}{k}\right)^{-(p-1)q}\le1
& \mbox{if}\ n\ge6,
\end{array}
\end{equation}
\begin{equation}
\begin{array}{lll}
\label{det_c_even9}
&\d 2^{2p}5\cdot7^{-1}pAk^2C_{4,p-1,p}C_2M_1^{p-2}\e^{p(p-1)}\\
\qquad&\d\times\left(\frac{2T+3k}{k}\right)^{1/2-(p-1)q}\le1
& \mbox{if}\ n=4.
\end{array}
\end{equation}
\par
Now, we are in a position to summarize all the conditions in (i)-(vii) above.
Set
\begin{equation}
\label{ep_even}
\e_0=1
\end{equation}
Then, (\ref{ep_even}) implies (\ref{det_c_even}),
(\ref{det_c_even3}) and (\ref{det_c_even6}). 
In order to make that (\ref{det_c_even_1}) includes (\ref{det_c_even4}), 
we employ Lemma \ref{lm:log} with $\delta=\gamma(p,4)/2>0$ and $X=(2T+3k)/k>1$. 
Then, if we set
\[
\begin{array}{ll}
c=\min
&\{(2^{2p-1}3^{p}pACk^2(M_1)^{p-1})^{-1},1,(2^{2p-1}3pAk^2C_{n,1,p}C_1^{p-1})^{-p},\\
&(2^{2p-1}3pAk^2C_{n,1,p}C_2^{p-1})^{-1},\gamma(4/3,4)/2,\\
&(2^{2p}5\cdot7^{-1}pAk^2C_{n,p-1,p}C_1M_1^{p-2})^{-p/(p-1)},\\
&(2^{2p}5\cdot7^{-1}pAk^2C_{n,p-1,p}C_2M_1^{p-2})^{-1}\}>0,
\end{array}
\]
the inequality $\e^{p(p-1)}D(T)\le c$ implies (\ref{det_c_even_1}), 
(\ref{det_c_even0}), (\ref{det_c_even1}), (\ref{det_c_even2}), (\ref{det_c_even4}), 
(\ref{det_c_even5}), (\ref{det_c_even7}), (\ref{det_c_even8}) and (\ref{det_c_even9})
because of $1-p/2\le\gamma(p,4)/2$ in (\ref{det_c_even0}),
$-2pq\le\gamma(p,n)/2$ in (\ref{det_c_even1}),
$n-(n-2)p\le \gamma(p,n)/2$ in (\ref{det_c_even2}), 
$2-3p/2\le\gamma(p,4)/2$ in (\ref{det_c_even5}),
$-pq\le\gamma(p,n)/2$ in (\ref{det_c_even7}),
$-(p-1)q\le\gamma(p,n)/2$ in (\ref{det_c_even8})
and $1/2-(p-1)q\le\gamma(p,4)/2$ in (\ref{det_c_even9}). 
\vskip10pt
\par\noindent
{\bf Conditions in the case of $\v{\d p=\frac{n+1}{n-1}}$.}
\par\noindent
(i) Conditions from $E_{0,1}(T)$.
\par
In (\ref{E(nu,t)_2}), setting $\nu=0$ and $a=1$, we have that 
$\sigma=-(n-3)(n+1)/2(n-1)<-1$ when $n\ge6$ and $\sigma>-1$ when $n=4$, 
which imply
\begin{equation}
\label{det_E_01_1}
E_{0,1}(T)=
\left\{
\begin{array}{lll}
\d 1 & 
\mbox{if}\ \d n\ge6,\\
\d \left(\frac{2T+3k}{k}\right)^{1-p/2} &
\mbox{if}\ \d n=4.
\end{array}
\right.
\end{equation}
Since $E_{0,1}(T)$ appears in (\ref{det_lifespan_even_0}), the conditions are
\begin{eqnarray}
\label{det_c_even10}
\e^{p(p-1)}\le1\ \mbox{if}\ n\ge6,\\
\label{det_c_even11}
\d \left(\frac{2T+3k}{k}\right)^{1-p/2}\e^{p(p-1)}\le 1\ \mbox{if}\ n=4. 
\end{eqnarray}
\par\noindent
(ii) Condition from $E_{1,0}(T)$.
\par
In (\ref{E(nu,t)_2}), setting $\nu=1$ and $a=0$,
we have $\sigma=-(n-1)(p-1)/2=-1$, which implies
\[
E_{1,0}(T)=\left(\frac{2T+3k}{k}\right)^{\delta}.
\]
Since $E_{1,0}(T)$ appears in (\ref{det_lifespan_even2}),
(\ref{det_lifespan_even5}) and (\ref{det_lifespan_even6}), 
the conditions is
\begin{equation}
\label{det_c_even12}
2^{2p-1}3pAk^2C_{n,1,p}C_1^{p-1}\e^{p-1}\left(\frac{2T+3k}{k}\right)^{\delta}\le1.
\end{equation}
\par\noindent
(iii) Condition from $E_{1,1}(T)$.
\par
In (\ref{E(nu,t)_2}), setting $\nu=a=1$, we have $\sigma=-(n-3)/(n-1)>-1$, which implies
\[
E_{1,1}(T)=\left(\frac{2T+3k}{k}\right)^{2/(n-1)}.
\] 
\par\noindent
Since $E_{1,1}(T)$ appears in (\ref{det_lifespan_even2})
 (\ref{det_lifespan_even5}) and (\ref{det_lifespan_even6}), 
the condition is 
\begin{equation}
\label{det_c_even13}
2^{2p-1}3pAk^2C_{n,1,p}C_2^{p-1}\e^{p(p-1)}\left(\frac{2T+3k}{k}\right)^{2/(n-1)}\le1.
\end{equation}
\par\noindent
(iv) Condition from $F_1(T)$.
\par
In (\ref{F(nu,t)_4}), setting $\nu=1$, we have $\kappa<-1$, which implies
$F_{1}(T)=1$. Since, $F_{1}(T)$ appears in (\ref{det_lifespan_even3}), 
the condition is 
\begin{equation}
\label{det_c_even14}
\e^{p(p-1)}\le1.
\end{equation}
\par\noindent
(v) Condition from $F_{p-1,0}(T)$.
\par
In (\ref{F(nu,t)_4}), setting $\nu=p-1$, we have $\kappa<-1$, which implies
$F_{p-1}(T)=1$. Since, $F_{p-1}(T)$ appears in (\ref{det_lifespan_even4}), 
the condition is (\ref{det_c_even14}). 
\par\noindent
(vi) Condition from $E_{p-1,0}(T)$.
\par
In (\ref{E(nu,t)_2}), setting $\nu=p-1$ and $a=0$, we have
$\sigma=-(n-1)/2<-1$, which implies
\[
E_{p-1,0}(T)=\left(\frac{2T+3k}{k}\right)^{(p-1)\delta}.
\]
Since $E_{p-1,0}(T)$ appears in (\ref{det_lifespan_even5}), the condition is 
\begin{equation}
\label{det_c_even15}
2^{2p}5\cdot7^{-1}pAk^2C_{n,p-1,p}C_1M_1^{p-2}\e^{(p-1)^2}
\left(\frac{2T+3k}{k}\right)^{(p-1)\delta}\le1.
\end{equation}
\par\noindent
(vii) Conditions from $E_{p-1,1}(T)$.
\par
In (\ref{E(nu,t)_2}), setting $\nu=p-1$ and $a=1$, we have that 
$\sigma=-(n-3)/2<-1$ when $n\ge 6$ and $\sigma>-1$ when $n=4$, 
which imply
\begin{equation}
\label{det_E_p-11_1}
E_{p-1,1}(T)=
\left\{
\begin{array}{ll}
\d \left(\frac{2T+3k}{k}\right)^{(p-1)\delta} &
\mbox{if}\ n\ge6,\\
\d \left(\frac{2T+3k}{k}\right)^{1/2} &
\mbox{if}\ n=4.
\end{array}
\right.
\end{equation}
Since $E_{p-1,1}(T)$ appears in (\ref{det_lifespan_even6}), the conditions are
\begin{equation}
\label{det_c_even16}
\begin{array}{lll}
&\d 2^{2p}5\cdot7^{-1}pAk^2C_{n,p-1,p}C_2M_1^{p-2}\e^{p(p-1)}\\
\qquad&\d\times\left(\frac{2T+3k}{k}\right)^{(p-1)\delta}\le1
& \mbox{if}\ n\ge6,\\
\end{array}
\end{equation}
\begin{equation}
\label{det_c_even17}
\begin{array}{lll}
&\d 2^{2p}5\cdot7^{-1}pAk^2C_{4,p-1,p}C_2M_1^{p-2}\e^{p(p-1)}\\
\qquad&\d\times\left(\frac{2T+3k}{k}\right)^{1/2}\le1.
& \mbox{if}\ n=4.
\end{array}
\end{equation}
\par
Now, we are in a position to summarize all the conditions in (i)-(vii) above.
First we note that (\ref{ep_even}) implies (\ref{det_c_even10}) and (\ref{det_c_even14}). 
Then, if we assume (\ref{delta_set_1}) and set
\[
\begin{array}{ll}
c=\min
&\{(2^{2p-1}3^{p}pACk^2(M_1)^{p-1})^{-1},1,(2^{2p-1}3pAk^2C_{n,1,p}C_1^{p-1})^{-p},\\
&(2^{2p-1}3pAk^2C_{n,1,p}C_2^{p-1})^{-1},\\
&(2^{2p}5\cdot7^{-1}pAk^2C_{n,p-1,p}C_1M_1^{p-2})^{-p/(p-1)},\\
&(2^{2p}5\cdot7^{-1}pAk^2C_{n,p-1,p}C_2M_1^{p-2})^{-1}
\}>0,
\end{array}
\]
the inequality $\e^{p(p-1)}D(T)\le c$ implies 
(\ref{det_c_even_1}), (\ref{det_c_even11}), (\ref{det_c_even12}), (\ref{det_c_even13}), 
(\ref{det_c_even15}), (\ref{det_c_even16}) and (\ref{det_c_even17})
because of $1-p/2\le1$ in (\ref{det_c_even11}),
$p\delta<1$ in (\ref{det_c_even12}) and (\ref{det_c_even15}), 
$2/(n-1)<1$ in (\ref{det_c_even13})
and $(p-1)\delta<1$ in (\ref{det_c_even16}). 
\vskip10pt
\par\noindent
{\bf Conditions in the case of $\v{\d p>\frac{n+1}{n-1}}$.}
\par\noindent
(i) Conditions from $E_{0,1}(T)$.
\par
In (\ref{E(nu,a)_1}), setting $\nu=0$ and $a=1$, we have that
$\mu=-(n-3)p/2<-1$ when $n\ge6$ and $\mu=-p/2\ge-1$ when $n=4$,
which imply
\begin{equation}
\label{det_E_01_2}
E_{0,1}(T)=
\left\{
\begin{array}{lll}
\d 1 & 
\mbox{if}\ \d n\ge6,\\
\d \left(\frac{2T+3k}{k}\right)^{1-p/2} &
\mbox{if}\ \d n=4,\ p<2,\\
\d \log\frac{2T+3k}{k} &
\mbox{if}\ \d n=4,\ p=2.
\end{array}
\right.
\end{equation}
Since $E_{0,1}(T)$ appear (\ref{det_lifespan_even_0}), the conditions are
\begin{eqnarray}
\label{det_c_even18}
\e^{p(p-1)}\le1\ & \mbox{if} & n\ge6,\\
\label{det_c_even19}
\d \left(\frac{2T+3k}{k}\right)^{1-p/2}\e^{p(p-1)}\le 1 & \mbox{if}& n=4,\ p<2,\\
\label{det_c_even20}
\d \log\frac{2T+3k}{k}\e^{p(p-1)}\le 1\ & \mbox{if} & n=4,\ p=2. 
\end{eqnarray}
\par\noindent
(ii) Condition from $E_{1,0}(T)$.
\par
In (\ref{E(nu,a)_1}), setting $\nu=1$ and $a=0$,
we have $\mu=-(n-1)(p-1)/2>-1$, which implies $\d E_{1,0}(T)=1$.
Since $E_{1,0}(T)$ appears in (\ref{det_lifespan_even2}),
(\ref{det_lifespan_even5}) and (\ref{det_lifespan_even6}),
the condition is 
\begin{equation}
\label{det_c_even21}
2^{2p-1}3pAk^2C_{n,1,p}C_1^{p-1}\e^{p-1}\le1.
\end{equation}
\par\noindent
(iii) Conditions from $E_{1,1}(T)$.
\par
In (\ref{E(nu,a)_1}), setting $\nu=a=1$,
we have $\mu=-(n-3)(p-1)/2-q=n-1-(n-2)p$, which implies
\begin{equation}
\label{det_E_11_2}
E_{1,1}(T)=
\left\{
\begin{array}{lll}
\d 1 & 
\mbox{if}\ \d p>\frac{n}{n-2},\\
\d \log\frac{2T+3k}{k} &
\mbox{if}\ \d p=\frac{n}{n-2},\\
\d \left(\frac{2T+3k}{k}\right)^{n-(n-2)p} &
\mbox{if}\ \d p<\frac{n}{n-2}.
\end{array}
\right.
\end{equation}
Since $E_{1,1}(T)$ appears in (\ref{det_lifespan_even2}),
(\ref{det_lifespan_even5}) and (\ref{det_lifespan_even6}),
the conditions are 
\begin{eqnarray}
\label{det_c_even22}
\left.
\begin{array}{l}
\d 2^{2p-1}3pAk^2C_{n,1,p}C_2^{p-1}\e^{p(p-1)}\times\\
\d \times\left(\frac{2T+3k}{k}\right)^{n-(n-2)p}
\end{array}
\right\}
\le1 & \mbox{if}\ \d p<\frac{n}{n-2},\\
\label{det_c_even23}
2^{2p-1}3pAk^2C_{n,1,p}C_2^{p-1}\e^{p(p-1)}\log\frac{2T+3k}{k}\le1
& \mbox{if}\ \d p=\frac{n}{n-2},\\
\label{det_c_even24}
2^{2p-1}3pAk^2C_{n,1,p}C_2^{p-1}\e^{p(p-1)}\le1
& \mbox{if}\ \d p>\frac{n}{n-2}.
\end{eqnarray}
\par\noindent
(iv) Condition from $F_1(T)$ and $F_{p-1}(T)$.
\par
In (\ref{F(nu,t)_4}), setting $\nu=1$ and $\nu=p-1$,
we have $\kappa<-1$, which implies $F_{1}(T)=F_{p-1}(T)=1$. 
Since $F_{1}(T)$ or $F_{p-1}(T)$ appears
in (\ref{det_lifespan_even3}) or (\ref{det_lifespan_even4})
respectively, the condition is 
\begin{equation}
\label{det_c_even25}
\e^{p(p-1)}\le1.
\end{equation} 
\par\noindent
(v) Condition from $E_{p-1,0}(T)$.
\par
In (\ref{E(nu,a)_1}), setting $\nu=p-1$ and $a=0$,
we have $\mu=-(n-1)/2-q(p-1)<-1$, which implies $E_{p-1,0}(T)=1$.
Since $E_{p-1,0}(T)$ appears in (\ref{det_lifespan_even5}), the condition is 
\begin{equation}
\label{det_c_even26}
2^{2p}5\cdot7^{-1}pAk^2C_{n,p-1,p}C_1M_1^{p-2}\e^{(p-1)^2}\le1.
\end{equation}
\par\noindent
(vi) Conditions from $E_{p-1,1}(T)$.
\par
In (\ref{E(nu,a)_1}), setting $\nu=p-1$ and $a=1$, we have that 
$\mu=-(n-3)/2-(p-1)q<-1$ when $n\ge 6$ and $\mu>-1$ when $n=4$, 
which imply 
\begin{equation}
\label{det_E_p-11_2}
E_{p-1,1}(T)=
\left\{
\begin{array}{ll}
\d 1 &
\mbox{if}\ n\ge6,\\
\d \left(\frac{2T+3k}{k}\right)^{1/2-(p-1)q} &
\mbox{if}\ n=4,\ p<2,\\
\d \log\frac{2T+3k}{k} &
\mbox{if}\ n=4,\ p=2.
\end{array}
\right.
\end{equation}
Since $E_{p-1,1}(T)$ appears in (\ref{det_lifespan_even6}), the condition are 
\begin{eqnarray}
\label{det_c_even27}
\d 2^{2p}5\cdot7^{-1}pAk^2C_{n,p-1,p}C_2M_1^{p-2}\e^{p(p-1)}\le1 & \mbox{if} & n\ge6,\\
\label{det_c_even28}
\left.
\begin{array}{r}
\d 2^{2p}5\cdot7^{-1}pAk^2C_{4,p-1,p}C_2M_1^{p-2}\e^{p(p-1)}\times\\
\d \times\left(\frac{2T+3k}{k}\right)^{1/2-(p-1)q}
\end{array}
\right\}
\le1 & \mbox{if} & n=4,\ p<2,\\
\label{det_c_even29}
\left.
\begin{array}{r}
\d 2^{2p}5\cdot7^{-1}pAk^2C_{4,p-1,p}C_2M_1^{p-2}\e^{p(p-1)}\times\\
\d \times\log\frac{2T+3k}{k}
\end{array}
\right\}
\le1 & \mbox{if} & n=4,\ p=2.
\end{eqnarray}
\par
Now, we are in a position to summarize all the conditions in (i)-(vi) above.
Set
\begin{equation}
\label{ep_det_fin}
\begin{array}{ll}
\e_0=\min
&\{1,(2^{2p-1}3pAk^2C_{n,1,p}C_1^{p-1})^{-1/(p-1)},\\
&\{2^{2p-1}3pAk^2C_{n,1,p}C_2^{p-1}\}^{-1/p(p-1)},\\
&\{2^{2p}5\cdot7^{-1}pAk^2C_{n,p-1,p}C_1M_1^{p-2}\}^{-1/(p-1)^2},\\
&\{2^{2p}5\cdot7^{-1}pAk^2C_{n,p-1,p}C_2M_1^{p-2}\}^{-1/p(p-1)}
\}>0.
\end{array}
\end{equation}
Then, (\ref{ep_det_fin}) implies (\ref{det_c_even18}), (\ref{det_c_even21}), 
(\ref{det_c_even24}), (\ref{det_c_even25}), (\ref{det_c_even26}) and (\ref{det_c_even27}).
In order to make that (\ref{det_c_even23}) includes (\ref{det_c_even_1}) when $n\ge6$, 
we employ the lemma \ref{lm:log} with $\delta=\gamma(p,n)/2>0$ and $X=(2T+3k)/k>1$. 
Then, if we set 
\[
\begin{array}{ll}
c=\min
&\{(2^{2p-1}3^{p}pACk^2(M_1)^{p-1})^{-1},1,(2^{2p-1}3pAk^2C_{n,1,p}C_2^{p-1})^{-1},\\
&(2^{2p-1}3pAk^2C_{n,1,p}C_2^{p-1})^{-1}\gamma(n/(n-2),n)/2,\\
&(2^{2p}5\cdot7^{-1}pAk^2C_{4,p-1,p}C_2M_1^{p-2})^{-1}
\}>0,
\end{array}
\]
the inequality $\e^{p(p-1)}D(T)\le c$ implies 
(\ref{det_c_even_1}), (\ref{det_c_even19}), (\ref{det_c_even20}), (\ref{det_c_even22}), 
(\ref{det_c_even23}) and (\ref{det_c_even28}) and (\ref{det_c_even29}), 
because of $1-p/2\le\gamma(p,4)/2$ in (\ref{det_c_even19}), 
$n-(n-2)p\le \gamma(p,n)/2$ in (\ref{det_c_even22}),  
and $1/2-(p-1)q\le\gamma(p,4)/2$ in (\ref{det_c_even28}).
Therefore the proof of proposition \ref{lower_lifespan_even} is now completed.
\hfill$\Box$
\section{Upper bound of the lifespan for the critical case in odd dimensions}
In this section, we prove Theorem \ref{thm:main2}
for the critical case in odd space dimensions. 
The proof is divided into three steps. 
In the first step, 
we get a point-wise estimate of the linear term $u^0$ from below
by means of the representation formula due to Rammaha \cite{R87, R88}. 
In the second step, we employ the comparison argument
between the solution of integral equations (\ref{IE})
and the blowing-up solution of ODE basically introduced by Zhou \cite{Z92_three},
in order to overcome the difficulty in the critical case. 
In the last step, we also employ the slicing method
introduced by Agemi, Kurokawa and Takamura \cite{AKT00}
which helps us to show the blow-up of the solution of such a modified ODE
arising from the high dimensional case.

\begin{prop}
\label{prop:lifespan_upper_odd}
Suppose that the assumptions of Theorem \ref{thm:main2} are fulfilled.
Let $u$ be a $C^0$-solution of (\ref{IE}) in $\R^n\times[0,T]$.
Then, there exists a positive constant $\e_0=\e_0(g,n,p,k)$
such that $T$ cannot be taken as 
\begin{equation}
\label{upper_lifespan_odd}
T>\d \exp\left(c\e^{-p(p-1)}\right)\ \mbox{if}\ p=p_0(n)
\end{equation}
for $0<\e \le \e_0$, where $c$ is a positive constant independent of $\e$.
\end{prop}

\par\noindent
{\bf Proof.} 
First we note that one may assume that the solution of (\ref{IE}) is radially symmetric
without loss of the generality.
To see this, we employ the spherical mean which is defined by 
\[
\wt{v}(r,t)=\frac{1}{\omega_n}\int_{|\omega|=1}v(r\omega,t)dS_{\omega}\ (r>0),
\]
for $v\in C(\R^n \times [0,\infty))$.
If we take the spherical mean of (\ref{IE}), we get
\[
\wt{u}=\e\wt{u^0}+\wt{L(|u|^p)}.
\]
Thanks to the fundamental identity for iterated spherical means, we have 
\[
\wt{L(|u|^p)}=L_{odd}(\wt{|u|^p}),
\]
where $L_{odd}$ is the one in (\ref{I_odd}).
See 78-81pp. of John \cite{J55} for details.
Thus, it follows from Jensen's inequality $\wt{|u|^p}\ge |\wt{u}|^p$ $(p>1)$
and the positivity of $L_{odd}$ that
\[
\wt{u}\ge\e\wt{u^0}+L_{odd}(|\wt{u}|^p).
\]
We estimate $\wt{u}$ from below all the time in this section,
so that we may assume that the equality holds here.
\par
Let $u=u(r,t)$ be a $C^0$-solution of 
\begin{equation}
\label{rad_IE}
u=\e u^0+L_{odd}(|u|^p)\quad\mbox{in}\quad(0,\infty)\times[0,T]
\end{equation}
which is associated by (\ref{IE}).
Note that $u^0=u^0(r,t)$ is a solution of
\begin{equation}
\label{u^0}
\left\{
\begin{array}{ll}
\d u^0_{tt}-\frac{n-1}{r}u^0_r-u^0_{rr}=0 & \mbox{in}\ (0,\infty)\times[0,\infty)\\
u^0(r,0)=0,\ u^0_t(r,0)=g(r), & r\in(0,\infty).
\end{array}
\right.
\end{equation}
\vskip10pt
\par\noindent
{\bf [The 1st step] Estimate of $\v{u^0}$.}
\par
We have the following representation of $u^0$. 
\begin{lem}[Rammaha \cite{R87}]
\label{lem:rad_u^0_odd}
Let $n=5,7,9\cdots$ and $u^0$ be a solution of (\ref{u^0}).
Then, $u^0$ is represented by
\[
u^0(r,t)=\frac{1}{2r^{(n-1)/2}}\int_{|r-t|}^{r+t}\lambda^{(n-1)/2}g(\lambda)P_{(n-3)/2}
\left(\frac{\lambda^2+r^2-t^2}{2r\lambda}\right)d\lambda,
\]
where $P_k$ is Legendre polynomial of degree $k$ defined by
\[
P_k(z)=\frac{1}{2^kk!}\frac{d^k}{dz^k}(z^2-1)^k.
\]
\end{lem}
See (6a) on 681p. in \cite{R87} for the proof.
This lemma implies the following estimate.
\begin{lem}[Rammaha \cite{R87}]
\label{lem:lower_v}
Let $n=5,7,9,\cdots$. Assume (\ref{blowup_asm}). 
Then there exists a positive constant $C_{g}$
such that for $t+k_0<r<t+k_1$ and $t\ge k_2$,
\begin{equation}
\label{lower_v}
u^0(r,t)\ge\frac{C_{g}}{r^{(n-1)/2}},
\end{equation}
where $\d k_2=k-k_0$. 
\end{lem}
See Lemma 2 on 682p. in \cite{R87} for the proof. 
\par
In order to prove the blow up result, we employ the iteration argument 
originally introduced by John \cite{J79}.
Our frame in the argument is obtained by the following lemma. 
\begin{lem}
\label{lm:frame_odd}
Let $u$ be a $C^0$-solution of (\ref{rad_IE}).
Assume (\ref{blowup_asm}). 
Then $u$ in $\Sigma_0=\left\{(r,t):\ 2k\le t-r\le r\right\}$ satisfies
\begin{equation}
\label{frame1_odd}
\begin{array}{l}
\d u(r,t)\ge\frac{C2^{(n-1)/2}(t-r)^{(n-1)/2}}{r^{(3n-7)/2}}\times\\
\d\times\int\!\!\!\int_{R(r,t)}\{(t-r-\tau+\lambda)(t+r-\tau-\lambda)\}^{(n-3)/2}
|u(\lambda,\tau)|^pd\lambda d\tau+\\
\qquad\d+\frac{E_1(t-r)^{(3n-5)/2-(n-1)p/2}}{{r}^{(3n-7)/2}}\e^p,
\end{array}
\end{equation}
where $C$ is the one in (\ref{I_odd}),
\[
E_1=\frac{CC_{g}^p(k_1-k_0)}{(n-1)2^{(n-1)p-(3n-9)/2}}
\]
and
\[
R(r,t)=\left\{(\lambda,\tau):\ t-r\le\lambda,
\tau+\lambda\le t+r, 2k\le \tau-\lambda\le t-r\right\}.
\]
\end{lem}
\par\noindent
{\bf Proof.} 
By virtue of Huygens' principle on $u^0$ and (\ref{I_odd}), we have 
\[
u\ge I_1+I_2\quad\mbox{in}\ \Sigma_0,
\]
where we set
\[
\begin{array}{l}
\d I_1(r,t)=Cr^{2-n}\int_{R(r,t)}(t-\tau)^{3-n}
h(\lambda,t-\tau,r)\lambda|u(\lambda,\tau)|^pd\lambda d\tau,\\
\d I_2(r,t)=Cr^{2-n}\int_{S(r,t)}(t-\tau)^{3-n}
h(\lambda,t-\tau,r)\lambda|u(\lambda,\tau)|^pd\lambda d\tau,\\
\d S(r,t)=\left\{(\lambda,\tau):\ t-r\le\lambda, \tau+\lambda\le t+r, 
-k_1\le \tau-\lambda\le -k_0\right\}.
\end{array}
\]
Changing variable by (\ref{alpha_beta}) in $I_1$, we have 
\[
\begin{array}{ll}
I_1(r,t)\ge
&\d \frac{Cr^{2-n}}{2}\int_{2k}^{t-r}
\left\{(t-r-\beta)(t+r-\beta)\right\}^{(n-3)/2} 
d\beta\times\\
&\d \times\int_{2(t-r)+\beta}^{t+r}
\left\{(\alpha-(t-r))(t+r-\alpha)\right\}^{(n-3)/2} \times\\
&\d \times(t-(\alpha+\beta)/2)^{3-n}(\alpha-\beta)|u(\lambda,\tau)|^p d\alpha
\end{array}
\]
in $\Sigma_0$. Noticing that 
\[
\begin{array}{lll}
\d t+r-\beta \ge 2r,\ t-\frac{\alpha+\beta}{2}\le r,\\  
\d \alpha-\beta \ge 2(t-r)\ \mbox{and}\ \alpha-(t-r)\ge t-r 
\end{array}
\]
hold in the domain of  the integral above, 
we have
\[
\begin{array}{ll}
I_1(r,t)\ge
&\d \frac{C2^{(n-3)/2}(t-r)^{(n-1)/2}}{r^{(3n-7)/2}}
\int_{2k}^{t-r}(t-r-\beta)^{(n-3)/2}d\beta \times\\
&\d \times\int_{2(t-r)+\beta}^{t+r}(t+r-\alpha)^{(n-3)/2}|u(\lambda,\tau)|^p d\alpha
\end{array}
\]
in $\Sigma_0$.
Hence, we have the first term of the right-hand side of (\ref{frame1_odd}). 
\par
Next, we shall show the second term of (\ref{frame1_odd}).
Similarly to the above, we have
\[
\begin{array}{ll}
I_2(r,t)\ge
&\d \frac{Cr^{2-n}}{2}\int_{-k_1}^{-k_0}
\left\{(t-r-\beta)(t+r-\beta)\right\}^{(n-3)/2}
d\beta\times\\
&\d \times\int_{2(t-r)+\beta}^{t+r}
\left\{(\alpha-(t-r))(t+r-\alpha)\right\}^{(n-3)/2} \times\\
&\d \times(t-(\alpha+\beta)/2)^{3-n}(\alpha-\beta)|u(\lambda,\tau)|^p d\alpha
\end{array}
\]
in $\Sigma_0$. 
Note that 
\[
\begin{array}{l}
\d t+r-\beta \ge r,\ t-\frac{\alpha+\beta}{2}\le r-\beta \le 2r,\\  
\d \alpha-(t-r) \ge t-r+\beta \ge t-r-k\ \mbox{and}\ t-r-\beta\ge t-r 
\end{array}
\]
hold in the domain of the integral above.
By making use of (\ref{lower_v}), we have
\[
\begin{array}{ll}
I_2(r,t)
&\ge\d\frac{CC_g^p(t-r)^{(n-3)/2}(t-r-k)^{(n-3)/2}}{2^{n-2}r^{(3n-7)/2}}\e^p
\int_{-k_1}^{-k_0}d\beta\times\\
&\quad\d\times\int_{2(t-r)+\beta}^{t+r}(\alpha-\beta)^{1-(n-1)p/2}
(t+r-\alpha)^{(n-3)/2}d\alpha\\
&\ge\d \frac{CC_{g}^p(t-r)^{(n-2)-(n-1)p/2}}{2^{(n-1)p-(3n-11)/2}r^{(3n-7)/2}}\e^p
\int_{-k_1}^{-k_0}d\beta\times\\
&\quad\d\times \int_{2(t-r)+\beta}^{3(t-r)}\{3(t-r)-\alpha\}^{(n-3)/2}d\alpha\\
\end{array}
\]
in $\Sigma_0$.
The second term of the right-hand side of (\ref{frame1_odd}) follows from this inequality.
Therefore, the proof of Lemma \ref{lm:frame_odd} is ended. 
\hfill$\Box$
\vskip10pt
\par\noindent
{\bf [The 2nd Step] Comparison argument.}
\par
Let us consider a solution $w$ of 
\begin{equation}
\label{v_equal}
\begin{array}{ll}
\d w(t-r)=
&\d\frac{C2^{(n-3)/2}(t-r)^{(n-1)/2}}{r^{(3n-7)/2}}
\int_{2k}^{t-r}(t-r-\beta)^{(n-3)/2}d\beta\\
&\d \times\int_{2(t-r)+\beta}^{t+r}(t+r-\alpha)^{(n-3)/2}|w(\beta)|^p d\alpha\\
&\d+\frac{E_1(t-r)^{(3n-5)/2-(n-1)p/2}}{2{r}^{(3n-7)/2}}\e^p.
\end{array}
\end{equation}
Then we have the following comparison lemma. 
\begin{lem}
\label{lem:comparison_arg}
Let $u$ be a solution of (\ref{rad_IE}) and $w$ be a solution of (\ref{v_equal}). 
Then, $u$ and $w$ satisfy 
\[
u>w\quad\mbox{in}\ \Sigma_0.
\]
\end{lem}
{\bf Proof.} 
Fix a point $(r_0,t_0)\in\Sigma_0$. Define
\[
\Lambda(r,t)=\left\{(\lambda,\tau)\in D(r,t):\ 2k\le \tau-\lambda\le \lambda\right\},
\]
where we set
\[
D(r,t)=\left\{(\lambda,\tau):\ t-r\le \tau+\lambda\le t+r, 
-k\le \tau-\lambda\le t-r\right\}
\]
which is the domain of the integral in (\ref{I_odd}). 
Let us consider $u$ and $v$ in $\Lambda(r_0,t_0)$. 
Note that $u>w$ on $\tau-\lambda=2k$
and at $(2k,4k)$ which is an edge of $\Sigma_0$. 
By compactness of the closure of $\Lambda(r_0,t_0)$, 
we have $u>w$ in a neighborhood of $\tau-\lambda=2k$ and $\lambda\ge 2k$. 
\par
Assume that there exist a point $(r_1,t_1)$ with $u(r_1,t_1)=w(t_1-r_1)$,
which is nearest to $(2k,4k)$ in such a neighborhood.
Since $u>w$ in $R(r_1,t_1)$,  we have 
\[
\begin{array}{l}
\d\frac{C2^{(n-3)/2}(t_1-r_1)^{(n-1)/2}}{r_1^{(3n-7)/2}}
\int\!\!\!\int_{R(r_1,t_1)}(t_1-r_1-\tau+\lambda)^{(n-3)/2}\times\\
\times(t_1+r_1-\tau-\lambda)^{(n-3)/2}\d|u(\lambda,\tau)|^pd\lambda d\tau
\d+\frac{E_1(t_1-r_1)^{(3n-5)/2-(n-1)p/2}}{{r_1}^{(3n-7)/2}}\e^p\\
>\d\frac{C2^{(n-3)/2}(t_1-r_1)^{(n-1)/2}}{r_1^{(3n-7)/2}}
\int\!\!\!\int_{R(r_1,t_1)}(t_1-r_1-\tau+\lambda)^{(n-3)/2}\\
\times(t_1+r_1-\tau-\lambda)^{(n-3)/2}|w(\tau-\lambda)|^pd\lambda d\tau
\d+\frac{E_1(t_1-r_1)^{(3n-5)/2-(n-1)p/2}}{{2r_1}^{(3n-7)/2}}\e^p,
\end{array}
\]
In view of (\ref{frame1_odd}) and (\ref{v_equal}),
this inequality implies that $u>w$ at $(r_1,t_1)$,
which is a contradiction to the definition of $(r_1,t_1)$.
Therefore, we have $u>w$ in $\Lambda(r_0,t_0)$.
$(r_0,t_0)$ stands for any point in $\Sigma_0$,
so that $\Lambda(r_0,t_0)$ covers all of $\Sigma_0$.
The proof is completed.
\hfill$\Box$
\par
We note that Lemma \ref{lem:comparison_arg} implies that
the lifespan of $w$ is greater than the one of $u$,
so that it is sufficient to look for the lifespan of $w$ in $\Sigma_0$.
By definition of $w$ in (\ref{v_equal}), we have
\[
\begin{array}{ll}
w(\xi)\ge
&\d \frac{C\xi^{3-n}}{2^{n-2}}\int_{2k}^{\xi}
(\xi-\beta)^{(n-3)/2}|w(\beta)|^pd\beta\\
&\d \times\int_{2\xi+\beta}^{3\xi}(3\xi-\alpha)^{(n-3)/2}d\alpha
+\frac{E_1}{2^{(3n-5)/2}}\xi^{-q-(n-1)/2}\e^p\\
\end{array}
\]
in $\Gamma_0$, where we set
\[
\xi=\frac{r}{2},\ \Gamma_0=\{t-r=\xi, r\ge 4k\}.
\]
Hence we obtain that
\[
w(\xi)\ge\frac{C\xi^{3-n}}{2^{n-3}(n-1)}\int_{2k}^{\xi}
(\xi-\beta)^{n-2}|w(\beta)|^pd\beta+\frac{E_1\xi^{-q-(n-1)/2}}{2^{(3n-5)/2}}\e^p
\]
for $\xi\ge2k$.
Then, it follows from the setting
\[
W(\xi)=\xi^{q+(n-1)/2}w(\xi)
\] 
that
\begin{equation}
\label{frame_sub_odd}
\d W(\xi)\ge D_n\xi^{q-(n-5)/2}\int_{2k}^{\xi}
\frac{(\xi-\beta)^{n-2}|W(\beta)|^pd\beta}{\beta^{(n-1)p/2+pq}}
+E_2\e^p\quad \mbox{for}\ \xi\ge 2k,
\end{equation}
where we set
\[
D_n=\frac{C}{2^{n-3}(n-1)},\ E_2=\frac{E_1}{2^{(3n-5)/2}}.
\]
Therefore we obtain the iteration frame in this section,
\begin{equation}
\label{frame3_odd}
W(\xi)\ge D_n\int_{2k}^{\xi}
\left(\frac{\xi-\beta}{\xi}\right)^{n-2}
\frac{|W(\beta)|^p}{\beta^{pq}}d\beta+E_2\e^p
\quad\mbox{for}\ \xi\ge 2k.
\end{equation}
\vskip10pt
\par\noindent
{\bf [The 3rd step] Slicing method in the iteration.} 
\par
Let us define a blow-up domain as follows.
Let us set
\[
\Gamma_j=\{\xi\ge l_{j}k\},\
l_j=2+\frac{1}{2}+\cdots+\frac{1}{2^j}\ (j\in\N).
\]
We shall use the fact that a sequence $\{l_j\}$
is monotonously increasing and bounded, $2<l_j<3$, 
so that $\Gamma_{j+1}\subset \Gamma_{j}$.
Assume an estimate of the form 
\begin{equation}
\label{j_times_odd}
W(\xi)\ge C_j \left(\log\frac{\xi}{l_jk}\right)^{a_j}\ \mbox{in}\ \Gamma_j
\end{equation}
where $a_j\ge 0$ and $C_j>0$.
Putting (\ref{j_times_odd}) into (\ref{frame3_odd}) and 
recalling that $pq=1$, we get an estimate in $\Gamma_{j+1}$ such as 
\[
\begin{array}{lll}
\d W(\xi)
&\ge&\d  D_nC_j^p\int_{l_jk}^{\xi}
\left(\frac{\xi-\beta}{\xi}\right)^{n-2}
\left(\log\frac{\beta}{l_jk}\right)^{pa_j}\frac{d\beta}{\beta}.
\end{array}
\]
Noting that $\d \frac{l_{j}}{l_{j+1}}\xi \ge l_jk$ in $\Gamma_{j+1}$, we have 
\[
\begin{array}{llll}
\d W(\xi)
&\ge&\d  D_nC_j^p\int_{l_jk}^{l_{j}\xi/l_{j+1}}
\left(\frac{\xi-\beta}{\xi}\right)^{n-2}
\left(\log\frac{\beta}{l_jk}\right)^{pa_j}\frac{d\beta}{\beta}\\
&\ge&\d D_nC_j^p\left(1-\frac{l_j}{l_{j+1}}\right)^{n-2}
\int_{l_jk}^{l_{j}\xi/l_{j+1}}
\left(\log\frac{\beta}{l_jk}\right)^{pa_j}\frac{d\beta}{\beta}\\
&=&\d \frac{D_nC_j^p}{pa_j+1}\left(1-\frac{l_j}{l_{j+1}}\right)^{n-2}
\d \left(\log\frac{\xi}{l_{j+1}k}\right)^{pa_j+1}.
\end{array}
\]
By monotonicity of $\{l_j\}$ and
\[
1-\frac{l_j}{l_{j+1}}
=\frac{l_{j+1}-l_{j}}{l_{j+1}}
=\frac{1}{2^{j+1}l_{j+1}} 
\ge\frac{1}{3\cdot 2^{j+1}},
\]
we finally obtain
\begin{equation}
\label{j+1_times_odd}
W(\xi)\ge C_{j+1}\left(\log\frac{\xi}{l_{j+1}k}\right)^{pa_j+1}\ \mbox{in}\ \Gamma_{j+1},
\end{equation}
where we set
\[
C_{j+1}=\frac{D_nC_j^p}{3^{n-2}\cdot2^{(j+1)(n-2)}(pa_j+1)}.
\]
\par
Now, we are in a position to define sequences in the iteration.
In view of (\ref{frame3_odd}), the first estimate is $W(\xi)\ge E_2 \e^p$,
so that, with the help of (\ref{j_times_odd}) and (\ref{j+1_times_odd}),
a sequence $\{a_j\}$ should be defined by 
\[
a_1=0,\ a_{j+1}=pa_j+1\ (j\in\N).
\]
Also a sequence $\{C_j\}$ should be defined by 
\[
C_1=E_2\e^p,\ C_{j+1}=\frac{D_nC_j^p}{3^{n-2}\cdot2^{(j+1)(n-2)}(pa_j+1)}\ (j\in\N).
\]
One can easily check that 
\[
a_j=\frac{p^{j-1}-1}{p-1}\ (j\in\N),
\]
which gives us
\[
\frac{1}{pa_j+1}\ge \frac{p-1}{p^{j}}.
\]
Thus one can find that 
\[
C_{j+1}\ge E\frac{C_j^p}{(2^{n-2}p)^{j}}\ (j\in\N),
\]
where $E$ is a positive constant defined by 
\[
E=\frac{D_n(p-1)}{6^{n-2}}.
\]
Hence we inductively obtain, for $j\ge2$, that 
\[
\log C_j\ge p^{j-1}\left\{\log C_1+\sum_{k=0}^{j-2}\frac{p^{k}\log E
-(j-1-k)p^{k}\log(2^{n-2}p)}{p^{j-1}}\right\}.
\]
The sum part of above inequality converges as $j\rightarrow \infty$ 
by d'Alembert's criterion.
It follows from this fact that there exist a constant $S$ 
independent of $j$ such that 
\[
C_j\ge \exp\{p^{j-1}(\log C_1+S)\}\quad\mbox{for}\ j\ge 2.
\]
Combining all the estimates above and making use of the monotonicity of $\Gamma_j$, 
we have the final inequality 
\[
\begin{array}{ll}
\d W(\xi)
&\ge\d \exp\{p^{j-1}(\log C_1+S)\}\left(\log\frac{\xi}{3k}\right)^{(p^{j-1}-1)/(p-1)}\\
&=\d \exp\{p^{j-1}I(\xi)\}\left(\log\frac{\xi}{3k}\right)^{-1/(p-1)}
\end{array}
\]
in $\d \Gamma_{\infty}=\{\xi \ge 3k\}$, where we set 
\[
I(\xi)=\log\left(e^SE_2\e^p\left(\log\frac{\xi}{3k}\right)^{1/(p-1)}\right).
\]
\par
If there exist a point $\xi_0 \in \Gamma_{\infty}\subset \Gamma_{j}\ (j\ge1)$
such that $I(\xi_0)>0$, 
we get $\d W(\xi_0)\rightarrow \infty$ as $j\rightarrow \infty$.
Note that $I(\xi_0)>0$ is equivalent to
\[
\xi_0>3k\exp\{(e^{S}E_2)^{-(p-1)}\e^{-p(p-1)}\}. 
\]
It is trivial that there exists a positive constant $\e_0=\e_0(g,n,p,k)$ such that
\[
\exp\{(e^{S}E_2)^{-(p-1)}\e^{-p(p-1)}\}\ge1
\quad\mbox{for}\ 0<\e\le\e_0.
\]
Since there exists $(r_0,t_0)\in\Sigma_0$ such that $t_0-r_0= \xi_0>3k$, 
we obtain the desired conclusion;
\[
T>3k\exp\{(e^{S}E_2)^{-(p-1)}\e^{-p(p-1)}\}.
\]
Therefore the proof of the critical case is now completed
with a minor modification on $\e_0$.
\hfill$\Box$
\section{Upper bound of the lifespan for the subcritical case in odd dimensions}
Similarly to the previous section,
we prove the blow-up result of solution for (\ref{IE})
in the subcritical case in odd space dimensions. 
Note that we do not have to make use of the slicing method.

\begin{prop}
\label{prop:lifespan_upper_sub_odd} 
Suppose that the assumptions of Theorem \ref{thm:main2} are fulfilled. 
Let $u$ be a $C^0$-solution of (\ref{IE}) in $\R^n\times[0,T]$. 
Then, there exists a positive constant $\e_0=\e_0(g,n,p,k)$ such that $T$ cannot 
be taken as 
\begin{equation}
\label{upper_lifespan_sub_odd}
T>\d c\e^{-2p(p-1)/\gamma(p,n)}\ \mbox{if}\ 1<p<p_0(n) \\
\end{equation}
for $0<\e \le \e_0$, where $c$ is a positive constant independent of $\e$.
\end{prop}
{\bf Proof.} 
Because of the fact that $-(n-1)p/2-pq<0$ for $n\ge5$, (\ref{frame_sub_odd}) yields 
\begin{equation}
\label{frame_sub_odd2}
W(\xi)\ge D_n\xi^{-(n-2)-pq}\int_{2k}^{\xi}
(\xi-\beta)^{n-2}|W(\beta)|^pd\beta
+E_2\e^p\quad \mbox{for}\ \xi\ge 2k.
\end{equation}
This is our iteration frame in this case.
\par
Assume an estimate of the form 
\begin{equation}
\label{j_times_sub}
W(\xi)\ge C_j \frac{(\xi-2k)^{a_j}}{\xi^{b_j}}\ \mbox{in}\ \Gamma_0,
\end{equation} 
where $a_j,b_j\ge 0$ and $C_j>0$.
Then, putting (\ref{j_times_sub}) into 
(\ref{frame_sub_odd2}), we get
\[
W(\xi)\ge\frac{D_nC_j^p}{\xi^{(n-2)+pq+pb_j}}\int_{2k}^{\xi}
(\xi-\beta)^{n-2}(\beta-2k)^{pa_j}d\beta
\quad\mbox{in}\ \Gamma_0.
\]
Applying the integration by parts to $\beta$-integral $(n-2)$ times,
we obtain that
\[
\begin{array}{lll}
&\d \frac{n-2}{pa_j+1}\int_{2k}^{\xi}(\xi-\beta)^{n-3}(\beta-2k)^{pa_j+1}d\beta\\
&\ge\d \frac{(n-2)(n-3)}{(pa_j+2)^2}\int_{2k}^{\xi}(\xi-\beta)^{n-4}
(\beta-2k)^{pa_j+2}d\beta\\
&\cdots\\
&\ge\d \frac{(n-2)!}{(pa_j+n-1)^{n-1}}(\xi-2k)^{pa_j+n-1}.
\end{array}
\]
Therefore we finally get 
\begin{equation}
\label{j+1_times_sub}
W(\xi)\ge C_{j+1}\frac{(\xi-2k)^{pa_j+n-1}}{\xi^{pb_j+n-2+pq}}\ \mbox{in}\ \Gamma_0,
\end{equation}
where we set
\[
C_{j+1}=\frac{D_nC_j^p(n-2)!}{(pa_j+n-1)^{n-1}}.
\]
\par
Now, we are in a position to define sequences in the iteration. In view 
of (\ref{frame_sub_odd2}), the first estimate is $W(\xi)\ge E_2 \e^p$,
so that, with the help of (\ref{j_times_sub}) and (\ref{j+1_times_sub}),
sequences $\{a_j\}$ and $\{b_j\}$ should be defined by 
\[
a_1=0,\ a_{j+1}=pa_j+n-1\ (j\in\N)
\]
and 
\[
b_1=0,\ b_{j+1}=pb_j+n-2+pq\ (j\in\N).
\]
Also a sequence $\{C_j\}$ should be defined by 
\[
C_1=E_2\e^p,\ C_{j+1}=\frac{D_nC_j^p(n-2)!}{(pa_j+n-1)^{n-1}}\ (j\in\N).
\]
One can readily check that 
\[
a_j=\frac{n-1}{p-1}(p^{j-1}-1),
\ b_j=\frac{pq+n-2}{p-1}(p^{j-1}-1)\ (j\ge\N),
\]
which gives us
\[
\frac{1}{pa_j+n-1}\ge \frac{p-1}{p^{j}(n-1)}.
\]
Hence one can find that 
\[
C_{j+1}\ge F\frac{C_j^p}{p^{(n-1)j}}\ (j\in\N),
\]
where $F$ is a positive constant defined by 
\[
F=\frac{D_n(n-2)!(p-1)^{n-1}}{(n-1)^{n-1}}.
\]
Due to the induction argument again, we obtain, for $j\ge2$, that
\[
\log C_j\ge p^{j-1}\left\{\log C_1+\sum_{k=0}^{j-2}\frac{p^{k}\log F
-(j-1-k)p^{k}\log(p^{n-1})}{p^{j-1}}\right\}.
\]
As before, this inequality yields
that there exist a constant $S$ independent of $j$ such that 
\[
C_j\ge \exp\{p^{j-1}(\log C_1+S)\quad\mbox{for}\ j\ge 2.
\]
Combining all the estimates above, we reach the the final inequality 
\[
\begin{array}{lll}
\d W(\xi)\ge \exp\{p^{j-1}I(\xi)\}\frac{\xi^{(pq+n-2)/(p-1)}}{(\xi-2k)^{(n-1)/(p-1)}}
\ \mbox{for}\ \xi \ge 2k,
\end{array}
\]
where we set 
\[
I(\xi)=\log\left(E_2e^S\e^p(\xi-2k)^{(n-1)/(p-1)}\xi^{-(pq+n-2)/(p-1)}\right).
\]
\par
Note that
\[
\frac{n-1}{p-1}-\frac{pq+n-2}{p-1}=\frac{1-pq}{p-1}.
\]
If there exist a point $\xi_0 \in \{\xi \ge 4k\}\subset\{\xi \ge 2k\}$ 
such that $I(\xi_0)>0$,
the desired conclusion can be established by the same argument as in the previous section.
$I(\xi_0)>0$ is equivalent to
\[
\xi_0>2^{(n-1)/(1-pq)}
\left(e^{S}E_2\right)^{-(p-1)/(1-pq)}\e^{-2p(p-1)/\gamma(p,n)}
\]
in this case, so that the proof of the subcritical case is now completed.
\hfill$\Box$
\section{Upper bound of the lifespan for the critical case in even dimensions}
In this section, we prove the blow-up theorem
in the critical case in even dimensions. 
The proof is based on the one in odd dimensional case. 
However, Huygens' principle for $u^0$ is no longer available.
Therefore the blow-up domain to ensure the positivity of the linear part is modified. 

\begin{prop}
\label{prop:lifespan_upper_critical_even} 
Suppose that of the assumption of Theorem \ref{thm:main2} are fulfilled. 
Let $u$ be a $C^0$-solution of (\ref{IE}) if $n>4$ and $p=p_0(n)$,
or $u$ be a classical solution of (\ref{NWIVP}) if $n=4$ and $p=p_0(4)$ 
in $\R^n\times[0,T]$. Then, there exists a positive constant $\e_0=\e_0(g,n,p,k)$ 
such that $T$ cannot be taken as 
\begin{equation}
\label{upper_lifespan_critical_even}
T>\exp\left(c\e^{-p(p-1)}\right)\ \mbox{if}\ p=p_0(n)
\end{equation}
for $0<\e \le \e_0$, where $c$ is a positive constant independent of $\e$.
\end{prop} 

\par\noindent
{\bf Proof.} 
Similarly to the odd dimensional case,
we may assume that the solution of (\ref{IE}) is radially symmetric 
without loss of the generality.
Let $u=u(r,t)$ be a $C^0$-solution of
\begin{equation}
\label{rad_IE_even}
u\ge\e u^0+L_{even,1}(|u|^p)\quad\mbox{in}\quad(0,\infty)\times[0,T],
\end{equation}
where $L_{even,1}$ is defined by (\ref{even_int_domain1})
and $u^0=u^0(r,t)$ is a solution of (\ref{u^0}). 
\vskip10pt
\par\noindent
{\bf [The 1st step] Estimate of $\v{u^0}$.}
\par
We shall employ the following representation of $u^0$. 
\begin{lem}[Rammaha \cite{R87}]
\label{lem:rad_u^0_even}
Let $n=4,6,8\cdots$ and $u^0$ be a solution of (\ref{u^0}).
Then, $u^0$ is represented by
\[
\begin{array}{l}
\d u^0(r,t)=\frac{2}{\pi r^{(n-2)/2}}\int_0^t\frac{\rho d\rho}{\sqrt{t^2-\rho^2}}\times\\ 
\d\quad\times\int_{|r-\rho|}^{r+\rho}\frac{\lambda^{(n-2)/2}g(\lambda)
T_{(n-4)/2}\left((\lambda^2+r^2-\rho^2)/(2r\lambda)\right)d\lambda} 
{\sqrt{\lambda^2-(r-\rho)^2}\sqrt{(r+\rho)^2-\lambda^2}},
\end{array}
\]
where $T_k$ is Tschebyscheff polynomials of degree $k$ defined by 
\[
T_k(z)=\frac{(-1)^k}{(2k-1)!!}(1-z^2)^{1/2}\frac{d^k}{dz^k}(1-z^2)^{k-(1/2)}. 
\]
\end{lem}
See (6b) on 681p. in \cite{R87} for the proof.
This lemma implies the following estimate.
\begin{lem}[Rammaha \cite{R87}]
\label{lem:lower_v_even}
Let $n=4,6,8,\cdots$. Assume (\ref{blowup_asm}). 
Then there exists a positive constant $C_{g}$
such that, for $t+k_0<r<t+k_1$ and $t\ge k_2$,
\begin{equation}
\label{lower_v_even}
u^0(r,t)\ge\frac{C_{g}}{r^{(n-1)/2}},
\end{equation}
where $\d k_2=k-k_0$. 
\end{lem}
See Lemma 2 on 682p. in \cite{R87} for the proof. 
\par
Our frame in the iteration argument is obtained by the following lemma. 
\begin{lem}
\label{lm:frame_even}
Let $u$ be a $C^0$-solution of (\ref{rad_IE_even}).
Assume (\ref{blowup_asm}). 
Then $u$ in $\Sigma_0=\{(r,t)\ :\ 2k\le t-r\le r\}$ satisfies
\begin{equation}
\label{frame1_even}
\begin{array}{l}
\d u(r,t)\ge \frac{C2^{(n-1)/2}(t-r)^{(n-1)/2}}{{(n-1)}r^{(3n-5)/2}}\times\\
\d\times\int\!\!\!\int_{R(r,t)}\{(t-r-\tau+\lambda)(t+r-\tau-\lambda)\}^{(n-2)/2}
|u(\lambda,\tau)|^pd\lambda d\tau+\\
\qquad\d+\frac{F_1(t-r)^{(3n-3)/2-(n-1)p/2}}{r^{(3n-5)/2}}\e^p+\e u^0(r,t),
\end{array}
\end{equation}
where $C$ is the one in (\ref{even_int_domain1}) and
\[
F_1=\frac{CC_g^p(k_1-k_0)2^{(11-3n)/2-(n-1)p}}{n(n-1)}.
\]
\end{lem}
\par\noindent
{\bf Proof.} 
In view of (\ref{even_int_domain1}), we have that
\[
\begin{array}{lll}
L_{even,1}(|u|^p)(r,t)
&\d \ge \frac{C}{r^{n-2}}\int_{0}^{t}(t-\tau)^{2-n}d\tau
\int_{|t-r-\tau|}^{t+r-\tau}\lambda |u(\lambda,\tau)|^pd\lambda\times&\\
&\d\quad\times \int_{|\lambda-r|}^{t-\tau}
\frac{\rho h(\lambda,\rho,r)}{\sqrt{(t-\tau)^2-\rho^2}}d\rho
\end{array}
\]
in $\Sigma_0$.
Noticing that $(\lambda+r)^2-\rho^2\ge (\lambda+r)^2-(t-\tau)^2$ 
for $\rho \le t-\tau$, we get
\[
\begin{array}{l}
\d L_{even,1}(|u|^p)(r,t)
\ge \frac{C}{r^{n-2}}\int_{0}^{t}(t-\tau)^{2-n}d\tau\times\\ 
\d \qquad\times\int_{|t-r-\tau|}^{t+r-\tau}\frac{\lambda |u(\lambda,\tau)|^p
\{(\lambda+r)^2-(t-\tau)^2\}^{(n-3)/2}d\lambda}{\{(t-\tau)^2-(\lambda-r)^2\}^{1/2}}\times\\
\d \qquad\times \int_{|\lambda-r|}^{t-\tau} \rho \{\rho^2-(\lambda-r)^2\}^{(n-3)/2} d\rho
\end{array}
\]
in $\Sigma_0$. 
Since the $\rho$-integral above is 
\[ 
\frac{1}{n-1}\{(t-\tau)^2-(\lambda-r)^2\}^{(n-1)/2},
\]
we obtain that
\[
\begin{array}{ll}
L_{even,1}(|u|^p)(r,t)
&\d \ge \frac{C}{r^{n-2}(n-1)}\int_{0}^{t}(t-\tau)^{2-n}d\tau\times\\
&\d \quad\times\int_{|t-r-\tau|}^{t+r-\tau}
\{(\lambda+r)^2-(t-\tau)^2\}^{(n-3)/2}\times\\
&\d \qquad\times\{(t-\tau)^2-(\lambda-r)^2\}^{(n-2)/2}\lambda
|u(\lambda,\tau)|^pd\lambda\\
&\ge J_1+J_2
\end{array}
\]
in $\Sigma_0$, where we set 
\[
\begin{array}{l}
\d J_1(r,t)=\frac{C}{(n-1)r^{n-2}}\int_{R(r,t)}(t-\tau)^{2-n}
\{(\lambda+r)^2-(t-\tau)^2\}^{(n-3)/2}\times\\
\d\qquad \times \{(t-\tau)^2-(\lambda-r)^2\}^{(n-2)/2}
\lambda|u(\lambda,\tau)|^pd\lambda d\tau,
\end{array}
\]
\[
\begin{array}{l}
\d J_2(r,t)=\frac{C}{(n-1)r^{n-2}}\int_{S(r,t)}(t-\tau)^{2-n}
\{(\lambda+r)^2-(t-\tau)^2\}^{(n-3)/2}\times \\
\d\qquad \times \{(t-\tau)^2-(\lambda-r)^2\}^{(n-2)/2}
\lambda|u(\lambda,\tau)|^pd\lambda d\tau.
\end{array}
\]
Changing variables by (\ref{alpha_beta}) in $J_1$, we have that
\[
\begin{array}{c}
\d J_1(r,t)\ge\frac{C}{2(n-1)r^{n-2}}
\int_{2k}^{t-r}(t-r-\beta)^{(n-2)/2}(t+r-\beta)^{(n-3)/2}d\beta\times\\
\d\times\int_{2(t-r)+\beta}^{t+r}
\{\alpha-(t-r)\}^{(n-3)/2}(t+r-\alpha)^{(n-2)/2}\times\\
\d\times\{t-(\alpha+\beta)/2\}^{2-n}(\alpha-\beta)|u(\lambda,\tau)|^p d\alpha
\end{array}
\]
in $\Sigma_0$.
Note that
\[
\begin{array}{l}
\d t+r-\beta \ge 2r,\ t-\frac{\alpha+\beta}{2}\le r,\\  
\d \alpha-\beta \ge 2(t-r),\ \alpha-(t-r) \ge t-r+\beta \ge t-r 
\end{array}
\]
hold in the domain of the integral above.
Hence we get
\[
\begin{array}{c}
J_1(r,t)
\ge\d \frac{C2^{(n-3)/2}(t-r)^{(n-1)/2}}{{(n-1)}r^{(3n-5)/2}}
\int_{2k}^{t-r}(t-r-\beta)^{(n-2)/2}d\beta\times\\
\d\times\int_{2(t-r)+\beta}^{t+r}(t+r-\alpha)^{(n-2)/2}|u(\lambda,\tau)|^p d\alpha
\end{array}
\]
in $\Sigma_0$. 
Therefore, we obtain the first term of the right-hand side in (\ref{frame1_even}). 
\par
Similarly to the above, $J_2(r,t)$ is bounded from below by
\[
\begin{array}{c}
\d \frac{C}{2(n-1)r^{n-2}}\int_{-k_1}^{-k_0}
(t-r-\beta)^{(n-2)/2}(t+r-\beta)^{(n-3)/2}d\beta\\
\d \times\int_{2(t-r)+\beta}^{t+r}\{\alpha-(t-r)\}^{(n-3)/2}
(t+r-\alpha)^{(n-2)/2}\\
\d\times\{t-(\alpha+\beta)/2\}^{2-n}(\alpha-\beta)|u(\lambda,\tau)|^p d\alpha
\end{array}
\]
in $\Sigma_0$.
Note that 
\[
\begin{array}{l}
\d t+r-\beta \ge r,\ t-\frac{\alpha+\beta}{2}\le 2r,\ \alpha-(t-r)\ge t-r-k\\ 
\d \mbox{and}\ t-r-\beta\ge t-r 
\end{array}
\]
hold in the domain of the integral above. 
Hence (\ref{lower_v_even}) yields that
$J_2(r,t)$ in $\Sigma_0$ is estimated from below by
\[
\begin{array}{l}
\ge\d \frac{\e^pCC_g^p2^{1-n}(t-r)^{(n-2)/2}
(t-r-k)^{(n-3)/2}}{(n-1)r^{(3n-5)/2}}\int_{-k_1}^{-k_0}d\beta\\
\d \quad\times\int_{2(t-r)+\beta}^{t+r}(\alpha-\beta)^{1-(n-1)p/2}
(t+r-\alpha)^{(n-2)/2}d\alpha\\
\ge\d \frac{\e^pCC_g^p2^{3-n-(n-1)p}
(t-r)^{n/2-(n-1)p/2}(t-r-k)^{(n-3)/2}}{(n-1)r^{(3n-5)/2}}\int_{-k_1}^{-k_0}d\beta\times\\
\d \quad\times\int_{2(t-r)+\beta}^{3(t-r)}\{3(t-r)-\alpha\}^{(n-2)/2}d\alpha.
\end{array}
\] 
The second term of the right-hand side of (\ref{frame1_even}) 
follows from this inequality. 
Therefore, the proof of Lemma \ref{lm:frame_even} is ended.
\hfill$\Box$
\par
Next, we shall show the positivity of the right-hand side of (\ref{frame1_even}). 
Under the condition (\ref{blowup_asm}), (\ref{decay_est_v_0even}) yields that
\[
\begin{array}{lll}
\d \e u^0(r,t)&\d \ge \frac{-C_{n,k,0,g}\e}{(t+r+2k)^{(n-1)/2}(t-r+2k)^{(n-1)/2}}&\\
&\d \ge \frac{-C_{n,k,0,g}\e}{r^{(n-1)/2}(t-r)^{(n-1)/2}}
\end{array}
\]
for $t-r\ge -k$.
Let us define a domain
\[
\Sigma_1
=\left\{(r,t)\in (0,\infty)^2:\ r \ge t-r\ge\frac{r}{2},\ r\ge K\e^{-L}\right\},
\]
where we set 
\[
\begin{array}{llll}
\d K=\left(2^{2n-(n-1)p/2-1}F_1^{-1}C_{n,k,0,g} \right)^{1/\left(n-(n-1)p/2\right)},\\
\d L=\frac{p-1}{n-(n-1)p/2}>0.
\end{array}
\]
Taking $\e$ to satisfy 
\[
K\e^{-L}\ge 4k
\]
and setting
\[
A(r,t)=\frac{F_1(t-r)^{(3n-3)/2-(n-1)p/2}}{r^{(3n-5)/2}}>0,
\]
we obtain that, in $\Sigma_1$,
\[
\begin{array}{lll}
\d \frac{A(r,t)}{2}\e^p+\e u^0(r,t)\\
\d\ge\frac{F_1\e^p(t-r)^{2n-2-(n-1)p/2}r^{(n-1)/2}
-2C_{n,k,0,g}\e r^{(3n-5)/2}}{2r^{2n-3}(t-r)^{(n-1)/2}}\ge0.
\end{array}
\]
Making use of this inequality, we obtain that
\[
\begin{array}{l}
u(r,t)\ge\d \frac{C2^{(n-3)/2}(t-r)^{(n-1)/2}}{{(n-1)}r^{(3n-5)/2}}
\int_{2k}^{t-r}(t-r-\beta)^{(n-2)/2}d\beta\times\\
\qquad\d \times\int_{2(t-r)+\beta}^{t+r}(t+r-\alpha)^{(n-2)/2}
|u(\lambda,\tau)|^p d\alpha+\frac{A(r,t)}{2}\e^p
\end{array}
\]
in $\Sigma_1$.
Cutting the domain of the integral, we get 
\[
\begin{array}{l}
u(r,t)>\d \frac{C2^{(n-3)/2}(t-r)^{(n-1)/2}}{{(n-1)}r^{(3n-5)/2}}
\int_{K\e^{-L}/2}^{t-r}(t-r-\beta)^{(n-2)/2}d\beta\times\\
\qquad\d \times\int_{3(t-r)}^{t+r}(t+r-\alpha)^{(n-2)/2}
|u(\lambda,\tau)|^p d\alpha+\frac{A(r,t)}{4}\e^p
\end{array}
\]
in $\Sigma_1$.
Here we introduce a change of variables $(\alpha,\beta)$ to $(\xi,\eta)$ by
\[
\xi=\alpha,
\ \eta=\frac{\alpha+\beta}{2}-\frac{3}{2}\cdot\frac{\alpha-\beta}{2}
=\frac{5\beta-\alpha}{4}.
\]
Then, cutting the domain of the integral again, we get 
\[
\begin{array}{l}
u(r,t)>
\d \frac{C2^{(n+1)/2}(t-r)^{(n-1)/2}}{{5(n-1)}r^{(3n-5)/2}}
\int_{K\e^{-L}/2}^{t-3r/2}
\left\{t-r-\left(\frac{4\eta+\xi}{5}\right)\right\}^{(n-2)/2}d\eta\times\\
\qquad\d \times\int_{3(t-r)}^{t+r}(t+r-\xi)^{(n-2)/2}|u(\lambda,\tau)|^p d\xi+B(r,t)\e^p
\end{array}
\]
in $\Sigma_2$, where we set
\[
\Sigma_2
=\left\{(r,t)\in (0,\infty)^2:\ \frac{r}{2}
\ge t-\frac{3}{2}r \ge \frac{K\e^{-L}}{2}\right\}
\]
and
\[
B(r,t)=\frac{A(r,t)}{4}.
\]
Therefore we obtain that, in $\Sigma_2$,
\begin{equation}
\label{u_even_equal}
\begin{array}{l}
\d u(r,t)
>\d \frac{C2^{(3n-3)/2}(t-r)^{(n-1)/2}}{5^{n/2}(n-1)r^{(3n-5)/2}}
\int_{K\e^{-L}/2}^{t-3r/2}\left(t-\frac{3}{2}r-\eta \right)^{(n-2)/2}d\eta\times\\
\d \qquad\times\int_{3(t-r)}^{t+r}(t+r-\alpha)^{(n-2)/2}
|u(\lambda,\tau)|^p d\alpha+B(r,t)\e^p.
\end{array}
\end{equation}
\vskip10pt
\par\noindent
{\bf [The 2nd Step] Comparison argument.}
\par
Let us consider a solution $y$ of
\begin{equation}
\label{y_equal_even}
\begin{array}{c}
\d y\left(t-\frac{3}{2}r\right)
=\frac{C2^{(3n-3)/2}(t-r)^{(n-1)/2}}{5^{n/2}(n-1)r^{(3n-5)/2}}\int_{K\e^{-L}/2}^{t-3r/2}
\left(t-\frac{3}{2}r-\eta\right)^{(n-2)/2}\hspace{-20pt}d\eta\times\\
\d\times\int_{3(t-r)}^{t+r}(t+r-\alpha)^{(n-2)/2}|y(\eta)|^p d\alpha+B(r,t)\e^p.
\end{array}
\end{equation}
Then we have the following comparison lemma. 
\begin{lem}
\label{lem:comparison_arg_even}
Let $u$ be a solution of (\ref{rad_IE_even}) and $y$ be a solution of (\ref{y_equal_even}). 
Then, $u$ and $y$ satisfy 
\[
u>y\quad \mbox{in}\ \Sigma_2.
\]
\end{lem}
{\bf Proof.} 
Fix a point for any $(r_0,t_0) \in \Sigma_2$. Define
\[
\Lambda(r,t)=\left\{(\lambda,\tau)\in D(r,t):\ \frac{K\e^{-L}}{2}\le 
\tau-\frac{3}{2}\lambda\le \frac{\lambda}{2}\right\},
\]
where
\[
D(r,t)=\left\{(\lambda,\tau):\ t-r\le \tau+\lambda\le t+r,
-k\le \tau-\lambda\le t-r\right\}.
\] 
Let us consider $u$ and $y$ in $\Lambda(r_0,t_0)$. 
Note that $u>y$ on $\d \tau-\frac{3}{2}\lambda=\frac{K\e^{-L}}{2}$
and at $(K\e^{-L},2K\e^{-L})$ 
which is an edge point of $\Sigma_2$. 
By compactness of the closure of $\Lambda(r_0,t_0)$, 
we have $u>y$ in a neighborhood of
$\d \tau-\frac{3}{2}\lambda=\frac{K\e^{-L}}{2}$
and $\lambda\ge K\e^{-L}$.
\par
Assume that there exist a point $(r_1,t_1)$ with $u(r_1,t_1)=y(t_1-3r_1/2)$ 
which is nearest to $(K\e^{-L},2K\e^{-L})$ in such a neighborhood.
Since $u>y$ in $R'(r_1,t_1)$, we have 
\[
\begin{array}{l}
\d \frac{C2^{(3n-5)/2}(t_1-r_1)^{(n-1)/2}}{5^{(n-2)/2}(n-1)r_1^{(3n-5)/2}}
\int\!\!\!\int_{R'(r_1,t_1)}
\left(t_1-\frac{3}{2}r_1-\tau+\frac{3}{2}\lambda\right)^{(n-2)/2}\times\\
\d \quad\times\left(t_1+r_1-\tau-\lambda\right)^{(n-2)/2}
|u(\lambda,\tau)|^pd\lambda d\tau+B(r_1,t_1)\e^p\\
>\d \frac{C2^{(3n-5)/2}(t_1-r_1)^{(n-1)/2}}{5^{(n-2)/2}(n-1)r_1^{(3n-5)/2}}
\int\!\!\!\int_{R'(r_1,t_1)}
\left(t_1-\frac{3}{2}r_1-\tau+\frac{3}{2}\lambda\right)^{(n-2)/2}\times\\
\quad \d\times\biggl(t_1+r_1-\tau-\lambda \biggr)^{(n-2)/2}
\left|y\left(\tau-\frac{3}{2}\lambda\right)\right|^pd\lambda d\tau+B(r_1,t_1)\e^p,
\end{array}
\]
where we set
\[
R'(r,t)=\left\{(\lambda,\tau):\ 3(t-r)\le \tau+\lambda\le t+r, 
\frac{K\e^{-L}}{2}\le \tau-\frac{3}{2}\lambda\le t-\frac{3}{2}r\right\}.
\]
In view of (\ref{u_even_equal}) and (\ref{y_equal_even}),
this inequality yield that $u>y$ at $(r_1,t_1)$,
which is a contradiction to the definition of $(r_1,t_1)$. 
Therefore, the proof of Lemma \ref{lem:comparison_arg_even} is now established
by the same argument as the one for Lemma \ref{lem:comparison_arg}.  
\hfill$\Box$
\par
We note that Lemma \ref{lem:comparison_arg_even} implies
that the lifespan of $y$ is greater than the lifespan of $u$,
so that it is sufficient to look for the lifespan of $y$ in $\Sigma_2$.
By definition $y$ in (\ref{y_equal_even}), we have
\[
\begin{array}{ll}
y(\xi)=
&\d \frac{C3^{(n-1)/2}\xi^{2-n}}{2^{(3n-7)/2}5^{n/2}(n-1)}
\int_{K\e^{-L}/2}^{\xi}(\xi-\eta)^{(n-2)/2}|y(\eta)|^pd\eta\\
&\d \times\int_{9\xi}^{11\xi}(11\xi-\alpha)^{(n-2)/2}d\alpha
+\frac{\e^pF_1\xi^{-q-(n-1)/2}}{2^{(3n-3)}3^{(n-1)p/2-(3n-3)/2}}
\end{array}
\]
in $\Gamma_1$, where we set
\[
\xi=\frac{r}{4},\ \Gamma_1=\left\{t-\frac{3}{2}r=\xi,
r\ge 2K\e^{-L}\right\}\subset\Sigma_2.
\]
Hence, we obtain that
\[
\begin{array}{ll}
\d y(\xi)\ge
&\d \frac{C3^{(n-1)/2}\xi^{(4-n)/2}}{2^{(2n-9)/2}5^{n/2}n(n-1)}
\int_{K\e^{-L}/2}^{\xi}(\xi-\eta)^{(n-2)/2}|y(\eta)|^pd\eta\\
&\d+\frac{\e^pF_1\xi^{-q-(n-1)/2}}{2^{(3n-3)}3^{(n-1)p/2-(3n-3)/2}}
\end{array}
\]
for $\xi\ge K\e^{-L}/2$.
Then it follows from the setting
\[
Y(\xi)=\xi^{q+(n-1)/2}y(\xi)
\] 
that
\begin{equation}
\label{frame_sub_even}
Y(\xi)\ge E_n\xi^{q+3/2}\int_{K\e^{-L}/2}^{\xi}
\frac{(\xi-\eta)^{(n-2)/2}|Y(\eta)|^p d\eta}{\eta^{(n-1)p/2+pq}}+F_2\e^p,
\end{equation}
where we set
\[ 
E_n=\frac{C3^{(n-1)/2}}{2^{(2n-9)/2}5^{n/2}n(n-1)},
\ F_2=\frac{F_1}{2^{(3n-3)}3^{(n-1)p/2-(3n-3)/2}}.
\]
Therefore we obtain the iteration frame in this section,
\begin{equation}
\label{frame3}
Y(\xi)\ge E_n\int_{K\e^{-L}/2}^{\xi}
\left(\frac{\xi-\eta}{\xi}\right)^{(n-2)/2}\frac{|Y(\eta)|^p}{\eta^{pq}}d\eta+F_2\e^p
\quad\mbox{for}\ \xi\ge\frac{K\e^{-L}}{2}.
\end{equation}
\vskip10pt
\par\noindent
{\bf [The 3rd step] Slicing method with the iteration.}
\par
Let us define a blow-up domain as follows. Let us set
\[
\Gamma_j=\{\xi\ge l_{j}K\e^{-L}\},\quad l_j=\frac{1}{2}+\cdots+\frac{1}{2^j}\ (j\in\N).
\]
We shall use the fact that a sequence $\{l_j\}$ is monotonously increasing
and bounded as $\d \frac{1}{2}<l_j<1$, so that $\Gamma_{j+1}\subset \Gamma_{j}$.
Assume an estimate of the form 
\begin{equation}
\label{j_times_even}
Y(\xi)\ge C_j \left(\log\frac{\xi}{l_jK\e^{-L}}\right)^{a_j}\ 
\mbox{in}\ \Gamma_j,
\end{equation}
where $a_j\ge 0$ and $C_j>0$. Putting (\ref{j_times_even}) 
into (\ref{frame3}) and recalling $pq=1$, 
we get an estimate in $\Gamma_{j+1}$ such as 
\[
Y(\xi)
\ge\d E_nC_j^p\int_{l_jK\e^{-L}}^{\xi}
\left(\frac{\xi-\eta}{\xi}\right)^{(n-2)/2}
\left(\log\frac{\eta}{l_jK\e^{-L}}\right)^{pa_j}\frac{d\eta}{\eta}.
\]
Noting that $\d \frac{l_{j}}{l_{j+1}}\xi\ge l_jK\e^{-L}$ in $\Gamma_{j+1}$,
we have 
\[
\begin{array}{lll}
\d Y(\xi)
&\ge&\d E_nC_j^p\left(1-\frac{l_j}{l_{j+1}}\right)^{(n-2)/2}
\int_{l_jK\e^{-L}}^{l_{j}\xi/l_{j+1}}
\d \left(\log\frac{\eta}{l_jk}\right)^{pa_j}\frac{d\eta}{\eta}\\
&=&\d \frac{E_nC_j^p}{pa_j+1}\left(1-\frac{l_j}{l_{j+1}}\right)^{(n-2)/2}
\left(\log\frac{\xi}{l_{j+1}K\e^{-L}}\right)^{pa_j+1}.
\end{array}
\]
By monotonicity of $\{l_j\}$ and
\[
1-\frac{l_j}{l_{j+1}}=\frac{l_{j+1}-l_{j}}{l_{j+1}}=
\frac{1}{2^{j+1}l_{j+1}} 
\ge\frac{1}{2^{j+1}},
\]
we finally obtain
\begin{equation}
\label{j+1_times_even}
Y(\xi)\ge C_{j+1}\left(\log\frac{\xi}{l_{j+1}K\e^{-L}}\right)^{pa_j+1}
\ \mbox{in}\ \Gamma_{j+1},
\end{equation}
where we set
\[
C_{j+1}=\frac{E_nC_j^p}{2^{(n-2)(j+1)/2}(pa_j+1)}.
\]
\par
Now, we are in a position to define sequences in the iteration. In view 
of (\ref{frame3}), the first estimate is $Y(\xi)\ge F_2 \e^p$,
so that, with the help of (\ref{j_times_even}) and (\ref{j+1_times_even}), a 
sequence $\{a_j\}$ should be defined by 
\[
\d a_1=0,\ a_{j+1}=pa_j+1\ (j\in\N).
\]
Also a sequence $\{C_j\}$ should be defined by 
\[
C_1=F_2\e^p,\ C_{j+1}=\frac{E_nC_j^p}{2^{(n-2)(j+1)/2}(pa_j+1)}\ (j\in\N).
\]
One can easily check that 
\[
a_j=\frac{p^{j-1}-1}{p-1}\ (j\in\N)
\]
which gives us
\[
\frac{1}{pa_j+1}\ge \frac{p-1}{p^{j}}.
\]
Thus one can find that 
\[
C_{j+1}\ge E\frac{C_j^p}{(2^{(n-2)/2}p)^j}\ (j\in\N),
\]
where $E$ is a positive constant defined by 
\[
E=\frac{E_n(p-1)}{2^{(n-2)/2}}.
\]
Hence, we inductively obtain, for $j\ge2$, that 
\[
\log C_j\ge p^{j-1}\left\{\log C_1+\sum_{k=0}^{j-2}\frac{p^{k}\log E
-(j-1-k)p^{k}\log(2^{(n-2)/2}p)}{p^{j-1}}\right\}.
\]
This inequality yields that there exist a constant $S$ independent of $j$ such that 
\[
C_j\ge \exp\{p^{j-1}(\log C_1+S)\}\quad\mbox{for}\ j\ge 2.
\]
Combining all the estimates above and making use of the monotonicity of $\Gamma_j$, 
we obtain the final inequality,
\[
\begin{array}{lll}
\d Y(\xi)
&\ge&\d \exp\{p^{j-1}(\log C_1+S)\}
\left(\log\frac{\xi}{K\e^{-L}}\right)^{(p^{j-1}-1)/(p-1)}\\
&=&\d \exp\{p^{j-1}I(\xi)\}\left(\log\frac{\xi}{K\e^{-L}}\right)^{-1/(p-1)}.
\end{array}
\]
in $\d \Gamma_{\infty}=\{\xi \ge K\e^{-L}\}$, where we set 
\[
I(\xi)=\log\left(e^SF_2\e^p\left(\log\frac{\xi}{K\e^{-L}}\right)^{1/(p-1)}\right).
\]
\par
If there exist a point $\xi_0 \in \Gamma_{\infty}\subset \Gamma_{j}\ (j\ge1)$ 
such that $I(\xi_0) >0$, we get the desired conclusion
by the same argument in the end of section 7.
In this case, $I(\xi_0)>0$ is equivalent to
\[
\xi_0>\exp\{(e^SF_2)^{-(p-1)}\e^{-p(p-1)}\}K\e^{-L}. 
\]
Therefore the proof of the critical case is now completed.
\hfill$\Box$
\section{Upper bound of the lifespan for the subcritical case in even dimensions}
Similarly to the previous section,
we prove the blow-up result of solution for (\ref{IE}) in 
the subcritical case in even space dimensions. 
Note that we do not have to make use of the slicing method.
\begin{prop}
\label{prop:lifespan_upper_sub_even}
Suppose that the same assumption of Theorem \ref{thm:main2} 
are fulfilled. Let $u$ be a $C^0$-solution of (\ref{IE}) $\R^n\times[0,T]$. 
Then there exists a positive constant $\e_0=\e_0(g,n,p,k)$
such that $T$ cannot be taken as 
\begin{equation}
\label{upper_lifespan_sub_even}
T>\d c\e^{-2p(p-1)/\gamma(p,n)}\ \mbox{if}\ 1<p<p_0(n)
\end{equation}
for $0<\e \le \e_0$, where $c$ is a positive constant independent of $\e$.
\end{prop}
{\bf Proof.} 
Because of the fact that $-(n-1)p/2-pq<0$ for $n\ge4$,
(\ref{frame_sub_even}) yields 
\begin{equation}
\label{frame_sub_even2}
Y(\xi)\ge E_n\xi^{-(n-2)/2-pq}\int_{K\e^{-L}/2}^{\xi}(\xi-\eta)^{(n-2)/2}
|Y(\eta)|^pd\eta+F_2\e^p
\end{equation}
for $\xi\ge K\e^{-L}/2$.
This is our iteration frame in this case.
\par
Assume an estimate of the form 
\begin{equation}
\label{j_times_sub_even}
Y(\xi)\ge C_j \frac{\left(\d \xi-\frac{K\e^{-L}}{2}\right)^{a_j}}{\xi^{b_j}}\ 
\mbox{in}\ \Gamma_1,
\end{equation} 
where $a_j,b_j\ge 0$ and $C_j>0$.
Then, putting (\ref{j_times_sub_even}) into (\ref{frame_sub_even2}),
we get the following estimate in $\Gamma_1$ of the form 
\[
Y(\xi)\ge\d E_nC_j^p\xi^{-(n-2)/2-pq-pb_j}\int_{K\e^{-L}/2}^{\xi}
(\xi-\eta)^{(n-2)/2}\left(\d \eta-\frac{K\e^{-L}}{2}\right)^{pa_j}d\eta.
\]
Applying the integration by parts to $\beta$-integral $(n-2)/2$ times,
we obtain that 
\[
\begin{array}{lll}
&\d \frac{(n-2)}{2(pa_j+1)}\int_{K\e^{-L}/2}^{\xi}(\xi-\eta)^{(n-4)/2}
\left(\d \eta-\frac{K\e^{-L}}{2}\right)^{pa_j+1}d\eta\\
&\ge\d \frac{(n-2)(n-4)}{(pa_j+2)^22\cdot 2}
\int_{K\e^{-L}/2}^{\xi}(\xi-\eta)^{(n-6)/2}
\left(\d \eta-\frac{K\e^{-L}}{2}\right)^{pa_j+2}d\eta\\
&\cdots\\
&\ge\d \frac{((n-2)/2)!}{\d \left(pa_j+\frac{n}{2}\right)^{n/2}}
\left(\d \xi-\frac{K\e^{-L}}{2}\right)^{pa_j+n/2}.
\end{array}
\]
Therefore we finally get 
\begin{equation}
\label{j+1_times_sub_even}
Y(\xi)\ge\frac{C_{j+1}
\left(\d \xi-\frac{K\e^{-L}}{2}\right)^{pa_j+n/2}}
{\xi^{pb_j+(n-2)/2+pq}}\ \mbox{in}\ \Gamma_1,
\end{equation}
where
\[
C_{j+1}=\frac{((n-2)/2)!E_nC_j^p}
{\d \biggl(pa_j+\frac{n}{2}\biggr)^{n/2}}.
\]
\par
Now, we are in a position to define sequences in the iteration. 
In view of (\ref{frame_sub_even2}), the first estimate is $Y(\xi)\ge F_2 \e^p$,
so that, with the help of (\ref{j_times_sub_even}) and 
(\ref{j+1_times_sub_even}),
sequences $\{a_j\}$ and $\{b_j\}$ should be defined by 
\[
a_1=0,\ a_{j+1}=pa_j+\frac{n}{2}\ (j\in\N)
\]
and
\[
b_1=0,\ b_{j+1}=pb_j+pq+\frac{n-2}{2}\ (j\in\N).
\]
Also a sequence $\{C_j\}$ should be defined by 
\[
C_1=F_2\e^p,\ C_{j+1}=\frac{((n-2)/2)!E_nC_j^p}
{\d\biggl(pa_j+\frac{n}{2}\biggr)^{n/2}}\ (j\in\N).
\]
One can readily check that 
\[
a_j=\frac{n}{2(p-1)}(p^{j-1}-1),
\ b_j=\frac{pq+(n-2)/2}{p-1}(p^{j-1}-1)\ (j\in\N),
\]
which gives us
\[
\frac{1}{pa_j+n/2}\ge \frac{2(p-1)}{p^{j}n}.
\]
Hence one can find that 
\[
C_{j+1}\ge F\frac{C_j^p}{p^{(n/2)j}}\ (j\in\N),
\]
where $F$ is a positive constant defined by 
\[
F=\frac{E_n((n-2)/2)!\{2(p-1)\}^{n/2}}{n^{n/2}}.
\]
Due to the induction argument again, we obtain, for $j\ge2$, that
\[
\log C_j\ge p^{j-1}\left\{\log C_1+\sum_{k=0}^{j-2}\frac{p^{k}\log F
-(j-1-k)p^{k}\log(p^{n/2})}{p^{j-1}}\right\}.
\]
This inequality yields that there exist a constant $S$ independent of $j$ such that 
\[
C_j\ge \exp\{p^{j-1}(\log C_1+S)\}\quad\mbox{for}\ j\ge 2.
\]
Combining all the estimates, we reach the final inequality 
\[
Y(\xi)\ge\exp\{p^{j-1}I(\xi)\}\frac{\xi^{(pq+(n-2)/2)/(p-1)}}
{\left(\d \xi-\frac{K\e^{-L}}{2}\right)^{n/2(p-1)}}
\quad\mbox{for}\ \xi \ge \frac{K\e^{-L}}{2},
\]
where we set 
\[
I(\xi)=
\log\left(F_2e^S\e^p\left(\d \xi-\frac{K\e^{-L}}{2}\right)^{n/2(p-1)}\xi^{-(pq+(n-2)/2)/(p-1)}\right).
\]
\par
Note that
\[
\frac{n}{2(p-1)}-\frac{pq+(n-2)/2}{p-1}=\frac{1-pq}{p-1}.
\]
If there exist a point $\xi_0 \in \{\xi \ge K\e^{-L}\} \subset 
\{\xi \ge (K\e^{-L})/2\}$ such that $\d I(\xi_0) >0$,
we get the desired conclusion as before.
$I(\xi_0)>0$ is equivalent to
\[
\xi_0>2^{n/2(1-pq)}\left(e^{S}F_2\right)^{-(1-p)/(1-pq)}
\e^{-2p(p-1)/\gamma(p,n)}.
\]
Therefore the proof of the subcritical case is completed.
\hfill$\Box$


\bibliographystyle{plain}

\end{document}